\documentclass{amsproc}
\usepackage{amsmath,amssymb,amsthm,mathtools}
\usepackage{booktabs,array,needspace}
\usepackage{microtype}
\usepackage{tikz}
\usetikzlibrary{arrows.meta}
\usepackage[hidelinks]{hyperref}

\newtheorem{theorem}{Theorem}[section]
\newtheorem{introtheorem}{Theorem}

\newtheorem{lemma}[theorem]{Lemma}
\newtheorem{corollary}[theorem]{Corollary}
\newtheorem*{corollary*}{Corollary}
\newtheorem{proposition}[theorem]{Proposition}
\newtheorem{conjecture}[theorem]{Conjecture}
\newtheorem*{firstcoefflemma}{Lemma}
\theoremstyle{remark}
\newtheorem{remark}[theorem]{Remark}
\newtheorem*{remark*}{Remark}
\newtheorem{example}[theorem]{Example}

\newcommand{\msn}{\mathfrak S_n}
\DeclareMathOperator{\des}{des}
\DeclareMathOperator{\excc}{exc}
\DeclareMathOperator{\lc}{lc}
\DeclareMathOperator{\sech}{sech}
\DeclareMathOperator{\sgn}{sgn}
\newcommand{\Lop}{\mathcal L}

\numberwithin{equation}{section}

\begin{document}

\title[Linear combinations of Eulerian-type polynomials]
{Sharp bounds and interlacing for linear combinations of Eulerian-type polynomials}
\author{Shi-Mei Ma}
\email{shimeimapapers@163.com}
\address{School of Mathematics and Statistics, Shandong University of Technology, Zibo, Shandong 255000, P.R. China}



\date{}

\dedicatory{}

\begin{abstract}
The Eulerian polynomials $A_n(x)$ count permutations of $\{1,2,\ldots,n\}$ by descents. For $n\ge3$, write
$$
A_n+bA_{n-1}+cA_{n-2}
=(A_n-cxA_{n-2})+\bigl(bA_{n-1}+c(1+x)A_{n-2}\bigr).
$$
The summands are symmetric polynomials in degrees $n-1$ and $n-2$. To study the zeros of
$A_n-cxA_{n-2}$, we use the operators
$\mathcal L_s=(z^2-1)\frac{d}{dz}+sz$.

Let $F$ and $G$ have the same degree $d\ge2$ and be both even or both
odd. Suppose $F(1)=0$ and their zeros are simple and lie in $[-1,1]$
and $(-1,1)$, respectively. If their positive zeros alternate strictly,
so do those of $\mathcal L_\mu F$ and $\mathcal L_\nu G$ whenever
$\nu>0$ and $0\le\nu-\mu\le2$. This range is sharp uniformly in
$\nu>0$.
We determine when $A_n-cxA_{n-2}$ has only simple real zeros, proving
Zhang's conjecture, and when $(1+x)A_n+kxA_{n-1}$ is Hurwitz stable. We also prove the
real-rootedness of the descent polynomial of $(n-3)$-stack-sortable
permutations for $n\ge4$. A residue criterion gives sharp
real-rootedness bounds for three-term combinations in five Eulerian
and Narayana families. For cycle-weighted Eulerian polynomials
$A_n(x;q)$, we determine when
$A_n(x;q)+b_iA_{n-1}(x;q)-cxA_{n-2}(x;q)$ forms a fully interlacing row
as $b_i$ increases. For $q>0$, the $q$-Eulerian transformation sends every
nonzero nonnegative combination of $t^{n-j}(1+t)^j$, $0\le j\le n$,
to a real-rooted polynomial for all $n\ge1$ if and only if $q\ge1$.
\end{abstract}
\maketitle

\tableofcontents
\section{Introduction}\label{sec:intro}
The symmetric group $\mathfrak{S}_n$ consists of all permutations of
$[n]=\{1,2,\ldots,n\}$.
Let $\pi=\pi(1)\pi(2)\cdots\pi(n)\in\msn$.
A \emph{descent} of \(\pi\) is an index \(i\in[n-1]\) such that
\(\pi(i)>\pi(i+1)\), and an \emph{excedance} is an index \(i\in[n]\)
such that \(\pi(i)>i\). Let \(\des(\pi)\) and \(\excc(\pi)\) denote
the numbers of descents and excedances of \(\pi\), respectively.
The classical \emph{Eulerian polynomials} are defined by
\[
A_n(x)=\sum_{\pi\in\msn}x^{\des(\pi)}=\sum_{\pi\in\msn}x^{\excc(\pi)},~A_0(x)=1.
\]
They satisfy the recurrence
\begin{equation}\label{eq:Arec}
 A_{n+1}(x)=(1+nx)A_n(x)+x(1-x)\frac{d}{dx}A_n(x).
\end{equation}
The real-rootedness of the Eulerian polynomials goes back to
Frobenius~\cite{Frobenius1910}. The zeros of consecutive Eulerian
polynomials strictly interlace (see~\cite{BrandenJochemko2022,Petersen2015}).
Hence \(A_n(x)+bA_{n-1}(x)\) has only real zeros for \(n\ge2\)
and \(b\in\mathbb R\). A combination of three consecutive Eulerian
polynomials can have nonreal zeros, as \(A_3(x)-4A_2(x)+4A_1(x)=x^2+1\) shows.
We determine when linear combinations of Eulerian-type polynomials have
only real zeros or have all their zeros in the open left half-plane.

Three-term combinations of Eulerian polynomials arise often in the
enumeration of restricted permutations. For example, if \(1\le i<j\le n\), then
\begin{equation}\label{inexAn}
  \sum_{\substack{\pi\in\mathfrak{S}_n\\ \pi(i)\ne i,\ \pi(j)\ne j}}
    x^{\excc(\pi)}
  = A_n(x)-2A_{n-1}(x)+A_{n-2}(x).
\end{equation}
To count permutations fixing \(i\), delete the entry \(i\) at position \(i\).
Relabel the remaining positions and values by the same increasing map
from \([n]\setminus\{i\}\) to \([n-1]\). This bijection with
\(\mathfrak S_{n-1}\) preserves excedances, since it preserves the comparison
of each entry with its position. The permutations fixing \(i\) thus
contribute \(A_{n-1}(x)\), as do those fixing \(j\). Deleting both fixed
points gives the contribution \(A_{n-2}(x)\) from permutations fixing
both \(i\) and \(j\). Inclusion--exclusion gives~\eqref{inexAn}.

Athanasiadis~\cite[Example~3.3]{Athanasiadis2025} used the Eulerian
transformation of Br\"and\'en and Jochemko~\cite{BrandenJochemko2022}
to show that \(A_n(x)+bA_{n-1}(x)+cA_{n-2}(x)\) has only real
zeros for \(n\ge3\) when \(0\le c\le1\) and \(2c\le b\le1+c\).
In~\eqref{inexAn}, the coefficients are \(b=-2\) and \(c=1\), so the
condition \(2c\le b\) fails.

We study the three-term sum by splitting it into two symmetric polynomials.
A polynomial \(p\) of degree at most \(d\) is \emph{symmetric in degree}
\(d\) if \(x^dp(1/x)=p(x)\).
Reversing the entries of a permutation of \([j]\) with \(k\) descents gives
one with \(j-1-k\) descents. Thus \(A_j\) is symmetric in degree \(j-1\)
for \(j\ge1\). For \(n\ge3\), \(xA_{n-2}\) is symmetric in degree \(n-1\), and \((1+x)A_{n-2}\)
is symmetric in degree \(n-2\). Writing \(A_{n-2}=(1+x)A_{n-2}-xA_{n-2}\)
therefore gives
\begin{equation}\label{sd}
 A_n+bA_{n-1}+cA_{n-2}
 =\bigl(A_n-cxA_{n-2}\bigr)+\bigl(bA_{n-1}+c(1+x)A_{n-2}\bigr).
\end{equation}
The two summands on the right are symmetric in degrees \(n-1\) and \(n-2\),
respectively. Reversing the sum in~\eqref{sd} in degree \(n-1\) leaves the
first summand unchanged and multiplies the second by \(x\), giving the
symmetric decomposition of the reversed sum in the sense of
Br\"and\'en and Solus~\cite[Lemma~2.1]{BrandenSolus2021}.
For each fixed \(n\ge3\), the combinations in~\eqref{sd} raise two questions.
\begin{enumerate}
  \item [$(i)$] For which values of \(c\) does \(A_n-cxA_{n-2}\) have only real zeros?
  \item [$(ii)$] For which values of \(c\) does \(A_n+bA_{n-1}+cA_{n-2}\)
  have only real zeros for every real \(b\)?
\end{enumerate}

The polynomial $A_n(x)-xA_{n-2}(x)$ is the descent polynomial of the
$(n-2)$-stack-sortable permutations.
Br\"and\'en~\cite[Theorem~5.2]{Branden2006} proved that \(A_n(x)-cxA_{n-2}(x)\) has only simple real zeros
when \(c<2\). The \emph{Euler zigzag number} $E_n$ is the number of alternating
permutations of $[n]$, i.e., permutations satisfying
$$\pi(1)<\pi(2)>\pi(3)<\cdots.$$
Zhang~\cite{Zhang2015} studied
\(A_n+kxA_{n-2}\) for real \(k\). We write \(c=-k\). In our notation,
Zhang's \(r\)th tangent number is \(E_{2r-1}\).
\begin{conjecture}[{\cite[Conjecture~3.2]{Zhang2015}}]\label{conj-01}
Let \(n\ge3\) and \(r=\lfloor n/2\rfloor\). The polynomial \(A_n(x)-cxA_{n-2}(x)\) has only simple real zeros
if and only if
\[
  c<\frac{E_{2r+1}}{E_{2r-1}}
  \qquad\text{or}\qquad
  c>n(n-1).
\]
\end{conjecture}

The ratio \(E_{2r+1}/E_{2r-1}\) is at least \(2\). For \(r\ge2\),
count the alternating permutations of \([2r+1]\) with maximum \(2r+1\)
in position \(2\). There are \(2r\) choices for the first entry and
\(E_{2r-1}\) alternating arrangements of the remaining entries.
Reversing these permutations preserves alternation and moves the maximum
to position \(2r\). The two sets are disjoint, so
\(E_{2r+1}\ge4rE_{2r-1}>2E_{2r-1}\).
For \(r=1\), \(E_3/E_1=2\).

A related question asks when all zeros of \((1+x)A_n+kxA_{n-1}\)
lie in the open left half-plane.
A nonzero polynomial is \emph{Hurwitz stable} if every zero has negative
real part. It is \emph{weakly Hurwitz stable} if every zero has
nonpositive real part.
In~\cite[(3)]{YangZhang2015}, Yang and Zhang index the type~\(A\) Eulerian polynomials by Coxeter
rank, so their \(A_{n-1}(x)\) is our \(A_n(x)\).
In our notation, their Theorem~2.3 gives weak Hurwitz stability of
\((1+x)A_n+kxA_{n-1}\) for \(n\ge2\) and \(k\ge-n\).
They used the case \(k=-n\) and the Hermite--Biehler theorem to give
another proof of Hyatt's interlacing result for the two half Eulerian
polynomials of type~\(D\). Yang and Zhang also made the following conjecture.
\begin{conjecture}[{\cite[Conjecture~3.2]{YangZhang2015}}]\label{conj-02}
For any \(n\ge3\), the polynomial
\((1+x)A_n(x)+kxA_{n-1}(x)\) is Hurwitz stable if and only if
\[
  k>-2\frac{E_n}{E_{n-1}}.
\]
\end{conjecture}

Both conjectures involve a ratio of Euler zigzag numbers.
Yang and Zhang~\cite[p.~472]{YangZhang2015} asked how the
two conjectures are related. We prove both conjectures using the operators
\[
 \Lop_s h=(z^2-1)\frac{d}{dz}h+szh,\qquad s\in\mathbb R,\quad h\in\mathbb R[z].
\]
For polynomials of consecutive degrees, we allow the coefficients of
\(z+1\) and \(z-1\) to vary independently. For \(a,b\in\mathbb R\), set
\[
 \Lop_{a,b}h=(z^2-1)h'+\bigl(a(z+1)+b(z-1)\bigr)h.
\]
Then \(\Lop_s=\Lop_{s/2,s/2}\).
Theorem~\ref{thm:two-endpoint-intro} gives the exact parameter range
for preserving strict interlacing between polynomials of consecutive degrees.

A nonzero polynomial is \emph{real-rooted} if all its zeros are real.
We use the terminology of~\cite[Section~2]{Branden2006} for interlacing.
Let \(f\) and \(g\) be real-rooted polynomials of degrees \(d\) and
\(m\), respectively, with zeros
\(a_1\ge\cdots\ge a_d\) and \(b_1\ge\cdots\ge b_m\), counted with
multiplicity. We say that \(g\) \emph{alternates left of} \(f\) if
\(d=m\) and
\[
 a_1\ge b_1\ge a_2\ge b_2\ge\cdots\ge a_m\ge b_m.
\]
If \(d=m+1\) and the chain of inequalities continues with
\(b_m\ge a_{m+1}\), we say that \(g\) \emph{interlaces} \(f\).
We write \(g\preceq f\) in either case and \(g\prec f\) when all
inequalities are strict. By convention, a nonzero constant strictly
interlaces any linear polynomial. For interlacing on an interval,
restrict each list to the zeros in that interval. An empty list
strictly interlaces a list consisting of one zero.

Suppose \(F\) and \(G\) each have \(m\) positive zeros.
Let \(0<\alpha_1<\cdots<\alpha_m\) and
\(0<\beta_1<\cdots<\beta_m\) be the positive zeros of \(F\) and \(G\),
respectively. We say that these zeros \emph{alternate strictly} if
\[
 \alpha_1<\beta_1<\cdots<\alpha_m<\beta_m
 \quad\text{or}\quad
 \beta_1<\alpha_1<\cdots<\beta_m<\alpha_m.
\]

\Needspace{12\baselineskip}
We first compare polynomials \(F\) and \(G\) of the same degree,
with \(F\) vanishing at \(\pm1\) and all zeros of \(G\) in \((-1,1)\).
\begin{introtheorem}\label{thm:param-interlace}\label{mainthm}
Let $F$ and $G$ be real polynomials of the same degree $d\ge2$,
both even or both odd. Suppose that $F$ has $d$ simple zeros in
$[-1,1]$, including both endpoints, whereas $G$ has $d$ simple
zeros in $(-1,1)$. Assume that the positive zeros of \(F\) and \(G\)
alternate strictly. If
\begin{equation}\label{eq:param-range}
  \nu>0, \qquad 0\le \nu-\mu\le 2,
\end{equation}
then $\Lop_\mu F$ has $d+1$ simple zeros in $[-1,1]$,
including both endpoints, whereas $\Lop_\nu G$ has $d+1$
simple zeros in $(-1,1)$. Moreover, the positive zeros of
$\Lop_\mu F$ and $\Lop_\nu G$ alternate strictly.
\end{introtheorem}

For each difference \(\nu-\mu\) outside \([0,2]\), there is a choice
of \(\nu>0\) and a pair of cubic polynomials whose images have a common
zero in \((0,1)\) (Proposition~\ref{prop:param-sharp}).

\begin{introtheorem}\label{thm:two-endpoint-intro}
Fix \(a,b,\alpha,\beta>0\). The relation
\(\Lop_{\alpha,\beta}F\prec\Lop_{a,b}G\) holds for every pair
\(F\prec G\) of real polynomials with positive leading coefficients,
\(\deg G=\deg F+1\), and all zeros in \((-1,1)\) if and only if
\(0\le\alpha-a\le1\) and \(0\le\beta-b\le1\).
For all positive parameters, both \(\Lop_{\alpha,\beta}F\) and
\(\Lop_{a,b}G\) have only simple zeros in \((-1,1)\).
\end{introtheorem}
For \(\mu,\nu>0\), setting \(\alpha=\beta=\mu/2\) and \(a=b=\nu/2\)
gives \(0\le\mu-\nu\le2\) as the exact condition for
\(\Lop_\mu F\prec\Lop_\nu G\) to hold for every consecutive-degree pair
in Theorem~\ref{thm:two-endpoint-intro}.

Theorem~\ref{thm:param-interlace} implies that for \(n\ge3\),
the zeros of \(A_{n-2}\) strictly interlace those of \(A_n\)
on each of \((-\infty,-1)\) and \((-1,0)\)
(Corollary~\ref{cor:typeA-two-step}).
Using this ordering, we show that the smaller endpoint in
Conjecture~\ref{conj-01} is the largest \(c\) for which
\((1+x)A_{n-2}(x)\) interlaces \(A_n(x)-cxA_{n-2}(x)\).
\begin{introtheorem}\label{thm:eulerian-fixed-interlacing}
Let \(n\ge3\), \(r=\lfloor n/2\rfloor\), and \(c\in\mathbb R\). Then
\[
 (1+x)A_{n-2}(x)\preceq A_n(x)-cxA_{n-2}(x)
 \quad\Longleftrightarrow\quad c\le\frac{E_{2r+1}}{E_{2r-1}}.
\]
\end{introtheorem}

For example, the upper bound on \(c\) in
Theorem~\ref{thm:eulerian-fixed-interlacing} is \(E_3/E_1=2\)
at \(n=3\) and \(E_5/E_3=8\) at \(n=4\). At these bounds,
\(A_3(x)-2xA_1(x)=(1+x)^2\) and \(A_4(x)-8xA_2(x)=(1+x)^3\).
In general, when \(c=E_{2r+1}/E_{2r-1}\),
\(-1\) is a double zero of \(A_n(x)-cxA_{n-2}(x)\) if \(n\) is odd
and a triple zero if \(n\) is even.
Corollary~\ref{cor:zhang-boundaries} also gives the second endpoint
\(c=n(n-1)\), at which \(1\) is a double zero. The polynomial
$A_n(x)-cxA_{n-2}(x)$ has only real zeros precisely when
$c\le E_{2r+1}/E_{2r-1}$ or $c\ge n(n-1)$, where $r=\lfloor n/2\rfloor$.

We now return to question~$(ii)$. For each fixed \(n\ge3\), we
seek the values of \(c\) for which \(A_n(x)+bA_{n-1}(x)+cA_{n-2}(x)\)
has only real zeros for every real \(b\). The largest zero of \(A_{n-1}\)
determines the sharp upper bound on \(c\). The largest zero of the middle
polynomial also determines the sharp bound in the other four families in
Table~\ref{tab:main-bounds}. In all five families we use the normalization
\(F_j(0)=1\). Let \(S_j\), \(N_j\), and \(M_j\) denote the second-order
Eulerian, ordinary Narayana, and type~\(B\) Narayana polynomials,
respectively, as defined in Section~\ref{sec:three-term}.

\begin{introtheorem}\label{thm:largest-intro}
Let \((F_j)\) be one of the five polynomial sequences listed in
Table~\ref{tab:main-bounds}.
Fix $n$ such that $F_{n-1}$ is nonconstant, and let $\rho$
denote its largest zero. For a fixed
real number $c$, the polynomial $F_n(x)+bF_{n-1}(x)+cF_{n-2}(x)$
has only real zeros for every $b\in\mathbb R$ if and only if
\[
  c\le -\frac{F_n(\rho)}{F_{n-2}(\rho)}.
\]
\end{introtheorem}
The denominator \(F_{n-2}(\rho)\) is nonzero because
\(F_{n-2}\prec F_{n-1}\).
The proof uses residues to test real-rootedness for all \(b\) by
checking inequalities at the zeros of \(F_{n-1}\).
The residue test also allows the coefficient of \(F_{n-2}\) to be \(c-kx\),
as needed for the combinations with cycle weights in
Section~\ref{sec:cycle-interlacing}.

\begin{table}[htbp]
\centering
\caption{Sharp upper bounds on \(c\) for real-rootedness of
\(F_n(x)+bF_{n-1}(x)+cF_{n-2}(x)\), uniform over \(b\in\mathbb R\)
and the indices \(n\) for which \(F_{n-1}\) is nonconstant.}
\label{tab:main-bounds}
\begin{tabular}{@{}lcc@{}}
\toprule
Family & Sharp upper bound on \(c\) & Index attaining the bound \\
\midrule
Eulerian \(A_n\)
  & \(12-6\sqrt3\) & \(4\) \\[3pt]
Type \(B\) Eulerian \(B_n\)
  & \(60(\sqrt6-\sqrt5)^2\) & \(4\) \\[3pt]
Second-order Eulerian \(S_n\)
  & \(\frac{5}{3}(4-\sqrt{10})\) & \(4\) \\[3pt]
Narayana \(N_n\)
  & \(\frac{35-7\sqrt{21}}{4}\) & \(5\) \\[3pt]
Type \(B\) Narayana \(M_n\)
  & \(\displaystyle
       \frac{210}{\bigl(3\sqrt7+\sqrt{35+2\sqrt{70}}\bigr)^2}\)
  & \(6\) \\
\bottomrule
\end{tabular}
\end{table}

For each family, Table~\ref{tab:main-bounds} gives the minimum of
the bounds in Theorem~\ref{thm:largest-intro} over all \(n\) for which
\(F_{n-1}\) is nonconstant, together with the unique minimizing index.
\begin{example}
Consider the Eulerian polynomials $A_n(x)$. For $n=3$, the
polynomial $A_2(x)=1+x$ has the single zero $-1$. Since
$A_1(x)=1$ and $A_3(x)=1+4x+x^2$, Theorem~\ref{thm:largest-intro}
gives $c\le -A_3(-1)/A_1(-1)=2$.
For $n=4$, the larger zero of $A_3(x)$ is
$\rho=-2+\sqrt3$. Using
$A_4(x)=A_2(x)\bigl(A_3(x)+6x\bigr)$,
we obtain
\[
 -\frac{A_4(\rho)}{A_2(\rho)}=-6\rho=12-6\sqrt3.
\]
Table~\ref{tab:main-bounds} shows that the bound at \(n=4\) is the
smallest. Thus \(A_n+bA_{n-1}+cA_{n-2}\) has only real zeros for all
\(n\ge3\) and \(b\in\mathbb R\) if and only if \(c\le12-6\sqrt3\).
\end{example}

We next give each cycle of a permutation weight \(q\). Let
\(\operatorname{cyc}(\pi)\) count the cycles of \(\pi\), including fixed
points, and define the \(q\)-Eulerian polynomials by
\[
 A_n(x;q)=\sum_{\pi\in\mathfrak S_n}
 x^{\excc(\pi)}q^{\operatorname{cyc}(\pi)},
 \qquad A_0(x;q)=1.
\]
In particular, \(A_n(x;1)=A_n(x)\).

An ordered list of values for the coefficient of \(A_{n-1}(x;q)\)
gives a row of three-term combinations. To study such rows, we use the notion of
full interlacing due to Athanasiadis and
Wagner~\cite[Definition~3.5]{AthanasiadisWagner2026}.
Let \((P_{ij}(x))\) be a \(p\times s\) matrix of polynomials.
For \(u,v\in\mathbb Z\) and \(1\le i\le p,\ 1\le j\le s\),
put \([x^{v-u}]P_{ij}(x)\) in row \(pu+i-1\) and column
\(sv+j-1\) of an infinite matrix. Coefficients of negative
powers are zero. The polynomial matrix is \emph{fully interlacing}
if this infinite matrix is totally nonnegative, meaning that all its
finite square minors are nonnegative. Athanasiadis and Wagner use the
term totally positive in this sense. Each nonzero entry of a fully
interlacing matrix has nonnegative coefficients and only real zeros.
In particular, for a fully interlacing row \((P_1(x),\ldots,P_s(x))\)
with nonzero entries, we have \(P_j\preceq P_i\) for \(i<j\).

At each fixed index, the bound on \(c\) in Theorem~\ref{thm:largest-intro}
is also the sharp bound for full interlacing of all rows with middle
coefficients \(-1-c\le b_1\le\cdots\le b_s\)
(Proposition~\ref{prop:five-family-full-interlacing}).

\begin{introtheorem}\label{thm:cycle-full-interlacing}
Fix \(q>0\) and \(c\in\mathbb R\). Given \(n\ge3\) and a nonempty
finite sequence \(-q\le b_1\le\cdots\le b_s\), form the row
\[
 \bigl(A_n(x;q)+b_iA_{n-1}(x;q)-cxA_{n-2}(x;q)\bigr)_{i=1}^{s}.
\]
All such rows are fully interlacing
if and only if \(c\le q+1\).
\end{introtheorem}
The lower bound \(-q\) on each \(b_i\) is necessary for nonnegative
coefficients. For \(n\ge4\) and \(c\le q+1\), the polynomial
\(A_n(x;q)+bA_{n-1}(x;q)-cxA_{n-2}(x;q)\) has a symmetric matrix
representation for every real \(b\). Varying \(b\) gives a rank-one perturbation of the
matrix, and each zero is a strictly decreasing function of \(b\).

Section~\ref{sec:two-parameter} proves the operator theorems.
Sections~\ref{sec:Zhang}--\ref{sec:ABD} apply them to the two
conjectures, stack sorting, and Eulerian polynomials of types~\(B\)
and~\(D\). In Sections~\ref{sec:eulerian-fixed-interlacing} and
\ref{sec:three-term}, we prove Theorems~\ref{thm:eulerian-fixed-interlacing}
and~\ref{thm:largest-intro}, respectively.
Section~\ref{sec:cycle-interlacing} treats
cycle weights and the $q$-Eulerian transformation.

\section{Proofs of Theorems~\ref{thm:param-interlace} and~\ref{thm:two-endpoint-intro}}\label{sec:two-parameter}
The Cayley transformation takes the negative real axis to \((-1,1)\)
and sends \(-1\) to \(0\). In the Cayley variable, the Eulerian recurrence
becomes an application of \(\Lop_2\), which we use to compare the zeros
of \(A_n\) and \(A_{n-2}\) on each side of \(-1\).

\subsection{The Cayley transformation}
%\hspace*{\parindent}

For \(n\ge1\) set \(C_n(z)=(z+1)^{n-1}A_n(x)\), where
\(x=(z-1)/(z+1)\). Since \(\deg A_n=n-1\), the factor
\((z+1)^{n-1}\) clears the denominator. The product thus defines
a polynomial \(C_n\), including at \(z=-1\). We have \(C_1=1\). Multiplying
\eqref{eq:Arec} by \((z+1)^n\) and using
\(dx/dz=2/(z+1)^2\) gives
\begin{align*}
 C_{n+1}(z)
 &= (z+1)^n\left((1+nx)A_n(x)+x(1-x)\frac{d}{dx}A_n(x)\right)\\
 &= \frac{d}{dz}\bigl((z-1)(z+1)^nA_n(x)\bigr)\\
 &= \frac{d}{dz}\bigl((z^2-1)C_n(z)\bigr)\\
 &= (z^2-1)\frac{d}{dz}C_n(z)+2zC_n(z).
\end{align*}
Starting with \(C_1=1\), the recurrence gives
\begin{align*}
 C_2(z)&=2z,~~C_3(z)=6z^2-2,~~C_4(z)=24z^3-16z,~~C_5(z)=120z^4-120z^2+16.
\end{align*}

In terms of the operator \(\Lop_s\) from the Introduction, this recurrence is
\begin{equation}\label{L2C2}
\Lop_2 C_n=C_{n+1}.
\end{equation}
Each application of \(\Lop_2\) reverses parity. Hence \(C_n\)
is even when \(n\) is odd and odd when \(n\) is even.
By the product rule,~\eqref{L2C2} also gives
\(\frac{d}{dz}\bigl((z^2-1)C_{n-2}(z)\bigr)=C_{n-1}(z)\)
for \(n\ge3\). Multiplying by \(z^2-1\) gives
\(\Lop_0\bigl((z^2-1)C_{n-2}\bigr)=(z^2-1)C_{n-1}\).
Since \(C_n\) and \((z^2-1)C_{n-2}\) have the same parity, we
compare their positive zeros and check the multiplicity at zero.
Reflection then gives the order of their negative zeros.

Since \((z+1)^2x=z^2-1\), the combination becomes
\begin{equation}\label{AC}
 (z+1)^{n-1}\bigl(A_n(x)-cxA_{n-2}(x)\bigr)
 =C_n(z)-c(z^2-1)C_{n-2}(z).
\end{equation}
We compare the positive zeros of \(C_n\) and \((z^2-1)C_{n-2}\)
by induction on \(n\).
For \(n=3\), the positive zeros of \((z^2-1)C_1(z)\) and \(C_3(z)\)
are \(1\) and \(1/\sqrt3\), respectively. To pass from \(n\) to
\(n+1\), we apply \(\Lop_0\) to \((z^2-1)C_{n-2}\) and \(\Lop_2\)
to \(C_n\), as in Theorem~\ref{thm:param-interlace} with
\((\mu,\nu)=(0,2)\).

\subsection{Proof of Theorem~\ref{thm:param-interlace}}
%\hspace*{\parindent}

We denote the leading coefficient of a nonzero polynomial \(f\) by
\(\lc(f)\) and the sign of a real number \(t\) by \(\sgn(t)\).
Write \(F(z)=z^\varepsilon f(y)\) and \(G(z)=z^\varepsilon g(y)\),
where \(y=z^2\) and \(\varepsilon=0\) in the even case and \(1\) in
the odd case.
Squaring the positive zeros preserves their order. The next lemma gives
the signs of the residues of \(g/f\) when both leading coefficients
are positive. We use these residues to determine the
sign of \(\Lop_\nu G\) at each zero of \(\Lop_\mu F\) in \((0,1)\).

\begin{lemma}\label{lem:residue-representation}
Let \(f,g\in\mathbb R[y]\) have degree \(m\ge1\), positive leading
coefficients, and only simple positive zeros. Suppose \(g\) strictly
alternates left of \(f\), whose zeros are \(1=a_1>\cdots>a_m>0\).
Then
\begin{equation}\label{eq:residue-representation}
 \frac{g(y)}{f(y)}
 =\frac{\lc(g)}{\lc(f)}
  +\sum_{i=1}^m\frac{g(a_i)}{f'(a_i)}\frac{1}{y-a_i},
 \qquad \frac{g(a_i)}{f'(a_i)}>0 \quad(1\le i\le m).
\end{equation}
In particular, for real \(y\) with \(f(y)\ne0\),
\begin{equation}\label{eq:Rprime-negative}
 \left(\frac{g}{f}\right)'(y)
 =-\sum_{i=1}^m\frac{g(a_i)}{f'(a_i)}\frac{1}{(y-a_i)^2}<0.
\end{equation}
\end{lemma}
\begin{proof}
Since the zeros of \(f\) are simple and \(g(y)/f(y)\) tends to
\(\lc(g)/\lc(f)\) as \(y\to\infty\), partial fractions give
\[
 \frac{g(y)}{f(y)}
 =\frac{\lc(g)}{\lc(f)}+\sum_{i=1}^m\frac{\rho_i}{y-a_i},
 \qquad
 \rho_i=\lim_{y\to a_i}(y-a_i)\frac{g(y)}{f(y)}
       =\frac{g(a_i)}{f'(a_i)}.
\]

Write \(b_1>\cdots>b_m\) for the zeros of \(g\). By hypothesis,
\(a_1>b_1>\cdots>a_m>b_m>0\). In each of the two products
\[
 f'(a_i)=\lc(f)\prod_{j\ne i}(a_i-a_j),
 \qquad
 g(a_i)=\lc(g)\prod_{j=1}^m(a_i-b_j),
\]
the negative factors are precisely those with \(j<i\).
Thus \(f'(a_i)\) and \(g(a_i)\) both have sign \((-1)^{i-1}\),
so every residue \(\rho_i\) is positive. Differentiation now
gives~\eqref{eq:Rprime-negative}.
\end{proof}

\begin{lemma}\label{lem:param-residue}
Let \(f\) and \(g\) satisfy the hypotheses of
Lemma~\ref{lem:residue-representation}. If \(0\le\nu-\mu\le2\), then
\[
 (\nu-\mu)\frac{g(y)}{f(y)}
 +2(y-1)\left(\frac{g}{f}\right)'(y)>0
 \qquad (0<y<1,\ f(y)\ne0).
\]
\end{lemma}
\begin{proof}
For fixed \(y\), the left-hand side is affine in \(\nu-\mu\).
At \(\nu=\mu\) it is positive because \(y<1\) and
\((g/f)'(y)<0\) by~\eqref{eq:Rprime-negative}.
At \(\nu=\mu+2\) it equals twice the derivative of
\((y-1)g(y)/f(y)\). Since the derivative of $(y-1)/(y-a_i)$ is
$(1-a_i)/(y-a_i)^2$, equation~\eqref{eq:residue-representation} gives
\[
 \frac{d}{dy}\left((y-1)\frac{g(y)}{f(y)}\right)
 =\frac{\lc(g)}{\lc(f)}
   +\sum_{i=1}^m\frac{g(a_i)}{f'(a_i)}
                 \frac{1-a_i}{(y-a_i)^2}>0.
\]
Each summand is nonnegative since the residues are positive and
\(a_i\le1\). The leading-coefficient ratio is positive. For each
fixed \(y\), positivity
at \(\nu-\mu=0\) and \(\nu-\mu=2\) gives positivity throughout
the interval of parameter differences.
\end{proof}

For \(f(y)=y-1\) and \(g(y)=y-\tfrac12\), the quotient \(g/f\) is
negative on \((\tfrac12,1)\). At \(\nu-\mu=2\), the left-hand side
in Lemma~\ref{lem:param-residue} is
\(2\frac{d}{dy}\bigl((y-1)g(y)/f(y)\bigr)=2\).

\hypersetup{bookmarksdepth=3}
\begin{proof}[Proof of Theorem~\ref{thm:param-interlace}]
Multiplying \(F\) or \(G\) by \(-1\) preserves the zeros of that
polynomial and of its image. We may therefore assume that both have
positive leading coefficients. Write \(d=2m+\varepsilon\), where
\(\varepsilon\in\{0,1\}\), and put
\(y=z^2\), \(F(z)=z^\varepsilon f(y)\), and \(G(z)=z^\varepsilon g(y)\).
The nonzero zeros of each polynomial occur in opposite pairs, and all
zeros are simple. Thus \(f(0)g(0)\ne0\), and the zeros of \(f\)
and \(g\) are the squares of the positive zeros of \(F\) and \(G\),
respectively. Squaring preserves the order of the positive zeros,
and \(1\) is the largest zero of \(F\). Therefore
\[
 1=a_1>b_1>a_2>b_2>\cdots>a_m>b_m>0,
\]
where the \(a_i\) are the zeros of \(f\) and the \(b_i\) are the
zeros of \(g\). In particular,
\(f'(1)=\lc(f)\prod_{i=2}^m(1-a_i)>0\), also when \(m=1\)
and the product is empty. Write \(R=g/f\). By~\eqref{eq:param-range},
\(\mu\ge\nu-2>-2\).

For each parity we first locate the positive zeros of \(\Lop_\mu F\).
We then use Lemma~\ref{lem:param-residue} to determine the sign of
\(\Lop_\nu G\) at each of these zeros in \((0,1)\).

\medskip
\noindent\emph{Even case.}
Here \(F(z)=f(y)\) and \(G(z)=g(y)\). Since \(F'(z)=2zf'(y)\),
we have
\[
\begin{aligned}
 \Lop_\mu F&=z\,p(y), & p(y)&=2(y-1)f'(y)+\mu f(y),\\
 \Lop_\nu G&=z\,q(y), & q(y)&=2(y-1)g'(y)+\nu g(y).
\end{aligned}
\]
The leading coefficients of \(p\) and \(q\) are
\((2m+\mu)\lc(f)\) and \((2m+\nu)\lc(g)\), respectively.
Both are positive, so \(p\) and \(q\) have degree \(m\).

Since \(f(1)=0\), the formula for \(p\) gives \(p(1)=0\).
Differentiation at \(y=1\) gives $p'(1)=2f'(1)+\mu f'(1)=(\mu+2)f'(1)>0$.
Thus \(1\) is a simple zero of \(p\), and the polynomial
\(p(y)/(y-1)\) has value \(p'(1)>0\) at \(1\).
At each remaining zero $a_i$ of \(f\), where $2\le i\le m$,
we have $p(a_i)/(a_i-1)=2f'(a_i)$.
Since precisely \(i-1\) zeros of \(f\) exceed \(a_i\),
\(\sgn f'(a_i)=(-1)^{i-1}\). Hence the values of
\(p(y)/(y-1)\) at \(a_1,\ldots,a_m\) alternate in sign.
The quotient \(p(y)/(y-1)\) has degree \(m-1\), so it has exactly one
simple zero \(\xi_i\) in \((a_i,a_{i-1})\) for each \(2\le i\le m\).
The zeros of \(p\) are \(\xi_1=1\) and the zeros of the quotient.

For \(f(y)\ne0\), substituting \(g=fR\) into the formula for \(q\)
gives
\[
 q(y)=R(y)p(y)+f(y)\bigl((\nu-\mu)R(y)+2(y-1)R'(y)\bigr).
\]
For \(2\le i\le m\), we have \(p(\xi_i)=0\) and
\(a_i<\xi_i<a_{i-1}\). In particular, \(0<\xi_i<1\) and
\(f(\xi_i)\ne0\). The expression in parentheses is positive at
\(\xi_i\) by Lemma~\ref{lem:param-residue}. Therefore
\(\sgn q(\xi_i) =\sgn f(\xi_i)=(-1)^{i-1}\).
At \(y=1\) we have \(q(1)=\nu g(1)>0\). At \(y=0\), the product
\(g(y)=\lc(g)\prod_{j=1}^m(y-b_j)\) gives
\[
 \frac{g'(0)}{g(0)}=-\sum_{j=1}^m\frac1{b_j},
 \qquad
 q(0)=g(0)\left(\nu+2\sum_{j=1}^m\frac1{b_j}\right).
\]
Thus \(\sgn q(0)=(-1)^m\). The values of \(q\) at
\(\xi_1,\ldots,\xi_m,0\) alternate in sign. Since \(\deg q=m\),
there is exactly one simple zero of \(q\) between each two
successive points. When \(m=1\), these two points are \(1\) and \(0\).

\medskip
\noindent\emph{Odd case.}
Here \(F(z)=zf(y)\) and \(G(z)=zg(y)\). Since
\(F'(z)=f(y)+2yf'(y)\), we have
\[
\begin{aligned}
 \Lop_\mu F&=p(y), & p(y)&=((\mu+1)y-1)f(y)+2y(y-1)f'(y),\\
 \Lop_\nu G&=q(y), & q(y)&=((\nu+1)y-1)g(y)+2y(y-1)g'(y).
\end{aligned}
\]
The leading coefficients are now \((2m+1+\mu)\lc(f)\) and
\((2m+1+\nu)\lc(g)\), both positive. Thus \(p\) and \(q\) have
degree \(m+1\).

Here again \(p(1)=0\). Differentiating the expression for \(p\)
at \(y=1\) gives \(p'(1)=\mu f'(1)+2f'(1)>0\), since \(f(1)=0\).
Thus \(1\) is a simple zero of \(p\) and \(p(y)/(y-1)\) is
positive at \(1\). At \(a_2,\ldots,a_m\) and at \(0\), we have
\[
 \frac{p(a_i)}{a_i-1}=2a_i f'(a_i)\quad(2\le i\le m),
 \qquad \frac{p(0)}{0-1}=f(0).
\]
Since \(a_i>0\) and \(\sgn f(0)=(-1)^m\), the
values of \(p(y)/(y-1)\) at \(a_1,\ldots,a_m,0\) alternate in
sign. The quotient has degree \(m\), so it has exactly one simple zero in
each of the \(m\) intervals between successive points. Denote
these zeros by \(\xi_1,\ldots,\xi_m\), where
\(a_i>\xi_i>a_{i+1}\) for \(1\le i<m\) and \(a_m>\xi_m>0\).

Substituting \(g=fR\) into the formula for \(q\) now gives
\[
 q(y)=R(y)p(y)+y f(y)\bigl((\nu-\mu)R(y)+2(y-1)R'(y)\bigr)
 \qquad(f(y)\ne0).
\]
For each \(1\le i\le m\), we have \(p(\xi_i)=0\),
\(0<\xi_i<1\), and exactly \(i\) zeros of \(f\) exceed \(\xi_i\).
At \(y=\xi_i\) the term \(R(y)p(y)\) vanishes. The factor
\(y\) is positive, and the factor in parentheses is positive by
Lemma~\ref{lem:param-residue}. Hence
\(\sgn q(\xi_i) =\sgn f(\xi_i)=(-1)^i\).
At the endpoints, \(q(1)=\nu g(1)>0\) and \(q(0)=-g(0)\) has
sign \((-1)^{m+1}\). Thus the values of \(q\) at
\(1,\xi_1,\ldots,\xi_m,0\) alternate in sign. Since
\(\deg q=m+1\), there is exactly one simple zero of \(q\) between
each two successive points.

In the variable \(y\), the zeros of \(p\) lie in \((0,1]\), those of
\(q\) lie in \((0,1)\), and \(q\) strictly alternates left of \(p\).
Taking positive square roots preserves the order. Each simple zero
\(y\) of \(p\) or \(q\) gives two simple zeros \(z=\pm\sqrt y\) of
\(\Lop_\mu F\) or \(\Lop_\nu G\), respectively.
When \(F\) and \(G\) are even, both images also have a simple zero
at \(z=0\), since \(p(0)q(0)\ne0\). Thus \(\Lop_\mu F\) has \(d+1\)
simple zeros in \([-1,1]\), including both endpoints, and
\(\Lop_\nu G\) has \(d+1\) simple zeros in \((-1,1)\).
Their positive zeros alternate strictly.
\end{proof}

\begin{example}\label{ex:smallest-operator-pair}
Take \((\mu,\nu)=(0,2)\), so that \(\nu-\mu=2\). The even polynomials
\((z^2-1)C_1=z^2-1\) and \(C_3=6z^2-2\) have positive zeros
\(1\) and \(1/\sqrt3\), respectively. Their images are
\[
 \begin{aligned}
 (z^2-1)C_2
 &=\Lop_0\bigl((z^2-1)C_1\bigr)
  =(z^2-1)\frac{d}{dz}(z^2-1)=2z(z^2-1),\\
 C_4
 &=\Lop_2C_3
  =(z^2-1)C_3'+2zC_3=8z(3z^2-2).
 \end{aligned}
\]
Both \((z^2-1)C_2\) and \(C_4\) have a simple zero at \(0\). Their positive
zeros are \(\sqrt{2/3}\) for \(C_4\) and \(1\) for \((z^2-1)C_2\),
so the positive zeros alternate strictly. Applying the operators
to this odd pair gives
\[
 \begin{aligned}
 (z^2-1)C_3
 &=\Lop_0\bigl((z^2-1)C_2\bigr)
  =(z^2-1)\frac{d}{dz}\bigl(2z(z^2-1)\bigr)=2(z^2-1)(3z^2-1),\\
 C_5
 &=\Lop_2C_4
  =(z^2-1)C_4'+2zC_4=120z^4-120z^2+16.
 \end{aligned}
\]
Here
\(C_5(0)=C_5(1)=16>0\) and \(C_5(1/\sqrt3)=-32/3<0\).
Since \(C_5\) is even of degree four, these signs give exactly
one simple zero in each of \((0,1/\sqrt3)\) and
\((1/\sqrt3,1)\), and the other two zeros are their negatives.
\end{example}

Neither bound on \(\nu-\mu\) in Theorem~\ref{thm:param-interlace}
can be improved uniformly in the parameters. For each difference
outside \([0,2]\), cubic polynomials give a counterexample for some
choice of parameters with that difference.
\begin{proposition}\label{prop:param-sharp}
Fix \(\Delta\in\mathbb R\setminus[0,2]\). There exist real numbers
\(\mu,\nu,b\) with \(\nu>0\), \(\nu-\mu=\Delta\), and \(0<b<1\)
such that $\Lop_\mu\bigl(z(z^2-1)\bigr)$
and $\Lop_\nu\bigl(z(z^2-b)\bigr)$
have a common zero in \((0,1)\). When \(\nu=0\), both images
vanish at \(1\) for \(\mu\in\mathbb R\) and \(b\in(0,1)\).
\end{proposition}

\begin{proof}
Let $F(z)=z(z^2-1)$ and $G(z)=z(z^2-b)$, where $0<b<1$.
Their positive zeros are \(1\) and \(\sqrt b\), respectively.
If \(\mu+3>0\), put \(a=(\mu+3)^{-1}\). Then
\(\Lop_\mu F=(z^2-1)((\mu+3)z^2-1)\) vanishes at \(\sqrt a\).
Using \((\nu+1)a-1=(\Delta-2)a\), we obtain
\[
 \begin{aligned}
 (\Lop_\nu G)(\sqrt a)
 &=((\nu+1)a-1)(a-b)+2a(a-1)\\
 &=a\bigl(\Delta a-(\Delta-2)b-2\bigr).
 \end{aligned}
\]
First fix \(\Delta>2\). With \(b=1/2\), \(\Lop_\nu G\)
vanishes at \(\sqrt a\) when \(a=1/2+1/\Delta\in(0,1)\).
This gives
\[
 \mu=\frac1a-3=-\frac{\Delta+6}{\Delta+2},
 \qquad
 \nu=\mu+\Delta
     =\frac{(\Delta-2)(\Delta+3)}{\Delta+2}>0.
\]
The common positive zero \(\sqrt a\) rules out strict alternation.

Now fix \(\Delta<0\). Choose any \(\nu>0\) and set \(\mu=\nu-\Delta>0\),
so \(0<a<1\). Solving \((\Lop_\nu G)(\sqrt a)=0\) for \(b\) gives
\(b=a+2(1-a)/(2-\Delta)\). Since \(0<2/(2-\Delta)<1\), we have
\(a<b<1\), and the two images again share the zero \(\sqrt a\).

Since \(F(1)=0\), we have \((\Lop_\mu F)(1)=0\) and
\((\Lop_\nu G)(1)=\nu G(1)\). Thus for \(\nu=0\) the two images
have the common positive zero \(1\). In particular, strict alternation
fails at \(\mu=\nu=0\), where \(\nu-\mu\) lies in \([0,2]\).
\end{proof}

\subsection{Polynomials of consecutive degrees}
\label{sec:consecutive-operators}

Theorem~\ref{thm:param-interlace} treats polynomials of the same
degree. We now let $\Lop_\mu$ act on the lower-degree member of
$f\prec g$.

\begin{theorem}\label{thm:parameter-interlacing}
Fix \(\mu,\nu>0\). Let \(f\) and \(g\) be real polynomials with all
zeros in \([-1,1]\). Assume \(f\prec g\) and \(\deg g=\deg f+1\). If
$0\le\mu-\nu\le2$, then \(\Lop_\mu f\prec\Lop_\nu g\).
In this range, all zeros of \(\Lop_\mu f\) are simple and lie in
\((-1,1)\). All zeros of \(\Lop_\nu g\) are simple and lie in \([-1,1]\).
\end{theorem}
\begin{proof}
Assume $0\le\mu-\nu\le2$ and put $d=\deg f$.
Scale $f$ and $g$ to have positive leading coefficients.
This does not change their zeros or the zeros of their images.
We first locate the zeros of $\Lop_\mu f$, then use partial
fractions to determine the signs of $\Lop_\nu g$ there.

Since $f\prec g$, 
all zeros of \(f\) and \(g\) are simple, and each zero of \(f\) lies
strictly between two consecutive zeros of \(g\). Thus for \(d\ge1\)
we may write the zeros of \(f\) as \(-1<a_1<\cdots<a_d<1\).
At each zero \(a_i\), the operator definition gives
$(\Lop_\mu f)(a_i)=(a_i^2-1)f'(a_i)$.
The product $f'(a_i)=\lc(f)\prod_{j\ne i}(a_i-a_j)$ has
exactly \(d-i\) negative factors. Thus
$\sgn f'(a_i)=(-1)^{d-i}$.
As \(a_i^2-1<0\), it follows that
$\sgn(\Lop_\mu f)(a_i)=(-1)^{d-i+1}$.

At the endpoints, \((\Lop_\mu f)(-1)=-\mu f(-1)\) and
\((\Lop_\mu f)(1)=\mu f(1)\). Since \(\mu>0\), \(\lc(f)>0\), and all
zeros of \(f\) lie in \((-1,1)\), these values have signs
\((-1)^{d+1}\) and \(+1\), respectively.
Set \(a_0=-1\) and \(a_{d+1}=1\).
Together with $\sgn(\Lop_\mu f)(a_i)=(-1)^{d-i+1}$,
these endpoint values show that \(\Lop_\mu f\) alternates
in sign at \(a_0,a_1,\ldots,a_{d+1}\). The intermediate value theorem
gives at least one zero in each interval $(a_{i-1},a_i)$ for
$1\le i\le d+1$.
The leading coefficient of \(\Lop_\mu f\) is
\((d+\mu)\lc(f)>0\), so its degree is \(d+1\).
These intervals therefore contain all its zeros, with exactly
one simple zero in each interval.
Write \(\xi_i\) for the zero in \((a_{i-1},a_i)\). For \(d\ge1\) we obtain
\[
 -1<\xi_1<a_1<\xi_2<\cdots<a_d<\xi_{d+1}<1.
\]
For \(d=0\), the identity \(\Lop_\mu f=\mu zf\) gives the
single simple zero \(\xi_1=0\).
Thus \(f\prec\Lop_\mu f\) in all cases, and
$\sgn f(\xi_i)=(-1)^{d+1-i}$ for
$1\le i\le d+1$.

We next determine the sign of \(\Lop_\nu g\) at each \(\xi_i\).
At every zero \(r\) of \(g\), strict interlacing gives $f(r)/g'(r)>0$.
Indeed, \(f(z)\) and \(g(z)/(z-r)\) have equally many zeros
to the right of \(r\). Their leading coefficients are positive,
so their values at \(r\), namely \(f(r)\) and \(g'(r)\),
have the same sign.
Since \(\deg f=\deg g-1\), we have the partial-fraction expansion
\[
 \frac{f(z)}{g(z)}
 =\sum_{g(r)=0}\frac{f(r)}{g'(r)}\frac1{z-r}.
\]
For \(g(z)\ne0\), the operator definition and differentiation
of this expansion give
\[
\begin{aligned}
 g\Lop_\mu f-f\Lop_\nu g
 &=g^2\left(
     (z^2-1)\left(\frac fg\right)'
     +(\mu-\nu)z\frac fg
   \right)\\
 &=g^2\sum_{g(r)=0}\frac{f(r)}{g'(r)}
   \frac{1-z^2+(\mu-\nu)z(z-r)}{(z-r)^2}.
\end{aligned}
\]
For \(-1<z<1\) and \(-1\le r\le1\), we have
\(zr\le|z|\) and \(z^2-|z|\le0\). Therefore
\[
\begin{aligned}
 1-z^2+(\mu-\nu)z(z-r)
 &\ge 1-z^2+(\mu-\nu)(z^2-|z|)\\
 &\ge 1-z^2+2(z^2-|z|)\\
 &=(1-|z|)^2>0.
\end{aligned}
\]
The first inequality uses \(\mu-\nu\ge0\), and the second
uses \(\mu-\nu\le2\).
Every summand is positive, so
$g\Lop_\mu f-f\Lop_\nu g>0$
whenever \(-1<z<1\) and \(g(z)\ne0\).
\Needspace{7\baselineskip}
If \(r\in(-1,1)\) is a zero of \(g\), direct substitution gives
\[
 \begin{aligned}
 (g\Lop_\mu f-f\Lop_\nu g)(r)
 &=-f(r)(r^2-1)g'(r)=(1-r^2)f(r)g'(r)>0.
 \end{aligned}
\]
Indeed, \(1-r^2>0\) and \(f(r)/g'(r)>0\).
Thus \(g\Lop_\mu f-f\Lop_\nu g>0\) throughout \((-1,1)\).
Evaluating at a zero \(\xi_i\) of \(\Lop_\mu f\) gives
$-f(\xi_i)(\Lop_\nu g)(\xi_i)>0$.
Consequently,
\(\sgn(\Lop_\nu g)(\xi_i)=(-1)^{d+2-i}\) for \(1\le i\le d+1\).
These values alternate in sign, so \(\Lop_\nu g\) has a zero
in each interval $(\xi_i,\xi_{i+1})$, where $1\le i\le d$.

At the left endpoint, \((\Lop_\nu g)(-1)=-\nu g(-1)\).
If \(g(-1)=0\), the image vanishes at \(-1\). Otherwise all \(d+1\)
zeros of \(g\) exceed \(-1\), so \((\Lop_\nu g)(-1)\) has sign
\((-1)^{d+2}\), opposite to the sign at \(\xi_1\). In either case
\(\Lop_\nu g\) has a zero in \([-1,\xi_1)\).
At the right endpoint, \((\Lop_\nu g)(1)=\nu g(1)\ge0\), while
\((\Lop_\nu g)(\xi_{d+1})<0\). This gives a zero in \((\xi_{d+1},1]\).
We find \(d+2\) distinct zeros of \(\Lop_\nu g\),
one in each interval
\([-1,\xi_1)\), \((\xi_i,\xi_{i+1})\) for \(1\le i\le d\), and \((\xi_{d+1},1]\).
Its leading coefficient is \((d+1+\nu)\lc(g)>0\), so its
degree is \(d+2\).
Each interval therefore contains exactly one simple zero of
\(\Lop_\nu g\), giving \(\Lop_\mu f\prec\Lop_\nu g\).
\end{proof}

For every \(\mu,\nu>0\) with
\(\mu-\nu\notin[0,2]\), a linear and a quadratic polynomial give
a counterexample to the interlacing conclusion of
Theorem~\ref{thm:parameter-interlacing}.

\begin{proposition}\label{prop:consecutive-sharpness}
Fix \(\mu,\nu>0\). If \(\mu<\nu\) or \(\mu>\nu+2\),
then there exist a linear polynomial \(f\) with its zero
in \((-1,1)\) and a quadratic polynomial \(g\) with two
simple zeros in \([-1,1]\) such that
\(f\prec g\) and \(\Lop_\mu f\not\prec\Lop_\nu g\).
\end{proposition}

\begin{proof}
For \(0<a<1\), take
$f(z)=z-a$.
The operator definition gives
\[
\begin{aligned}
 \Lop_\mu f
 &=(z^2-1)+\mu z(z-a)=(\mu+1)z^2-\mu az-1.
\end{aligned}
\]
Its constant term is negative and its leading coefficient is
positive, so it has one negative zero and one positive zero.
The positive zero \(\xi\) lies in \((a,1)\), since $(\Lop_\mu f)(a)=a^2-1<0$ and
 $(\Lop_\mu f)(1)=\mu(1-a)>0$.
Evaluating \(\Lop_\mu f\) at \(\xi\) gives
$\xi^2-1=-\mu\xi(\xi-a)$.

For a monic quadratic \(g\), the polynomial \(\Lop_\nu g\)
has degree three and leading coefficient \(\nu+2>0\).
If \(\Lop_\mu f\prec\Lop_\nu g\), then \(\xi\) lies
between the two largest zeros of \(\Lop_\nu g\).
Thus strict interlacing would require
\((\Lop_\nu g)(\xi)<0\).
In each of the two parameter ranges we choose \(a\) and \(g\) so that
\(f\prec g\) and \((\Lop_\nu g)(\xi)>0\).

Suppose first that \(\mu<\nu\).
Choose $\frac1{1+\nu-\mu}<a<1$.
For \(0<\eta<1-a\), put $g(z)=(z+1)(z-a-\eta)$. Its zeros are \(-1\)
and \(a+\eta\). Since $-1<a<a+\eta<1$, the zero \(a\) of \(f\) lies
strictly between them. Hence \(f\prec g\).

Differentiating \(g\) and factoring out \(z+1\), we obtain
\[
\begin{aligned}
 \Lop_\nu g
 &=(z^2-1)\bigl((z-a-\eta)+(z+1)\bigr)
   +\nu z(z+1)(z-a-\eta)\\
 &=(z+1)\left(
   ((\nu+1)z-1)(z-a-\eta)+z^2-1
   \right).
\end{aligned}
\]
At \(z=\xi\), substitute
\(\xi^2-1=-\mu\xi(\xi-a)\) to get
\[
\begin{aligned}
 (\Lop_\nu g)(\xi)
 &=(\xi+1)\left(
   ((\nu+1)\xi-1)(\xi-a-\eta)
   -\mu\xi(\xi-a)
   \right)\\
 &=(\xi+1)\left(
   (\xi-a)((1+\nu-\mu)\xi-1)
   -\eta((\nu+1)\xi-1)
   \right).
\end{aligned}
\]
Since $\xi>a>1/(1+\nu-\mu)$, both factors in $(\xi-a)((1+\nu-\mu)\xi-1)$
are positive.
With \(a\) and \(\xi\) fixed, the expression in parentheses
is positive at \(\eta=0\), and hence remains positive for
all sufficiently small \(\eta>0\).
As \(\xi+1>0\), we obtain
$(\Lop_\nu g)(\xi)>0$.
This contradicts the sign required by strict interlacing.

Now suppose that \(\mu>\nu+2\).
Choose $a$ with $1/(\mu-\nu-1)<a<1$.
For \(0<\eta<a\), put
$g(z)=(z-1)(z-a+\eta)$.
Since \(0<a-\eta<a<1\), the zero \(a\) of \(f\) lies strictly
between the zeros \(a-\eta\) and \(1\) of \(g\). Thus \(f\prec g\).
This time, differentiating \(g\) and factoring out \(z-1\)
gives
\[
\begin{aligned}
 \Lop_\nu g
 &=(z^2-1)\bigl((z-a+\eta)+(z-1)\bigr)
   +\nu z(z-1)(z-a+\eta)\\
 &=(z-1)\left(
   ((\nu+1)z+1)(z-a+\eta)+z^2-1
   \right).
\end{aligned}
\]
Using the same identity for \(\xi^2-1\), we obtain
\[
\begin{aligned}
 (\Lop_\nu g)(\xi)
 &=(\xi-1)\left(
   ((\nu+1)\xi+1)(\xi-a+\eta)
   -\mu\xi(\xi-a)
   \right)\\
 &=(1-\xi)\left(
   (\xi-a)((\mu-\nu-1)\xi-1)
   -\eta((\nu+1)\xi+1)
   \right).
\end{aligned}
\]
Here \(1-\xi>0\), and $\xi>a>1/(\mu-\nu-1)$
implies $(\xi-a)((\mu-\nu-1)\xi-1)>0$.
The expression in parentheses is therefore positive at
\(\eta=0\), and remains positive for all sufficiently small
\(\eta>0\).
Thus \((\Lop_\nu g)(\xi)>0\), again contradicting the sign
required by strict interlacing.
\end{proof}

\begin{example}\label{ex:parameter-interlacing}
Take $f(z)=z$ and $g(z)=(z+1)(8z-3)$.
The zero of \(f\) is \(0\), and the zeros of \(g\) are
\(-1\) and \(3/8\). Since $-1<0<3/8<1$, we have \(f\prec g\),
and the zeros satisfy the hypotheses of
Theorem~\ref{thm:parameter-interlacing}.

Choose \((\mu,\nu)=(3,1)\), so that \(\mu-\nu=2\).
Using \(f'(z)=1\) and \(g'(z)=16z+5\), we obtain
\[
\begin{aligned}
 \Lop_3 f
 &=(z^2-1)+3z^2=4z^2-1
  =(2z+1)(2z-1),\\
 \Lop_1 g
 &=(z^2-1)(16z+5)+z(z+1)(8z-3)=(z+1)(4z+1)(6z-5).
\end{aligned}
\]
Thus the zeros of \(\Lop_3 f\) are \(-1/2\) and \(1/2\),
whereas those of \(\Lop_1 g\) are \(-1,-1/4,5/6\).
Their order is
\(-1<-\frac12<-\frac14<\frac12<\frac56<1\).
Hence $\Lop_3 f\prec\Lop_1 g$ at the upper bound \(\mu-\nu=2\).
The polynomial \(g\) is neither even nor odd, and \(-1\) is a
common zero of \(g\) and \(\Lop_1 g\).
\end{example}

\subsection{The two-parameter operator}

We now let the coefficients of \(z+1\) and \(z-1\) vary independently
in the operator \(\Lop_{a,b}\) from the Introduction. Taking \(a=b=s/2\) gives
\(\Lop_s=\Lop_{s/2,s/2}\).

\begin{theorem}\label{thm:two-endpoint-interlacing}
Fix \(a,b,\alpha,\beta>0\). The relation
\(\Lop_{\alpha,\beta}F\prec\Lop_{a,b}G\) holds for every pair
\(F\prec G\) of real polynomials with positive leading coefficients,
\(\deg G=\deg F+1\), and all zeros in \((-1,1)\) if and only if
\(0\le\alpha-a\le1\) and \(0\le\beta-b\le1\).
For all positive parameters, both \(\Lop_{\alpha,\beta}F\) and
\(\Lop_{a,b}G\) have only simple zeros in \((-1,1)\).
\end{theorem}

\begin{proof}
Let \(d=\deg F\). If \(d\ge1\), write the zeros of \(F\) as
\(a_1<\cdots<a_d\). At each \(a_i\), the value
\((\Lop_{\alpha,\beta}F)(a_i)=(a_i^2-1)F'(a_i)\)
has sign \((-1)^{d-i+1}\). At the endpoints,
\((\Lop_{\alpha,\beta}F)(-1)=-2\beta F(-1)\) and \((\Lop_{\alpha,\beta}F)(1)=2\alpha F(1)\).
These values have signs \((-1)^{d+1}\) and \(+1\), respectively.
The intermediate value theorem gives one zero in each of the
\(d+1\) intervals determined by \(-1,a_1,\ldots,a_d,1\).
The leading coefficient of \(\Lop_{\alpha,\beta}F\) is
\((d+\alpha+\beta)\lc(F)>0\), so these are all its zeros and they are
simple. In increasing order they satisfy
\[
 -1<\xi_1<a_1<\xi_2<\cdots<a_d<\xi_{d+1}<1.
\]
If \(d=0\), then \((\Lop_{\alpha,\beta}F)(z)=F(0)\bigl((\alpha+\beta)z+\alpha-\beta\bigr)\).
Its unique zero is \(\xi_1=(\beta-\alpha)/(\alpha+\beta)\in(-1,1)\).
These conclusions hold for all \(\alpha,\beta>0\).
The same argument shows that \(\Lop_{a,b}G\) has \(d+2\) simple
zeros in \((-1,1)\) for all \(a,b>0\).

Now assume \(0\le\alpha-a\le1\) and \(0\le\beta-b\le1\).
To place the zeros of \(\Lop_{a,b}G\) relative to those of
\(\Lop_{\alpha,\beta}F\), we determine the sign of
\(\Lop_{a,b}G\) at each \(\xi_i\).
The difference \(\Lop_{\alpha,\beta}-\Lop_{a,b}\) is multiplication
by \(\Delta(z)=(\alpha-a)(z+1)+(\beta-b)(z-1)\).
At each zero \(r\) of \(G\), strict interlacing gives
\(\lambda_r=F(r)/G'(r)>0\).
Since \(\deg G=\deg F+1\), partial fractions give
\[
 \frac{F(z)}{G(z)}
 =\sum_{G(r)=0}\frac{\lambda_r}{z-r}.
\]
The quotient rule and the partial fraction expansion of \(F/G\) give,
for \(G(z)\ne0\),
\[
 \begin{aligned}
 \frac{G\Lop_{\alpha,\beta}F-F\Lop_{a,b}G}{G^2}
 &=(z^2-1)\left(\frac{F}{G}\right)'+\Delta(z)\frac{F}{G}\\
 &=\sum_{G(r)=0}\lambda_r
   \frac{1-z^2+\Delta(z)(z-r)}{(z-r)^2}.
 \end{aligned}
\]
For \(-1<z<1\), the parameter assumptions give
\((\alpha-a)(z+1)\in[0,z+1]\) and \((\beta-b)(z-1)\in[z-1,0]\).
Adding these bounds gives \(z-1\le\Delta(z)\le z+1\).
When \(\Delta(z)>0\), the inequality \(r<1\) gives
\(1-z^2+\Delta(z)(z-r)>(1-z)(1+z-\Delta(z))\ge0\).
When \(\Delta(z)<0\), the inequality \(r>-1\) gives
\(1-z^2+\Delta(z)(z-r)>(1+z)(1-z+\Delta(z))\ge0\).
If \(\Delta(z)=0\), the numerator is \(1-z^2>0\).
Thus every summand is positive for \(G(z)\ne0\), and
\[
 G(z)\Lop_{\alpha,\beta}F(z)
 -F(z)\Lop_{a,b}G(z)>0
 \qquad(-1<z<1,\ G(z)\ne0).
\]
At a zero \(r\) of \(G\), evaluation gives
\[
 \bigl(G\Lop_{\alpha,\beta}F-F\Lop_{a,b}G\bigr)(r)
 =(1-r^2)F(r)G'(r)>0.
\]

At a zero \(\xi_i\) of \(\Lop_{\alpha,\beta}F\), we obtain
\(-F(\xi_i)(\Lop_{a,b}G)(\xi_i)>0\).
Since
\(\sgn F(\xi_i)=(-1)^{d+1-i}\), it follows that
\(\sgn(\Lop_{a,b}G)(\xi_i)=(-1)^{d+2-i}\).
Also,
\((\Lop_{a,b}G)(-1)=-2bG(-1)\) and \((\Lop_{a,b}G)(1)=2aG(1)\).
The first value has sign \((-1)^{d+2}\), and the second is positive.
Thus \(\Lop_{a,b}G\) has a zero in each of the intervals
\((-1,\xi_1)\), \((\xi_i,\xi_{i+1})\) for \(1\le i\le d\), and \((\xi_{d+1},1)\).
Its degree is \(d+2\), with leading coefficient
\((d+1+a+b)\lc(G)>0\). These \(d+2\) zeros are therefore simple and give
\(\Lop_{\alpha,\beta}F\prec\Lop_{a,b}G\).

We now prove that both parameter bounds are necessary. If
\(\alpha-a\notin[0,1]\) or \(\beta-b\notin[0,1]\), the values
\(\Delta(1)=2(\alpha-a)\) and \(\Delta(-1)=-2(\beta-b)\)
show that there is a point \(z_0\in(-1,1)\), sufficiently close to
\(1\) or \(-1\), for which either
\[
 \Delta(z_0)>z_0+1
 \qquad\text{or}\qquad
 \Delta(z_0)<z_0-1.
\]
The first inequality gives \(1-z_0^2+\Delta(z_0)(z_0-1)<0\).
The second gives \(1-z_0^2+\Delta(z_0)(z_0+1)<0\).

We choose a linear polynomial \(F(z)=z-s\) so that
\(\Lop_{\alpha,\beta}F\) vanishes at \(z_0\). Put
\(P(z)=\alpha(z+1)+\beta(z-1)\) and \(s=z_0-(1-z_0^2)/P(z_0)\).
Since \(P(1)=2\alpha>0\), we have \(s<z_0<1\) for \(z_0\) near \(1\),
and \(s\to1\) as \(z_0\uparrow1\). Similarly, \(P(-1)=-2\beta<0\)
gives \(-1<z_0<s\) near \(-1\), and \(s\to-1\) as \(z_0\downarrow-1\).
Thus \(s\in(-1,1)\) in either case.
The definition of \(s\) gives \((\Lop_{\alpha,\beta}F)(z_0)=0\).

If \(\Delta(z_0)>z_0+1\), choose \(r_2\in(z_0,1)\) with \(r_2>s\)
sufficiently close to \(1\) that
\(1-z_0^2+\Delta(z_0)(z_0-r_2)<0\), and then let \(r_1\) approach
\(s\) from below. If \(\Delta(z_0)<z_0-1\), choose
\(r_1\in(-1,z_0)\) with \(r_1<s\) sufficiently close
to \(-1\) that
\(1-z_0^2+\Delta(z_0)(z_0-r_1)<0\), and then let \(r_2\) approach
\(s\) from above. In both cases we have
\(-1<r_1<s<r_2<1\).
Set \(G(z)=(z-r_1)(z-r_2)\). Then \(F\prec G\), and
\[
 \frac{F(z)}{G(z)}
 =\frac{\lambda_1}{z-r_1}+\frac{\lambda_2}{z-r_2},
 \qquad
 \lambda_1=\frac{s-r_1}{r_2-r_1},\quad
 \lambda_2=\frac{r_2-s}{r_2-r_1}.
\]
As \(r_1\to s^-\) in the first case, \(\lambda_1\to0\) and
\(\lambda_2\to1\). As \(r_2\to s^+\) in the second, \(\lambda_1\to1\)
and \(\lambda_2\to0\). Write \(r=r_2\) in the first case and \(r=r_1\)
in the second.
At \(z=z_0\), the summand for the fixed pole \(r\) in
\((G\Lop_{\alpha,\beta}F-F\Lop_{a,b}G)/G^2\) tends to
\[
 \frac{1-z_0^2+\Delta(z_0)(z_0-r)}{(z_0-r)^2}<0.
\]
The other summand tends to zero because the corresponding residue
of \(F/G\) tends to zero and its pole tends to \(s\ne z_0\).
It follows that
\(\bigl(G\Lop_{\alpha,\beta}F-F\Lop_{a,b}G\bigr)(z_0)<0\)
for suitable \(r_1,r_2\). Since \((\Lop_{\alpha,\beta}F)(z_0)=0\), we have \(F(z_0)(\Lop_{a,b}G)(z_0)>0\).

The polynomial \(\Lop_{\alpha,\beta}F\) is quadratic. Its two zeros
lie on opposite sides of \(s\), since its values at \(-1,s,1\)
have signs \(+,-,+\). The cubic \(\Lop_{a,b}G\) has positive leading
coefficient \(2+a+b\). Strict interlacing would make \(\Lop_{a,b}G\)
positive at the smaller zero of \(\Lop_{\alpha,\beta}F\) and negative
at the larger one.
Since \(F(z)=z-s\), the product \(F(z)(\Lop_{a,b}G)(z)\) would then be
negative at both zeros. Its positive value at \(z_0\) proves necessity.
\end{proof}

To apply Theorem~\ref{thm:two-endpoint-interlacing} to polynomials with
negative zeros, we work in the variable \(x\). For a polynomial \(f\) of degree \(d\)
and \(a,b\in\mathbb R\), put
\[
 T_d^{a,b}f(x)=\bigl(a+(d+b)x\bigr)f(x)+x(1-x)f'(x).
\]
Set \(x=(z-1)/(z+1)\). Since \((z+1)\,dx/dz=1-x\), the product rule gives
\[
 \frac{d}{dz}\bigl((z+1)^d f(x)\bigr)
 =(z+1)^{d-1}\bigl(d f(x)+(1-x)f'(x)\bigr).
\]
Substituting this derivative into \(\Lop_{a,b}\) gives
\begin{equation}\label{eq:cayley-two-endpoint}
 (z+1)^{d+1}T_d^{a,b}f\!\left(\frac{z-1}{z+1}\right)
 =\Lop_{a,b}\!\left((z+1)^d
       f\!\left(\frac{z-1}{z+1}\right)\right).
\end{equation}
The substitution \(x=(z-1)/(z+1)\) is strictly increasing from
\((-1,1)\) onto \((-\infty,0)\). Suppose \(f\) has degree \(d\)
and a positive leading coefficient and that all its zeros are negative.
The polynomial \((z+1)^d f((z-1)/(z+1))\) has degree \(d\) with
leading coefficient \(f(1)>0\). Its zeros all lie in \((-1,1)\).
Returning to \(x\) in Theorem~\ref{thm:two-endpoint-interlacing}
gives the following statement. Let \(f\) and \(g\) have degrees \(d\)
and \(d+1\), respectively, positive leading coefficients, and only
simple negative zeros. Assume \(f\prec g\). For positive \(a,b,\alpha,\beta\),
\[
 T_d^{\alpha,\beta}f\prec T_{d+1}^{a,b}g
 \qquad\text{if }0\le\alpha-a\le1,\quad0\le\beta-b\le1.
\]
Each image has only simple negative zeros.
The coefficient of \(x^{d+1}\) in \(T_d^{a,b}f\) is
\((d+b)\lc(f)-d\lc(f)=b\lc(f)>0\) for \(b>0\).
For the inverse substitution \(z=(1+x)/(1-x)\), we have
\[
 \frac{1-x}{2}(z-r)=\frac{(1-r)+(1+r)x}{2}.
\]
For \(-1<r<1\), this linear factor has
positive leading coefficient and zero \((r-1)/(r+1)<0\).
In the counterexamples for Theorem~\ref{thm:two-endpoint-interlacing},
the inverse substitution sends the input pair \(F,G\) to a strictly
interlacing pair with positive leading coefficients and simple negative
zeros. By~\eqref{eq:cayley-two-endpoint}, the corresponding operator
images also fail to interlace strictly in \(x\).
Thus the parameter bounds are necessary for strict interlacing of the
images over all strictly interlacing input pairs of degrees \(d,d+1\),
with \(d\ge0\), positive leading coefficients, and only simple negative
zeros. We use the interlacing criterion for \(T_d^{a,b}\) in
Section~\ref{sec:cycle-interlacing}.

The Cayley transformation also lets us compare the zeros of the
Eulerian polynomials of types \(A\) and \(B\). We first write the
type~\(B\) recurrence in the variable \(z\).
A signed permutation of order $n$ is a word \(w_1\cdots w_n\)
whose absolute values form a permutation of \([n]\).
Set \(w_0=0\), and call an index \(0\le i<n\) a descent
if \(w_i>w_{i+1}\). Let \(B_n(x)\) count signed permutations
by these descents. Brenti~\cite{Brenti1994} proved that
\(B_n(x)\) has only real zeros.
These polynomials satisfy
\begin{equation}\label{eq:Brec}
 B_{j+1}(x)
 =\bigl(1+(2j+1)x\bigr)B_j(x)
  +2x(1-x)B_j'(x),
\end{equation}
with \(B_0(x)=1\).
With \(x=(z-1)/(z+1)\), define
\(H_j(z)=2^{-j}(z+1)^jB_j(x)\).
Then \(H_0(z)=1\). Since \(dx/dz=2/(z+1)^2\), differentiation
gives
\[
 H_j'(z)-\frac{j}{z+1}H_j(z)
 =2^{1-j}(z+1)^{j-2}B_j'(x).
\]
Multiplying \eqref{eq:Brec} by
\(2^{-j-1}(z+1)^{j+1}\) therefore yields
\begin{equation}\label{eq:Hrec}
 \begin{aligned}
 H_{j+1}(z)
 &=\bigl((j+1)z-j\bigr)H_j(z)+(z^2-1)
   \left(H_j'(z)-\frac{j}{z+1}H_j(z)\right)\\
 &=(z^2-1)H_j'(z)+zH_j(z)\\
 &=\Lop_1H_j(z).
 \end{aligned}
\end{equation}
Both sides of \(H_{j+1}=\Lop_1H_j\) are polynomials, so the identity
also holds at \(z=-1\).
By~\eqref{L2C2} and~\eqref{eq:Hrec}, increasing \(j\) by one
applies \(\Lop_2\) to \(C_j\) and \(\Lop_1\) to \(H_j\).
These parameters satisfy \(\mu-\nu=1\), so
Theorem~\ref{thm:parameter-interlacing} applies.
\begin{corollary}[{\cite[Corollary~6.8]{Branden2006}}]
For every \(n\ge1\), we have $A_n(x)\prec B_n(x)$.
\end{corollary}
\begin{proof}
We have \(C_1(z)=1\prec z=H_1(z)\).
Theorem~\ref{thm:parameter-interlacing} with
\((\mu,\nu)=(2,1)\) gives by induction
\(C_n(z)\prec H_n(z)\) for \(n\ge1\). The theorem places the
zeros of \(C_n\) in \((-1,1)\) and those of \(H_n\) in \([-1,1]\),
and all these zeros are simple.
Evaluating \eqref{eq:Hrec} at \(z=\pm1\) and using \(H_0(z)=1\),
we obtain $H_n(-1)=(-1)^n$ and $H_n(1)=1$.
Thus the zeros of \(H_n\) also lie in \((-1,1)\).
The substitution \(x=(z-1)/(z+1)\) maps \((-1,1)\)
strictly increasingly onto \((-\infty,0)\).
For \(-1<z<1\), both \((z+1)^{n-1}\) and
\(2^{-n}(z+1)^n\) are positive. Hence, with
\(x=(z-1)/(z+1)\), we have
\[
 \begin{aligned}
 C_n(z)=0&\quad\Longleftrightarrow\quad A_n(x)=0,\\
 H_n(z)=0&\quad\Longleftrightarrow\quad B_n(x)=0.
 \end{aligned}
\]
The strict order of the zeros is therefore preserved, giving
\(A_n\prec B_n\).
\end{proof}

\section{Zhang's conjecture on simple real zeros}
\label{sec:Zhang}

\label{sec:applications-a}

For \(n\ge3\), the descent polynomial of the \((n-2)\)-stack-sortable
permutations is \(A_n(x)-xA_{n-2}(x)\). Zhang's Theorem~1.3, with
\(k=-c\), shows that \(A_n(x)-cxA_{n-2}(x)\) has only real zeros
for \(n>3\) and \(c\le n\)~\cite{Zhang2015}.
Conjecture~\ref{conj-01} asks for the exact range of \(c\) for which
all zeros of \(A_n-cxA_{n-2}\) are real and simple. We isolate a
reciprocal quadratic factor. The quotient by that factor has \(n-3\)
simple negative zeros for every \(c\). The quadratic has a double zero
at \(-1\) at the smaller endpoint and a double zero at \(1\) at the
larger endpoint. Between the endpoints the quadratic has two nonreal zeros.
When \(n\) is even, the quotient
has a simple zero at \(-1\), so the original polynomial has a triple
zero there at the smaller endpoint.

\begin{theorem}\label{thm:Zhang-main}
Let \(n\ge3\), \(r=\lfloor n/2\rfloor\), and \(c\in\mathbb R\).
The polynomial \(A_n(x)-cxA_{n-2}(x)\) has only simple real zeros
if and only if
\[
 c<\frac{E_{2r+1}}{E_{2r-1}}
 \qquad\text{or}\qquad c>n(n-1).
\]
\end{theorem}

\subsection{The zeros of a linear combination}

The polynomial \(A_n-cxA_{n-2}\) is symmetric in degree \(n-1\)
and has constant term one, so the reciprocal of every zero is a zero.
For \(x\ne0,1\), replacing \(x\) by \(1/x\) changes the sign of
\(z=(1+x)/(1-x)\). Thus \(y=z^2\) takes the same value at both
members of a reciprocal pair. We transform both \(A_n\) and \(xA_{n-2}\)
in degree \(n-1\). By~\eqref{AC},
\[
 \begin{aligned}
 (z+1)^{n-1}A_n(x)&=C_n(z),\\
 (z+1)^{n-1}xA_{n-2}(x)&=(z^2-1)C_{n-2}(z).
 \end{aligned}
\]
By~\eqref{L2C2},
\[
 \frac{d}{dz}\bigl((z^2-1)C_{n-2}(z)\bigr)
 =(z^2-1)C'_{n-2}(z)+2zC_{n-2}(z)=C_{n-1}(z).
\]
The transforms \(C_n\) and \((z^2-1)C_{n-2}\) are even for odd \(n\)
and odd for even \(n\). For even \(n\), divide both by \(z\), removing
the zero that corresponds to \(x=-1\). In the variable \(y=z^2\),
write \(g(y)\) for the resulting transform of \(A_n\) and \(f(y)\)
for that of \(xA_{n-2}\). The zeros of \(g\) and
\(f\) are the squares of the positive zeros of the respective
transforms. The proof of Corollary~\ref{cor:typeA-two-step} shows
that \(g\) strictly alternates left of \(f\), whose largest zero is \(1\).

\begin{lemma}\label{lem:remaining}
Let \(f,g\in\mathbb R[y]\) have degree \(m\ge1\), positive
leading coefficients, and only simple positive zeros. Suppose that
\(g\) strictly alternates left of \(f\), whose zeros are
\(1=a_1>\cdots>a_m>0\). Then $0<g(0)/f(0)<\lc(g)/\lc(f)$.
For every real \(c\), the polynomial
\(g-cf\) has one simple zero in each \((a_{i+1},a_i)\),
\(1\le i<m\). The table locates the remaining zero, which is simple whenever
it exists.
\[
\begin{array}{c|c}
 \text{Value of }c & \text{Remaining zero of }g-cf\\
 \hline
 c<g(0)/f(0) &(0,a_m)\\[1mm]
 c=g(0)/f(0) &0\\[1mm]
 g(0)/f(0)<c<\lc(g)/\lc(f) &(-\infty,0)\\[1mm]
 c=\lc(g)/\lc(f) &\text{none }(\deg(g-cf)=m-1)\\[1mm]
 c>\lc(g)/\lc(f) &(1,\infty)
\end{array}
\]
\end{lemma}
\begin{proof}
Write \(R(y)=g(y)/f(y)\) and \(R(\infty)=\lc(g)/\lc(f)\).
For real \(y\) with \(f(y)\ne0\), Lemma~\ref{lem:residue-representation} gives
\[
 R(y)=R(\infty)+\sum_{i=1}^m\frac{\rho_i}{y-a_i},
 \qquad \rho_i=\frac{g(a_i)}{f'(a_i)}>0,
 \qquad R'(y)=-\sum_{i=1}^m\frac{\rho_i}{(y-a_i)^2}<0.
\]
Both \(f(0)\) and \(g(0)\) have sign \((-1)^m\), since each
polynomial has positive leading coefficient and \(m\) positive zeros.
Thus $R(0)>0$, and the expansion gives
$R(0)=R(\infty)-\sum_{i=1}^m\rho_i/a_i<R(\infty)$.
At each pole \(a_i\), we have \((g-cf)(a_i)=g(a_i)\ne0\).
Elsewhere, the zeros of \(g-cf=f(R-c)\) are precisely the solutions
of \(R(y)=c\). At any such real \(y\) we have $(g-cf)'(y)=f(y)R'(y)\ne0$,
so every real solution is a simple zero of \(g-cf\).

As \(y\) approaches \(a_i\), the term \(\rho_i/(y-a_i)\) tends to
\(-\infty\) from the left and to \(+\infty\) from the right. All other
terms remain bounded, so $R$ has the same one-sided limits at $a_i$.
As $y\to\pm\infty$, each fraction tends to zero and $R(y)\to R(\infty)$.
On each \((a_{i+1},a_i)\), the function \(R\) decreases from
\(+\infty\) to \(-\infty\), giving exactly one solution.
On \((-\infty,a_m)\), \(R\) decreases from \(R(\infty)\) to
\(-\infty\) and takes the value \(R(0)\) at zero. On \((1,\infty)\),
\(R\) decreases from \(+\infty\) to \(R(\infty)\). For \(c<R(\infty)\),
the remaining zero lies to the left of \(a_m\). It is positive at
\(c<R(0)\), zero at \(c=R(0)\), and negative at \(R(0)<c<R(\infty)\).
For \(c>R(\infty)\), the remaining zero lies to the right of \(1\).
For \(c\ne R(\infty)\), the leading coefficient
of \(g-cf\) is \(\lc(f)(R(\infty)-c)\ne0\). The \(m-1\) zeros in the
bounded gaps and the remaining zero exhaust its degree.

When \(c=R(\infty)\), the partial fractions give
\[
 g(y)-cf(y)=\sum_{i=1}^m\rho_i\frac{f(y)}{y-a_i}.
\]
Each quotient \(f(y)/(y-a_i)\) has degree \(m-1\) and leading
coefficient \(\lc(f)\). Thus the coefficient of \(y^{m-1}\)
in the sum is \(\lc(f)\sum_i\rho_i>0\), so the degree is exactly \(m-1\).
The \(m-1\) zeros in the bounded gaps are all the zeros of \(g-cf\).
For \(m=1\), the polynomial is a nonzero constant.
\end{proof}

When \(c=\lc(g)/\lc(f)\), the degree of \(g-cf\) drops by one in
\(y=z^2\), so the Cayley transform of \(A_n-cxA_{n-2}\) drops by
two degrees in \(z\). The next lemma relates the degree drop to a
double zero at \(x=1\).

\begin{lemma}\label{lem:cayley-mult}
Let \(d\ge0\) be an integer and let \(p\) be a nonzero polynomial
of degree at most \(d\). Write \(q\) for its order of vanishing at
\(x=1\), with \(q=0\) if \(p(1)\ne0\). The expression
\[
 (z+1)^d p\!\left(\frac{z-1}{z+1}\right)
\]
extends to a polynomial with values \(2^dp(0)\) at \(z=1\) and
\((-2)^d[x^d]p(x)\) at \(z=-1\). For \(z_0\ne-1\), its order of vanishing at \(z_0\)
equals that of \(p\) at \((z_0-1)/(z_0+1)\). In particular,
\(z=0\) corresponds to \(x=-1\) with the same multiplicity.
The transformed polynomial has degree \(d-q\).
\end{lemma}
\begin{proof}
Expanding \(p\) in powers of \(x\) gives the polynomial extension
\[
 P(z)=\sum_{j=0}^d\bigl([x^j]p(x)\bigr)(z-1)^j(z+1)^{d-j}.
\]
At \(z=1\) only the term with \(j=0\) can contribute, and at
\(z=-1\) only the term with \(j=d\) can contribute. For \(z\ne-1\),
the factor \((z+1)^d\) is nonzero and \(dx/dz=2/(z+1)^2\ne0\),
so the substitution preserves multiplicities. By the definition of
\(q\), we can write
\(p(x)=(x-1)^q h(x)\) with \(h(1)\ne0\). Thus
\[
 \frac{P(z)}{(z+1)^{d-q}}
 =(-2)^qh\!\left(\frac{z-1}{z+1}\right)
 \longrightarrow (-2)^qh(1)\ne0
 \qquad(z\to\infty).
\]
Hence \(\deg P=d-q\).
\end{proof}

\subsection{The two endpoint values}

To apply Lemma~\ref{lem:remaining} to \(g\) and \(f\), we first show
by induction on \(n\) that the positive zeros of \(C_n(z)\) and
\((z^2-1)C_{n-2}(z)\) alternate strictly. We then return to the zeros
of \(A_n\) and \(A_{n-2}\) on each side of \(-1\).

\begin{corollary}\label{cor:typeA-two-step}
Let \(n\ge3\). On each of \((-\infty,-1)\) and \((-1,0)\), the zeros
of \(A_{n-2}(x)\) strictly interlace those of \(A_n(x)\).
More precisely, the zeros of \(A_{n-2}(x)\) divide each interval
into open subintervals, each containing exactly one zero of \(A_n(x)\).
When \(n\) is even, \(-1\) is their unique common zero and is simple
in both polynomials. When \(n\) is odd, neither polynomial vanishes
at \(-1\), and they have no common zero.
\end{corollary}
\begin{proof}
We show that \((z^2-1)C_{n-2}(z)\) and \(C_n(z)\) have only simple
zeros in \([-1,1]\) and \((-1,1)\), respectively. Their positive zeros
alternate strictly and the largest is the zero \(1\) of \((z^2-1)C_{n-2}(z)\).

For \(j\ge1\), the definition of \(C_j\) gives
\[
 \lim_{z\to\infty}\frac{C_j(z)}{(z+1)^{j-1}}
 =\lim_{z\to\infty}A_j\!\left(\frac{z-1}{z+1}\right)=A_j(1)=j!,
\]
since \(A_j\) counts all permutations of \([j]\) at \(x=1\).
Thus \(C_j\) has degree \(j-1\) and leading coefficient \(j!\). In particular,
\begin{equation}\label{eq:FGdegree}
 \begin{aligned}
 \deg\bigl((z^2-1)C_{n-2}\bigr)&=\deg C_n=n-1,\\
 \lc\bigl((z^2-1)C_{n-2}\bigr)&=(n-2)!,\qquad \lc(C_n)=n!.
 \end{aligned}
\end{equation}
For \(n=3\), the pair is \(z^2-1\), \(6z^2-2\), whose positive
zeros satisfy \(1>1/\sqrt3>0\). By \eqref{L2C2} and the product rule,
\(\Lop_0((z^2-1)C_{n-2})=(z^2-1)C_{n-1}\) and \(\Lop_2C_n=C_{n+1}\).
The polynomials \((z^2-1)C_{n-2}\) and \(C_n\) have the same parity,
and \((z^2-1)C_{n-2}\) vanishes at \(z=-1\) and \(z=1\).
Assume the asserted zero locations and strict alternation hold at index
\(n\), with all zeros simple. Theorem~\ref{thm:param-interlace}
with \((\mu,\nu)=(0,2)\) gives the same conclusions at index \(n+1\).
In particular, \(1\) is the largest zero of \((z^2-1)C_{n-1}\).

The positive zeros of \((z^2-1)C_{n-2}\) and \(C_n\) alternate
strictly, with \(1\) last. Removing \(1\) leaves the zeros of
\(C_{n-2}\) and \(C_n\) in \((0,1)\). The resulting list begins
and ends with a zero of \(C_n\), so every subinterval cut out by
the zeros of \(C_{n-2}\) contains exactly one zero of \(C_n\).
Reflection about zero reverses the list and gives the same
interlacing on \((-1,0)\).
The map \(x=(z-1)/(z+1)\) is strictly increasing and sends
\((-1,0)\) and \((0,1)\) onto \((-\infty,-1)\) and \((-1,0)\),
respectively. The factors \((z+1)^{j-1}\) in the definition of \(C_j\) are
nonzero on both intervals. Thus the substitution preserves the
order of the zeros of \(A_{n-2}\) and \(A_n\).
For \(n=3,4\), each of \((-\infty,-1)\) and \((-1,0)\) contains one
zero of \(A_n\) and none of \(A_{n-2}\).
For even \(n\), both \(C_n\) and \(C_{n-2}\) are odd with a simple
zero at \(0\). For odd \(n\), both are even and nonzero at \(0\).
Since \(z=0\) corresponds to \(x=-1\), the two Eulerian polynomials
share a zero there exactly when \(n\) is even.
\end{proof}

\begin{proof}[Proof of Theorem~\ref{thm:Zhang-main}]
Write \(n-1=2m+\varepsilon\), where \(\varepsilon=0\) or \(1\).
Put \(y=z^2\) and define the polynomials \(f\) and \(g\) by
\((z^2-1)C_{n-2}(z)=z^\varepsilon f(y)\) and \(C_n(z)=z^\varepsilon g(y)\).
By~\eqref{AC},
\[
 (z+1)^{n-1}\bigl(A_n(x)-cxA_{n-2}(x)\bigr)
 =z^\varepsilon\bigl(g(y)-cf(y)\bigr),
 \qquad x=\frac{z-1}{z+1}.
\]
The zeros of \(f\) and \(g\) are the squares of the positive zeros of
\((z^2-1)C_{n-2}\) and \(C_n\), respectively. By the proof of
Corollary~\ref{cor:typeA-two-step}, these positive zeros alternate strictly.
Thus \(f\) and \(g\) have degree \(m\) and only simple positive zeros.
Their leading coefficients are positive by~\eqref{eq:FGdegree}.
The largest zero of \(f\) is \(1\), whereas every zero of \(g\)
is less than \(1\). Thus \(g\) strictly alternates left of \(f\),
and Lemma~\ref{lem:remaining} applies to the quotient
\begin{equation}\label{rhotau}
 R(y)=\frac{g(y)}{f(y)}
 =\frac{C_n(z)}{(z^2-1)C_{n-2}(z)}
 =\frac{A_n(x)}{xA_{n-2}(x)},\qquad
 y=\left(\frac{1+x}{1-x}\right)^2.
\end{equation}
Strict alternation gives \(f(y)\ne0\) at every zero of \(g-cf\).
Thus the zeros of \(g-cf\) are exactly the solutions of \(R(y)=c\).
For even \(n\), canceling \(z\) in the quotient of the Cayley transforms
corresponds to canceling \(1+x\) in \(A_n(x)/(xA_{n-2}(x))\). At \(y=0\),
which corresponds to \(z=0\) and \(x=-1\), we use the reduced
quotient \(R(0)=g(0)/f(0)\), with \(f(0)\ne0\) in both parities.

The poles of \(R\) lie in \((0,1]\). Between any two consecutive poles,
Lemma~\ref{lem:remaining} gives a solution of \(R(y)=c\) for every \(c\).
Each solution in a gap satisfies \(0<y<1\), so \(z=\pm\sqrt y\)
gives two negative zeros of \(A_n-cxA_{n-2}\). At \(c=R(0)\), the
remaining solution of \(R(y)=c\) is \(y=0\).
At \(c=R(\infty)\), the degree of \(g-cf\) drops by one.
The leading coefficients in \eqref{eq:FGdegree} give
\(R(\infty)=n!/(n-2)!=n(n-1)\). If \(n=2r+1\), then
\(R(0)=-C_{2r+1}(0)/C_{2r-1}(0)\). If \(n=2r\), both
\(C_{2r}\) and \(C_{2r-2}\) have simple zeros at \(0\).
Canceling their factors \(z\) gives
\[
 R(0)=-\frac{C'_{2r}(0)}{C'_{2r-2}(0)}
      =-\frac{C_{2r+1}(0)}{C_{2r-1}(0)}.
\]
For the second equality, set \(z=0\) in~\eqref{L2C2} to get
\(C_{j+1}(0)=-C'_j(0)\). Since \(C_j(0)=A_j(-1)\), in both parities
\(R(0)=-A_{2r+1}(-1)/A_{2r-1}(-1)\).
The Eulerian generating function~\cite[Theorem~1.6]{Petersen2015} is
\begin{equation}\label{eq:eulerian-egf-evaluation}
 \sum_{j\ge0}A_j(x)\frac{w^j}{j!}
 =\frac{1-x}{e^{(x-1)w}-x}.
\end{equation}
At \(x=-1\) the generating function reduces to $2/(1+e^{-2w})=1+\tanh w$.
The odd part of the generating function
\(\sum_{j\ge0}E_jw^j/j!=\sec w+\tan w\) is \(\tan w\).
Substituting \(\tanh w=-i\tan(iw)\) and comparing the coefficients
of \(w^{2j+1}\) gives
\(A_{2j+1}(-1)=(-1)^jE_{2j+1}\)
\cite[Problem~4.2(3) and its solution]{Petersen2015}.
Consequently,
\begin{equation}\label{eq:Zhang-endpoints}
 R(0)=\frac{E_{2r+1}}{E_{2r-1}},\qquad
 R(\infty)=n(n-1),\qquad 0<R(0)<R(\infty).
\end{equation}

Let \(y_0\) be a zero of \(g-cf\). The factor \(y-y_0\) gives a
quadratic factor of \(A_n-cxA_{n-2}\) because \(y=z^2\) and
\begin{equation}\label{eq:cayley-quadratic-factor}
 \frac{z^2-y_0}{(z+1)^2}
 =\frac{(1+x)^2-y_0(1-x)^2}{4}.
\end{equation}
Indeed, \(z/(z+1)=(1+x)/2\) and \(1/(z+1)=(1-x)/2\).
The factor \(z^\varepsilon\) similarly contributes
\(((1+x)/2)^\varepsilon\). When \(\deg(g-cf)=m\), we factor
\(g-cf\) into linear factors and apply~\eqref{eq:cayley-quadratic-factor}
to each factor.
This gives
\[
 A_n(x)-cxA_{n-2}(x)
 =\lc(g-cf)\left(\frac{1+x}{2}\right)^{\!\varepsilon}
  \prod_{(g-cf)(y_0)=0}\frac{(1+x)^2-y_0(1-x)^2}{4},
\]
since \(2m+\varepsilon=n-1\).

The numerator on the right of~\eqref{eq:cayley-quadratic-factor}
is $(1-y_0)x^2+2(1+y_0)x+(1-y_0)$.
The value \(y_0=1\) cannot occur, since \(1\) is a pole of \(R\).
The two zeros of this quadratic have product one and sum
\(-2(1+y_0)/(1-y_0)\). Its discriminant is
\(4(1+y_0)^2-4(1-y_0)^2=16y_0\).
For \(0<y_0<1\) the zeros are real and distinct and their sum is
less than \(-2\). Both are negative, one on either side of \(-1\).
For \(y_0>1\) the sum exceeds \(2\), giving two positive reciprocal
zeros, one on either side of \(1\). For \(y_0<0\) the discriminant
is negative. The two zeros are conjugates with product one, so
both have modulus one. For each finite \(y_0\), the numerator
\((1+x)^2-y_0(1-x)^2\) equals \(4\) at \(x=1\). Hence each of its zeros
\(x\) differs from \(1\) and
determines \(y_0=((1+x)/(1-x))^2\). Thus quadratics for distinct
values of \(y_0\) have no common zero. At \(c=n(n-1)\), the zero
\(x=1\) comes from the drop in degree of \(g-cf\).

For every \(c\), Lemma~\ref{lem:remaining} gives one simple
zero of \(g-cf\) between each pair of consecutive poles of \(R\).
These \(m-1\) zeros all lie in \((0,1)\).
Their quadratic factors give \(2(m-1)\) distinct negative zeros,
none equal to \(-1\). The product of these quadratics and
\((1+x)^\varepsilon\) therefore has \(2(m-1)+\varepsilon=n-3\)
simple negative zeros.

The polynomial \(A_n-cxA_{n-2}\) is monic of degree \(n-1\), since
\(cxA_{n-2}\) has degree at most \(n-2\). For \(c\ne R(\infty)\),
its other two zeros come from the quadratic factor corresponding
to the remaining zero \(y_0\) of \(g-cf\). For \(c<R(0)\), this
zero \(y_0\) lies in \((0,1)\) and gives two distinct negative
zeros in \(x\). For \(R(0)<c<R(\infty)\),
it lies in \((-\infty,0)\) and gives a nonreal conjugate pair on the
unit circle. For \(c>R(\infty)\), it lies in \((1,\infty)\) and
gives two distinct positive zeros.

At \(c=R(0)\), the remaining zero \(y_0=0\) gives the factor
\((1+x)^2/4\) and is simple because $(g-cf)'(0)=f(0)R'(0)\ne0$.
For every other zero \(y_0\) of \(g-cf\), the factor
$\bigl((1+x)^2-y_0(1-x)^2\bigr)/4$ has value $-y_0\ne0$ at $x=-1$.
Thus \(-1\) has multiplicity exactly \(2+\varepsilon\), namely two
for odd \(n\) and three for even \(n\). At all other values of \(c\),
\((g-cf)(0)=f(0)(R(0)-c)\ne0\), so the multiplicity at \(-1\)
is exactly \(\varepsilon\).

As \(c\uparrow R(\infty)\), the remaining solution \(y_0\) in the
leftmost interval tends to \(-\infty\). As \(c\downarrow R(\infty)\),
the solution in the rightmost interval tends to \(+\infty\).
In either limit, dividing \((1+x)^2-y_0(1-x)^2\) by \(-y_0\) gives
the limiting quadratic \((1-x)^2\), so both zeros tend to \(1\).
At \(c=R(\infty)\) itself, the degree of \(g-cf\) is \(m-1\).
To return to \(x\), divide each of the \(m-1\) factors \(z^2-y_0\)
by \((z+1)^2\), and divide \(z^\varepsilon\) by \((z+1)^\varepsilon\).
Since \(2(m-1)+\varepsilon=n-3\), the two remaining powers of \(z+1\)
in the denominator give the factor \(1/(z+1)^2=(1-x)^2/4\).
For each remaining quadratic factor,
\(\left.\tfrac14\bigl((1+x)^2-y_0(1-x)^2\bigr)\right|_{x=1}=1\).
Since \((1+x)/2\) also equals one at \(x=1\), we obtain
\[
 \lim_{x\to1}\frac{A_n(x)-cxA_{n-2}(x)}{(x-1)^2}
 =\tfrac14\lc(g-cf)
 =\tfrac14\lc\bigl(C_n-c(z^2-1)C_{n-2}\bigr)\ne0.
\]
Thus \(1\) is a double zero, as the degree formula in
Lemma~\ref{lem:cayley-mult} also shows. At either endpoint all
zeros other than the multiple zero are simple and negative.
\end{proof}

The proof also gives the location of the zeros for every real \(c\).

\begin{corollary}\label{cor:zhang-boundaries}
Let \(n\ge3\), \(r=\lfloor n/2\rfloor\), and \(c\in\mathbb R\).
The polynomial \(A_n(x)-cxA_{n-2}(x)\) has only real zeros if and only if
\[
 c\le \frac{E_{2r+1}}{E_{2r-1}}\qquad\text{or}\qquad c\ge n(n-1).
\]
Below the smaller endpoint, all \(n-1\) zeros are simple and negative.
Above the larger endpoint, all are simple, with two positive and
\(n-3\) negative. Strictly between the endpoints there are \(n-3\)
simple negative zeros and one nonreal conjugate pair on the unit circle.

At \(c=E_{2r+1}/E_{2r-1}\), the zero \(-1\) is double for odd \(n\)
and triple for even \(n\). At \(c=n(n-1)\), the zero \(1\) is double.
All other zeros at either endpoint are simple and negative.
\end{corollary}

For example, \(A_4(x)-cxA_2(x)=(1+x)\bigl(x^2+(10-c)x+1\bigr)\).
The quadratic has discriminant \((c-8)(c-12)\), and its zeros have
product one. At \(c=8,10,12\), respectively, we obtain
\((1+x)^3\), \((1+x)(x^2+1)\), and \((1+x)(x-1)^2\).
At \(c=8\) the quadratic contributes a double zero at \(-1\), and
the factor \(1+x\) raises its multiplicity to three. For \(8<c<12\),
the two nonreal zeros have real part \((c-10)/2\). Thus they cross
the imaginary axis at \(c=10=2E_4/E_2\), the Hurwitz stability
boundary for \(A_4-cxA_2\) in Section~\ref{sec:hurwitz}.

We now ask for which \(c\) the polynomial \(A_n-cxA_{n-2}\) has
only simple negative zeros for all \(n\ge3\).
The lower endpoint \(E_{2r+1}/E_{2r-1}\) is at least \(2\) for
\(r\ge1\), with \(E_3/E_1=2\). Hence \(c<2\) is sufficient.
This recovers the assertion on simple real zeros in Br\"and\'en's
Theorem~5.2~\cite{Branden2006}. In our normalization his polynomial is
\(A_n(t,x)=x\bigl(A_n(x)+txA_{n-2}(x)\bigr)\), and \(t>-2\) corresponds to \(c<2\).
Conversely, for \(n=3\) the polynomial $A_3(x)-cxA_1(x)=x^2+(4-c)x+1$
has discriminant $(4-c)^2-4=(c-2)(c-6)$.
Its zeros are distinct and real exactly when \(c<2\) or \(c>6\).
Their product is one and their sum is \(c-4\), so they are distinct
and negative exactly when \(c<2\). At \(c=2\) the zero \(-1\) is double.
Thus \(A_n(x)-cxA_{n-2}(x)\) has only simple negative zeros for
every \(n\ge3\) if and only if \(c<2\).
In particular, the descent polynomial
\(A_n-xA_{n-2}\) has \(n-1\) distinct negative zeros.

\section{Descent polynomials of \texorpdfstring{\((n-3)\)}{(n-3)}-stack-sortable permutations}
\label{sec:stacksorting}

\label{sec:bona}

The stack-sorting map \(s\) acts on words with distinct letters and
fixes the empty word. Write a nonempty word as \(uav\), where \(a\)
is its largest letter and \(u\) and \(v\) are the subwords before and
after it. Set \(s(uav)=s(u)s(v)a\).
A permutation \(\pi\in\mathfrak S_n\)
is \(t\)-stack-sortable if \(s^t(\pi)=12\cdots n\). Put
\[
 W_{n,t}(x)=
 \sum_{\substack{\pi\in\mathfrak S_n\\s^t(\pi)=12\cdots n}}
 x^{\des(\pi)}.
\]
For example, \(s(231)=213\) and \(s(213)=123\). The permutation
\(231\) has one descent, so it contributes \(x\) to \(W_{3,2}\).
It does not contribute to \(W_{3,1}\).
B\'ona's conjecture~\cite{Bona2002} states that \(W_{n,t}\) has only
real zeros for every \(1\le t\le n-1\). The cases \(t=2\) and
\(t=n-2\) follow from Br\"and\'en's Theorem~5.1 and
Corollary~5.3~\cite{Branden2006}. For \(t=n-2\), we have
\(W_{n,n-2}=A_n-xA_{n-2}\). Zhang and
Kitaev~\cite[Theorem~1]{ZhangKitaev2025} gave the following formula for
\(W_{n,n-3}\) and left its real-rootedness open.
We show that all its zeros are simple and negative.

\begin{theorem}\label{thm:stack-n-3}
For every \(n\ge4\), the descent polynomial of the
\((n-3)\)-stack-sortable permutations,
\begin{equation}\label{eq:KZ-stack-formula}
 W_{n,n-3}(x)=A_n(x)-\frac52xA_{n-2}(x)
 -\left(n-\frac32\right)x(1+x)A_{n-3}(x),
\end{equation}
has \(n-1\) distinct negative zeros.
\end{theorem}

The operator proof compares the zeros of \(xA_{n-3}\) and
\(A_{n-1}-(n-\tfrac32)xA_{n-3}\). The latter is the combination
in Theorem~\ref{thm:Zhang-main} at index \(n-1\) with \(c=n-\tfrac32\).
The next lemma gives the ordering needed to apply
Theorem~\ref{thm:param-interlace}.

\begin{lemma}\label{lem:stack-auxiliary-zeros}
Let \(n\ge5\). The polynomial
\(A_{n-1}(x)-\left(n-\tfrac32\right)xA_{n-3}(x)\)
has \(n-2\) simple negative zeros. Its zeros below \(-1\) are the
reciprocals of its zeros in \((-1,0)\). The zeros of \(A_{n-3}\) in
\((-1,0)\) divide this interval into subintervals, each containing
exactly one zero of \(A_{n-1}-(n-\tfrac32)xA_{n-3}\).
For odd \(n\), the only common zero of \(A_{n-1}-(n-\tfrac32)xA_{n-3}\)
and \(A_{n-3}\) is \(-1\). For even \(n\), they have no common zero.
\end{lemma}

\begin{proof}
Put \(r=\lfloor(n-1)/2\rfloor\). At index \(n-1\),
equations~\eqref{rhotau} and~\eqref{eq:Zhang-endpoints} give
\[
 y=\left(\frac{1+x}{1-x}\right)^2,\qquad
 R(y)=\frac{A_{n-1}(x)}{xA_{n-3}(x)},\qquad
 R(0)=\frac{E_{2r+1}}{E_{2r-1}}.
\]
For odd \(n\), cancel the common factor \(1+x\) in
\(A_{n-1}(x)/(xA_{n-3}(x))\) before evaluating at \(x=-1\).
For \(-1<x<0\) with \(A_{n-3}(x)\ne0\), the polynomial in the
lemma vanishes exactly when \(R(y)=n-\tfrac32\).

Since \(r\ge2\) and \(n\le2r+2\), the alternating-permutation
count in the Introduction gives
$R(0)=E_{2r+1}/E_{2r-1}\ge4r>2r+\tfrac12\ge n-\tfrac32$.

The map \(x\mapsto y\) is strictly increasing from \((-1,0]\) onto
\((0,1]\). By Corollary~\ref{cor:typeA-two-step}, each zero of
\(A_{n-3}\) in \((-1,0)\) gives a pole of \(R\) in \((0,1)\),
since it is not a zero of \(A_{n-1}\).
The factor \(x\) gives the pole \(y=1\), since \(A_{n-1}(0)=1\).
The zeros of \(A_{n-3}\) other than \(-1\) occur in reciprocal
pairs, with one zero of each pair in \((-1,0)\). Since
\(\deg A_{n-3}=n-4\), this interval contains
\(\lfloor(n-4)/2\rfloor\) zeros of \(A_{n-3}\).
In the variable \(z\), the denominator before cancellation is
\((z^2-1)C_{n-3}(z)\), of degree \(n-2\). For odd \(n\), cancel \(z\)
from both transforms.
Setting \(y=z^2\) then gives a denominator of degree
\(\lfloor(n-2)/2\rfloor\). The
\(1+\lfloor(n-4)/2\rfloor=\lfloor(n-2)/2\rfloor\) distinct poles
just found exhaust that degree.

The proof of Corollary~\ref{cor:typeA-two-step} at index \(n-1\)
shows that the reduced numerator of \(R\) strictly alternates left of
its denominator, whose largest zero is \(1\). Since \(n-\tfrac32<R(0)\),
Lemma~\ref{lem:remaining} gives one simple solution of
\(R(y)=n-\tfrac32\) between \(0\) and the smallest pole and one
between each pair of consecutive poles. The total is
\(\lfloor(n-2)/2\rfloor\) solutions in \((0,1)\).
For each such \(y\), taking \((1+x)/(1-x)=\pm\sqrt y\) gives
\[
 \frac{\sqrt y-1}{\sqrt y+1}\in(-1,0),\qquad
 \frac{\sqrt y+1}{\sqrt y-1}\in(-\infty,-1).
\]
Neither is a zero of \(xA_{n-3}\). Distinct solutions in \(y\) give
disjoint pairs.
Since \(y\) increases with \(x\) on \((-1,0]\), the solutions between
\(0\) and the smallest pole and between successive poles give exactly one
zero in each of the asserted subintervals of \((-1,0)\).

For odd \(n\), write
\[
 A_{n-1}(x)-(n-\tfrac32)xA_{n-3}(x)
 =xA_{n-3}(x)\bigl(R(y)-(n-\tfrac32)\bigr).
\]
At \(x=-1\), \(xA_{n-3}(x)\) has a simple zero and
\(R(y)-(n-\tfrac32)\) has the nonzero limit \(R(0)-(n-\tfrac32)>0\).
Thus \(A_{n-1}-(n-\tfrac32)xA_{n-3}\) has a simple zero at \(-1\).
For even \(n\), the reciprocal pairs give \(n-2\) distinct zeros.
For odd \(n\), counting \(-1\) as well gives the same total.
The combination is monic of degree \(n-2\), so all \(n-2\) zeros are simple.
\end{proof}

\begin{proof}[Proof of Theorem~\ref{thm:stack-n-3}]
At \(n=4\), formula~\eqref{eq:KZ-stack-formula} gives
$W_{4,1}(x)=1+6x+6x^2+x^3=(x+1)(x^2+5x+1)$.
The zeros \(-1\) and \((-5\pm\sqrt{21})/2\) are distinct
and negative.

Now suppose \(n\ge5\). Put \(x=(z-1)/(z+1)\) and define
$\widehat W_n(z)=(z+1)^{n-1}W_{n,n-3}(x)$.
The identities \((z+1)^2x=z^2-1\) and
\((z+1)^3x(1+x)=2z(z^2-1)\) give
\[
 \widehat W_n
 =C_n-\frac52(z^2-1)C_{n-2}
       -(2n-3)z(z^2-1)C_{n-3}.
\]
Set \(F=(z^2-1)C_{n-3}\) and \(J=C_{n-1}-(n-\tfrac32)F\).
By~\eqref{L2C2}, \(F'=C_{n-2}\) and \(C_n=\Lop_2C_{n-1}\).
Using \(2zF=\Lop_2F-\Lop_0F\), we obtain
\begin{equation}\label{eq:stack-interlacing-decomp}
 \begin{aligned}
 \widehat W_n
 &=\Lop_2 C_{n-1}-\frac52\Lop_0F
   -\left(n-\frac32\right)(\Lop_2F-\Lop_0F)\\
 &=\Lop_2J+(n-4)\Lop_0F.
 \end{aligned}
\end{equation}
Both \(F\) and \(J\) are even when \(n\) is even and odd when \(n\) is odd.
With \(x=(z-1)/(z+1)\), the polynomials \(F\) and \(J\) are given by
\[
 F(z)=(z+1)^{n-2}xA_{n-3}(x),\qquad
 J(z)=(z+1)^{n-2}\bigl(A_{n-1}(x)-(n-\tfrac32)xA_{n-3}(x)\bigr).
\]
The strictly increasing inverse map \(z=(1+x)/(1-x)\) sends
\((-\infty,0)\) onto \((-1,1)\) and \((-1,0]\) onto \((0,1]\).
Thus the \(n-2\) simple negative zeros in
Lemma~\ref{lem:stack-auxiliary-zeros} give \(n-2\) simple zeros of \(J\)
in \((-1,1)\).
The positive zeros of \(F\) are the images of the zeros of
\(A_{n-3}\) in \((-1,0)\), together with \(1\), which comes from
the factor \(x\). By Lemma~\ref{lem:stack-auxiliary-zeros}, \(J\) has
one zero between \(0\) and the smallest positive zero of \(F\), and
one between each pair of consecutive positive zeros of \(F\).
Thus the positive zeros of \(F\) and \(J\) alternate strictly.
The proof of Corollary~\ref{cor:typeA-two-step} at index \(n-1\)
shows that \(F=(z^2-1)C_{n-3}\) has \(n-2\) simple zeros in \([-1,1]\),
including both endpoints. By~\eqref{eq:FGdegree},
the leading coefficients of \(F\) and \(J\) are
\(\lc(F)=(n-3)!\) and \(\lc(J)=(n-3)!((n-1)(n-2)-(n-\tfrac32))>0\).
Apply Theorem~\ref{thm:param-interlace} to \(F,J\) with
\((\mu,\nu)=(0,2)\). The polynomials \(\Lop_0F\) and \(\Lop_2J\)
have degree \(n-1\) and only simple zeros in \([-1,1]\).
Every zero of \(\Lop_2J\) lies in \((-1,1)\).
The leading coefficients are
\((n-2)\lc(F)>0\) for \(\Lop_0F\) and \(n\lc(J)>0\) for \(\Lop_2J\).
Since \(n-4>0\), \eqref{eq:stack-interlacing-decomp} gives
\(\deg\widehat W_n=n-1\).

Put \(m=\lfloor(n-1)/2\rfloor\). By parity each image has \(m\)
positive zeros, and the two sets of positive zeros alternate strictly.
The largest of these zeros is \(1\), which belongs to \(\Lop_0F\).
In increasing order the positive zeros are
\[
 \beta_1<\alpha_1<\cdots<\beta_m<\alpha_m=1,
 \qquad (\Lop_2J)(\beta_j)=(\Lop_0F)(\alpha_j)=0.
\]
There are \(m-j+1\) zeros of \(\Lop_0F\) to the right of \(\beta_j\)
and \(m-j\) zeros of \(\Lop_2J\) to the right of \(\alpha_j\).
Both leading coefficients are positive, so
\((\Lop_0F)(\beta_j)\) and \((\Lop_2J)(\alpha_j)\) have opposite signs.
By~\eqref{eq:stack-interlacing-decomp},
\[
 \widehat W_n(\beta_j)=(n-4)(\Lop_0F)(\beta_j),\qquad
 \widehat W_n(\alpha_j)=(\Lop_2J)(\alpha_j).
\]
Since \(n-4>0\), \(\widehat W_n\) has a zero in each of the \(m\)
disjoint intervals \((\beta_j,\alpha_j)\subset(0,1)\).
The polynomial \(\widehat W_n\) is even for odd \(n\) and odd for
even \(n\). Reflection gives \(m\) distinct negative zeros.
If \(n\) is even, \(0\) is also a zero. The \(n-1\) distinct zeros
in \((-1,1)\) exhaust the degree of \(\widehat W_n\) and are simple.

The map \(x=(z-1)/(z+1)\) is strictly increasing on \((-1,1)\)
and has image \((-\infty,0)\). Since \((z+1)^{n-1}\ne0\) for
\(-1<z<1\), the \(n-1\) zeros of \(\widehat W_n\) give \(n-1\)
distinct negative zeros of \(W_{n,n-3}\). Formula~\eqref{eq:KZ-stack-formula}
shows that \(W_{n,n-3}\) has leading term \(x^{n-1}\), since both
subtracted terms have degree \(n-2\). The \(n-1\) negative zeros exhaust
the degree of \(W_{n,n-3}\) and are simple.
\end{proof}

At \(n=6\), we have \(2<t=n-3<n-2\) for the first time.
Exactly \(606\) of the
\(720\) permutations of \([6]\) are \(3\)-stack-sortable. Their
descent polynomial is
\(W_{6,3}(x)=(1+x)(1+49x+203x^2+49x^3+x^4)\).
The quartic is positive near \(-\infty\) and takes the values
\(-669\) at \(-5\) and \(107\) at \(-1\). It has a zero in each
of \((-\infty,-5)\) and \((-5,-1)\). Since the quartic is
symmetric in degree four, it has two further zeros in \((-1,0)\).
Together with \(-1\), these give five
distinct negative zeros of \(W_{6,3}\).

\section{The Hurwitz stability conjecture of Yang and Zhang}
\label{sec:YangZhang}

\label{sec:hurwitz}

By Yang and Zhang's Theorem~2.3~\cite{YangZhang2015}, the polynomial
\((1+x)A_n+kxA_{n-1}\) is weakly Hurwitz stable for \(n\ge2\)
and \(k\ge-n\). Their Conjecture~3.2 predicts the exact range
of Hurwitz stability for \(n\ge3\). They also ask how this range
relates to Zhang's conjecture on \(A_n-cxA_{n-2}\)
\cite[p.~472]{YangZhang2015}.
At index \(n-1\), Yang and Zhang's polynomial is
\((1+x)A_{n-1}+kxA_{n-2}\), so the parameters in both conjectures
multiply \(xA_{n-2}\). We obtain a common stability criterion
by taking any nonzero nonnegative linear combination of the other
terms, \(A_n\) and \((1+x)A_{n-1}\).

\begin{theorem}\label{thm:combined-hurwitz}
Let \(n\ge3\), \(a,b\ge0\) with \(a+b>0\), and \(c\in\mathbb R\).
The polynomial
\[
 aA_n(x)+b(1+x)A_{n-1}(x)+cxA_{n-2}(x)
\]
is Hurwitz stable if and only if
\(2aE_n+2bE_{n-1}+cE_{n-2}>0\).
It is weakly Hurwitz stable if and only if
\(2aE_n+2bE_{n-1}+cE_{n-2}\ge0\).
\end{theorem}

Each summand in Theorem~\ref{thm:combined-hurwitz} is symmetric
in degree \(n-1\). The sum has positive constant term \(a+b\),
so every zero \(x\) is nonzero and \(1/x\) is also a zero.
Put \(u=x+x^{-1}\). Since
\(\operatorname{Re}u=(1+|x|^{-2})\operatorname{Re}x\), the real parts
of \(u\) and \(x\) have the same sign. For a quadratic example, take
\(f(x)=x^2+\tau x+1\) with \(\tau\in\mathbb R\) and \(g(x)=1\).
Then \(f(x)/(xg(x))=u+\tau\).
At a zero of \(f(x)+cxg(x)\), we have \(u=-\tau-c\), so
the polynomial is Hurwitz stable exactly when \(\tau+c>0\).

\begin{lemma}\label{lem:reciprocal-pair}
Let \(f\) and \(g\) be real polynomials of degrees \(d\) and \(d-2\),
respectively, where \(d\ge2\), and assume that both are symmetric in
their respective degrees. Suppose \(\lc(f)>0\), and let \(g\)
be monic with only simple negative zeros. Assume \(f(r)/g'(r)>0\)
at every zero \(r<-1\) of \(g\).
The number \(f(i)/(i\,g(i))\) is real. For every real \(c\),
the polynomial \(f(x)+cxg(x)\) is Hurwitz stable if and only if
$c>-f(i)/(i\,g(i))$. It is weakly Hurwitz stable if and only if
\(c\ge -f(i)/(i\,g(i))\).

Write \(d=2m+\varepsilon\), where \(\varepsilon\in\{0,1\}\).
The quotient $f(x)/(xg(x))$ is a rational function of \(u=x+x^{-1}\).
Its finite poles are the distinct numbers
\(r+r^{-1}<-2\), where \(r\) ranges over the zeros of \(g\) below \(-1\).
As a rational function of the real variable \(u\), the quotient increases
strictly from \(-\infty\) to \(+\infty\) on each interval cut out by its finite poles.
If there are no finite poles, the interval is all of \(\mathbb R\).
Thus the quotient takes the value \(-c\) at \(m\) distinct real points
\(t_1<\cdots<t_m\), one in each interval. Each \(t_j\) is a simple
solution in \(u\), and \(t_j<-2\) for \(j<m\). The equations
\(x+x^{-1}=t_j\) give the quadratic factors in
\[
 f(x)+cxg(x)=\lc(f)(1+x)^\varepsilon
             \prod_{j=1}^m(x^2-t_jx+1).
\]
\end{lemma}

\begin{proof}
When \(xg(x)\ne0\), division by \(xg(x)\) gives \(f(x)/(xg(x))=-c\)
at each zero of \(f+cxg\). Since \(f\) and \(xg\) are symmetric in
degree \(d\), their quotient has the same value at \(x\) and \(1/x\).
We write \(Q(u)=f(x)/(xg(x))\) with \(u=x+x^{-1}\) and show that
\(Q\) is rational in \(u\).

If \(h\) is symmetric in odd degree \(d\), then
\(h(-1)=-h(-1)\), so \(1+x\) divides \(h\).
The quotient \(h(x)/(1+x)\) is symmetric in degree \(d-1\), since
\[
 x^{d-1}\frac{h(1/x)}{1+1/x}
 =\frac{x^dh(1/x)}{1+x}=\frac{h(x)}{1+x}.
\]
Now let \(h(x)=\sum_{j=0}^{2m}h_jx^j\) be symmetric in degree \(2m\).
Since \(h_{m-j}=h_{m+j}\), division by \(x^m\) gives
\[
 x^{-m}h(x)=h_m+\sum_{j=1}^m h_{m+j}(x^j+x^{-j}).
\]
Each sum \(x^j+x^{-j}\) is a polynomial in \(u\), by the recurrence
\(x^{j+1}+x^{-j-1} =u(x^j+x^{-j})-(x^{j-1}+x^{-j+1})\).
The initial terms are \(2\) and \(u\), so for \(j\ge1\)
the polynomial representing \(x^j+x^{-j}\) is monic of degree
\(j\). In particular, the coefficient of \(u^m\) in \(x^{-m}h(x)\)
is the leading coefficient of \(h\).
The zeros of \(g\) other than \(-1\) come in reciprocal pairs.
All zeros of \(g\) are negative, so neither \(0\) nor \(1\) is a zero.
Each reciprocal pair contributes two to the degree. Since all zeros
are simple, \(-1\) occurs exactly when \(d-2\) is odd, that is, when
\(\varepsilon=1\). List the numbers \(r+r^{-1}\), for zeros \(r<-1\)
of \(g\), in increasing order as \(s_1<\cdots<s_{m-1}<-2\) when
\(m>1\). For \(s_j=r+r^{-1}\), the corresponding pair gives the factor
$(x-r)(x-r^{-1})=x(u-s_j)$.
The numbers \(s_j\) are distinct because
\(d(r+r^{-1})/dr=1-r^{-2}>0\) for \(r<-1\).
Since \(g\) is monic, it follows that
\[
 g(x)=(1+x)^\varepsilon x^{m-1}\prod_{j=1}^{m-1}(u-s_j).
\]
The quotient \(f(x)/(1+x)^\varepsilon\) is symmetric in degree \(2m\).
When \(d\) is odd, \(f\) and \(g\) both contain \(1+x\).
Applying the even-degree calculation to \(f(x)/(1+x)^\varepsilon\)
shows that \(f(x)/(x^m(1+x)^\varepsilon)\) is a polynomial in \(u\) of degree
\(m\) with leading coefficient \(\lc(f)\). Thus \(f/(xg)\) is a rational function of \(u\), and
\[
 Q(u)\prod_{j=1}^{m-1}(u-s_j)=\frac{f(x)}{x^m(1+x)^\varepsilon}.
\]
Since all \(s_j<-2\), the denominator of \(Q\) in this expression
is nonzero at \(u=-2\).

The polynomial part of \(Q\) is linear with slope \(\lc(f)\).
For each \(s_j=r+r^{-1}\) corresponding to a zero \(r<-1\) of \(g\),
we have $u-s_j=(x-r)(1-1/(xr))$.
Since \(g\) has a simple zero at \(r\),
\[
 \lim_{x\to r}(u-s_j)\frac{f(x)}{xg(x)}
 =\left(1-\frac{1}{r^2}\right)\frac{f(r)}{rg'(r)}<0.
\]
Here \(1-r^{-2}>0\) because \(r<-1\), and
\(f(r)/(rg'(r))<0\) because \(f(r)/g'(r)>0\) and \(r<0\).
Thus each \(s_j\) is a simple pole with negative residue,
and no factor \(u-s_j\) cancels.
It remains to find the constant term of the linear polynomial part.
Symmetry in their respective degrees gives \(f(0)=\lc(f)\) and
\(g(0)=\lc(g)=1\). Let \(x\to0^+\), so that \(u=x+x^{-1}\to+\infty\). Then
\[
 Q(u)-\lc(f)u
 =\frac{f(x)/g(x)-f(0)}{x}-f(0)x
 \longrightarrow
 \left(\frac fg\right)'(0)=f'(0)-f(0)g'(0).
\]
Each partial fraction tends to zero at infinity, so the limit gives
the constant term. Write \(-w_j\) for the residue at \(s_j\), with
\(w_j>0\). For real \(u\) away from the poles we obtain
\begin{equation}\label{eq:reciprocal-quotient}
 \begin{aligned}
 Q(u)&=\lc(f)u+f'(0)-f(0)g'(0)
        -\sum_{j=1}^{m-1}\frac{w_j}{u-s_j},\\
 Q'(u)&=\lc(f)+\sum_{j=1}^{m-1}\frac{w_j}{(u-s_j)^2}>0.
 \end{aligned}
\end{equation}
At each pole \(s_j\), the left and right limits of \(Q\) are
\(+\infty\) and \(-\infty\), respectively. Since
\(Q(u)/u\to\lc(f)>0\) as \(u\to\pm\infty\), we have
\(Q(u)\to-\infty\) as \(u\to-\infty\) and \(Q(u)\to+\infty\)
as \(u\to+\infty\). Hence \(Q(u)=-c\) has exactly one
solution in each of the \(m\) intervals cut out by the poles.
This also covers \(m=1\), when there are no poles and \(Q\) is linear.

Write \(t_1<\cdots<t_m\) for the solutions of \(Q(u)=-c\).
Clearing the denominator gives a polynomial of degree \(m\)
with leading coefficient \(\lc(f)\), since the polynomial part
of \(Q\) is linear with slope \(\lc(f)\). Factoring at the \(m\)
distinct zeros \(t_j\) gives
\[
 (Q(u)+c)\prod_{j=1}^{m-1}(u-s_j)
 =\lc(f)\prod_{j=1}^m(u-t_j).
\]
For \(1\le j<m\), the solution in the \(j\)th interval satisfies
\(t_j<s_j<-2\). For \(xg(x)\ne0\), the identity
\(x(u-t_j)=x^2-t_jx+1\) gives
\[
 \begin{aligned}
 f(x)+cxg(x)
 &=xg(x)(Q(u)+c)\\
 &=(1+x)^\varepsilon x^m(Q(u)+c)\prod_{j=1}^{m-1}(u-s_j)\\
 &=\lc(f)(1+x)^\varepsilon
   \prod_{j=1}^m(x^2-t_jx+1).
 \end{aligned}
\]
The first and last expressions are polynomials in \(x\), so the
factorization holds for every \(x\). Since \(cxg\) has degree at most
\(d-1\), the degree of \(f+cxg\) is \(d=2m+\varepsilon\), as is the
degree of the product. In odd degree the explicit factor \(1+x\)
restores the one copy canceled in forming \(f/(xg)\).
Finally, \((f+cxg)(0)=f(0)=\lc(f)>0\),
so the substitution \(u=x+x^{-1}\) omits no zero at the origin.

For \(j<m\), the quadratic \(x^2-t_jx+1\) has discriminant
\(t_j^2-4>0\). Its zeros have product one and sum \(t_j<-2\),
so both are negative and distinct. Neither equals \(-1\), since
the quadratic has value \(2+t_j<0\) there. Distinct quadratic factors
have no common zero, since each zero \(x\) determines its parameter
as \(x+x^{-1}\).
Only the last solution \(t_m\) can equal \(-2\) or \(2\).
At \(t_m=-2\), the last quadratic is \((x+1)^2\), so \(-1\) has
multiplicity \(2+\varepsilon\), since the other quadratics
are nonzero there. At \(t_m=2\), the last quadratic is \((x-1)^2\).
The other quadratics take values \(2-t_j>0\) at \(1\).
Since \((1+x)^\varepsilon\) is also nonzero at \(1\), the multiplicity
is exactly two.
In both cases all other zeros are simple and negative.

Thus only \(x^2-t_m x+1\) can affect Hurwitz stability.
For either zero \(x\) of \(x^2-t_mx+1\), the constant term gives
\(x\ne0\). Taking real parts in
\(t_m=x+x^{-1}=x+\overline{x}/|x|^2\) gives
$t_m=(1+|x|^{-2})\operatorname{Re}x$.
Since \(1+|x|^{-2}>0\), each zero has real part of the same sign
as \(t_m\). Thus both zeros have negative real part exactly when \(t_m<0\).
The points \(t_m\) and \(0\) lie in the rightmost interval, on which
\(Q\) is strictly increasing. Therefore
\[
 t_m<0\quad\Longleftrightarrow\quad
 -c=Q(t_m)<Q(0)\quad\Longleftrightarrow\quad c>-Q(0).
\]
The same comparison gives \(t_m\le0\) if and only if \(c\ge -Q(0)\).
At equality \(t_m=0\) and the last quadratic is \(x^2+1\).
Finally, \(i+i^{-1}=0\) and \(g(i)\ne0\), so
\(f(i)/(i\,g(i))=Q(0)\), which is real by~\eqref{eq:reciprocal-quotient}.
\end{proof}

\begin{proof}[Proof of Theorem~\ref{thm:combined-hurwitz}]
Put \(f=aA_n+b(1+x)A_{n-1}\) and \(g=A_{n-2}\).
We check the hypotheses of Lemma~\ref{lem:reciprocal-pair} and then
evaluate \(Q(u)=f(x)/(xg(x))\) at \(x=i\), where
\(u=x+x^{-1}=0\).

The polynomials \(A_n\) and \((1+x)A_{n-1}\) are monic and
symmetric in degree \(n-1\). Since \(a+b>0\),
the numerator has degree \(n-1\) and leading coefficient \(a+b\).
The polynomial \(A_{n-2}\) is monic of degree \(n-3\), is symmetric
in that degree, and has only simple negative zeros.
Since \(A_1=1\) and \(A_2=1+x\), there is no zero below \(-1\)
to check for \(n=3,4\).

For \(n\ge5\), let \(r<-1\) be the \(j\)th zero of \(A_{n-2}\) from
the left. By Corollary~\ref{cor:typeA-two-step}, \(r\) lies between the
\(j\)th and \((j+1)\)st zeros of \(A_n\). To its right lie \(n-1-j\) zeros
of \(A_n\) and \(n-3-j\) zeros of \(A_{n-2}\). Since both polynomials
are monic, these counts give
\(\sgn A_n(r)=(-1)^{n-1-j}=(-1)^{n-3-j}=\sgn A'_{n-2}(r)\).
Also, the recurrence~\eqref{eq:Arec}, evaluated at \(r\), gives
\(A_{n-1}(r)=r(1-r)A'_{n-2}(r)\). It follows that
\[
 \frac{aA_n(r)+b(1+r)A_{n-1}(r)}{A'_{n-2}(r)}
 =a\frac{A_n(r)}{A'_{n-2}(r)}+br(1-r^2)>0.
\]
Here \(A_n(r)/A'_{n-2}(r)>0\) because \(A_n(r)\) and \(A'_{n-2}(r)\) have the same sign,
and \(r(1-r^2)>0\) because \(r<-1\). Since \(a,b\ge0\)
are not both zero, the sum is positive.

The lemma now gives Hurwitz stability exactly when \(c>-Q(0)\).
To evaluate \(Q(0)\), take \(x=i\) and replace \(w\) by
\(w/(1+i)\) in~\eqref{eq:eulerian-egf-evaluation}.
Subtract the constant term \(1\), multiply the remaining series by
\(1+i\), and restore the constant term. Since \((i-1)/(1+i)=i\), we obtain
\[
 \begin{aligned}
 1+\sum_{j\ge1}\frac{A_j(i)}{(1+i)^{j-1}}\frac{w^j}{j!}
 &=1+(1+i)\left(\frac{1-i}{e^{iw}-i}-1\right)\\
 &=\frac{2}{e^{iw}-i}-i
  =\frac{2(e^{-iw}+i)}{(e^{iw}-i)(e^{-iw}+i)}-i\\
 &=\frac{e^{-iw}+i}{1-\sin w}-i
  =\frac{\cos w}{1-\sin w}=\sec w+\tan w.
 \end{aligned}
\]
Comparing coefficients with
\(\sec w+\tan w=\sum_{j\ge0}E_jw^j/j!\)
\cite[Problem~4.2(3)]{Petersen2015} yields
\begin{equation}\label{eq:eulerian-at-i}
 A_j(i)=(1+i)^{j-1}E_j\qquad(j\ge1).
\end{equation}
Since \(A_{n-2}\) has only real zeros, \(A_{n-2}(i)\ne0\).
Using~\eqref{eq:eulerian-at-i}, \(i+i^{-1}=0\), and \((1+i)^2=2i\), we obtain
\[
 \begin{aligned}
 Q(0)&=\frac{aA_n(i)+b(1+i)A_{n-1}(i)}{iA_{n-2}(i)}\\
 &=\frac{(1+i)^2(aE_n+bE_{n-1})}{iE_{n-2}}
  =2\frac{aE_n+bE_{n-1}}{E_{n-2}}.
 \end{aligned}
\]
As \(E_{n-2}>0\), the sign of \(c+Q(0)\) is the sign of
\(2aE_n+2bE_{n-1}+cE_{n-2}\).
The strict and weak inequalities in Lemma~\ref{lem:reciprocal-pair}
therefore give the two criteria in the theorem.
\end{proof}

\Needspace{9\baselineskip}
\begin{corollary}\label{cor:hurwitz-boundary}
Let \(n\ge3\), \(a,b\ge0\) with \(a+b>0\), and \(c\in\mathbb R\).
The polynomial
\[
 aA_n(x)+b(1+x)A_{n-1}(x)+cxA_{n-2}(x)
\]
has simple zeros at \(i\) and \(-i\) if \(2aE_n+2bE_{n-1}+cE_{n-2}=0\).
If \(2aE_n+2bE_{n-1}+cE_{n-2}<0\), it has exactly two zeros
in the open right half-plane, counted with multiplicity. In both
cases its other \(n-3\) zeros are simple and negative.
\end{corollary}

\begin{proof}
For \(u=x+x^{-1}\), write
\[
 Q(u)=\frac{aA_n(x)+b(1+x)A_{n-1}(x)}{xA_{n-2}(x)}.
\]
The evaluation at \(x=i\) in the proof of Theorem~\ref{thm:combined-hurwitz}
gives \(Q(0)=2(aE_n+bE_{n-1})/E_{n-2}\). The largest solution \(t_m\)
of \(Q(u)=-c\) lies in the same interval as \(0\), on which \(Q\)
is strictly increasing. Thus \(t_m=0\) when \(c+Q(0)=0\), and
\(t_m>0\) when \(c+Q(0)<0\). For \(t_m=0\), the quadratic
\(x^2-t_mx+1\) is \(x^2+1\), with simple zeros \(i\) and \(-i\).
For \(t_m>0\), both zeros of \(x^2-t_mx+1\) have positive real part, since
\(t_m=(1+|x|^{-2})\operatorname{Re}x\) at either zero.

By Lemma~\ref{lem:reciprocal-pair}, \(t_j<-2\) for \(j<m\), so
\(x^2-t_jx+1\) contributes two distinct negative zeros, neither
equal to \(-1\). Distinct quadratics have no common zero,
and the factor \(1+x\) occurs once when \(n-1\) is odd. This gives
the remaining \(n-3\) simple negative zeros.
\end{proof}

Applying Theorem~\ref{thm:combined-hurwitz} at index \(n+1\)
with \((a,b,c)=(0,1,k)\)
gives, for \(n\ge2\),
\[
 (1+x)A_n(x)+kxA_{n-1}(x)\text{ is Hurwitz stable}
 \quad\Longleftrightarrow\quad k>-2\frac{E_n}{E_{n-1}}.
\]
This proves Conjecture~\ref{conj-02}. At \(k=-2E_n/E_{n-1}\),
Corollary~\ref{cor:hurwitz-boundary} at index \(n+1\) shows that \(i\) and \(-i\) are
simple zeros and that the other \(n-2\) zeros are simple and negative.

At \(n=3\) we have
\((1+x)A_3(x)+kxA_2(x)=(1+x)\bigl(x^2+(4+k)x+1\bigr)\).
For \(-6<k<-2\), the quadratic has two nonreal zeros with real part
\(-(4+k)/2\). At \(k=-4\), the two zeros are \(i\) and \(-i\).
They lie in the open left half-plane for \(-4<k<-2\) and in the open
right half-plane for \(-6<k<-4\). In particular, \(-4<k<-3\) gives Hurwitz stability
beyond Yang and Zhang's sufficient range \(k\ge-3\) for weak Hurwitz stability.
At \(k=-6\), the quadratic is \((x-1)^2\), so \(1\) is a double zero
in the open right half-plane.

In Theorem~\ref{thm:combined-hurwitz}, take \(n\ge3\), set
\((a,b)=(1,0)\), and replace \(c\) by \(-c\). The polynomial
\(A_n-cxA_{n-2}\) is then Hurwitz stable exactly when
\(2E_n-cE_{n-2}>0\), or \(c<2E_n/E_{n-2}\).

To compare the bound \(2E_n/E_{n-2}\) for Hurwitz stability with the
endpoints of the real-rootedness range, write \(A_n(x)/(xA_{n-2}(x))\)
as \(Q(u)\) for \(u=x+x^{-1}\)
and as \(R(y)\) for \(y=((1+x)/(1-x))^2\), as in~\eqref{rhotau}.
For real \(u\ne2\) away
from the poles of \(Q\),
\[
 Q(u)=R\!\left(\frac{u+2}{u-2}\right),\qquad
 Q'(u)=-\frac{4}{(u-2)^2}
 R'\!\left(\frac{u+2}{u-2}\right)>0.
\]
Indeed, dividing the numerator and denominator of $y=(1+x)^2/(1-x)^2$
by $x$ gives $y=(u+2)/(u-2)$, since \(u=x+x^{-1}\).
This substitution gives the first identity, and the chain rule gives
the second. The inequality follows from \(R'(y)<0\)
in~\eqref{eq:Rprime-negative}. At \(u=2\) we have \(x=1\), so direct
substitution gives \(Q(2)=n(n-1)=R(\infty)\). At \(u=0\) we have
\(y=-1\) and \(x=\pm i\). Hence \(Q(0)=R(-1)\), and
\(A_n-cxA_{n-2}\) is Hurwitz stable exactly when \(c<R(-1)\).

Since the numerator and denominator of \(R\) have the same degree,
$R(y)$ tends to $R(\infty)=n(n-1)$ as $y\to\pm\infty$.
Since \(R\) strictly
decreases on \((-\infty,0]\), we obtain
\(R(0)<R(-1)<R(\infty)\). Using~\eqref{eq:Zhang-endpoints}
with \(r=\lfloor n/2\rfloor\) gives
\[
 \frac{E_{2r+1}}{E_{2r-1}}<\frac{2E_n}{E_{n-2}}<n(n-1).
\]
For \(E_{2r+1}/E_{2r-1}<c<n(n-1)\), the polynomial
\(A_n-cxA_{n-2}\) has exactly one nonreal conjugate pair, which
lies on the unit circle by Corollary~\ref{cor:zhang-boundaries}.
Apply Lemma~\ref{lem:reciprocal-pair} with \(f=A_n\), \(g=A_{n-2}\),
and the lemma's parameter replaced by \(-c\). The equation is therefore
\(Q(u)=c\), with solutions \(t_1<\cdots<t_m\) for
\(m=\lfloor(n-1)/2\rfloor\). For \(j<m\), we have \(t_j<-2\).
Each such \(t_j\) gives two negative zeros, so the nonreal pair comes from
\(x^2-t_mx+1\). All poles of \(Q\) lie below \(-2\).
The rightmost interval therefore contains \(0\) and \(t_m\).
Since \(Q\) increases there and \(Q(t_m)=c\), the sign of \(t_m\)
is the sign of \(c-Q(0)\). Both zeros have real part \(t_m/2\), so they lie in the open left
half-plane when
\(E_{2r+1}/E_{2r-1}<c<2E_n/E_{n-2}\),
and in the open right half-plane when
\(2E_n/E_{n-2}<c<n(n-1)\).
At \(c=2E_n/E_{n-2}\), the two zeros are \(i\) and \(-i\).

\section{A combination of Eulerian polynomials of types \texorpdfstring{\(A\)}{A} and \texorpdfstring{\(B\)}{B}}
\label{sec:ABD}
\label{sec:mixed-typeD}

For \(n\ge2\) and \(t\in\mathbb R\), consider
\begin{equation}\label{eq:mixed-typeD}
 B_n(x)-2^{n-1}t\,xA_{n-1}(x).
\end{equation}
At \(t=0\) this is the type~\(B\) Eulerian polynomial. At \(t=n\),
Stembridge's identity gives the type~\(D\) Eulerian polynomial
\(D_n(x)=B_n(x)-n2^{n-1}xA_{n-1}(x)\)~\cite[(16)]{YangZhang2022}. Here \(D_n\) counts
signed permutations \(w=w_1\cdots w_n\) with an even number of negative entries
by their type~\(D\) descents. A type~\(D\) descent is an index \(0\le i<n\) for which
\(w_i>w_{i+1}\), with \(w_0=-w_2\).
Br\"and\'en proved real-rootedness for \(t\le0\)
\cite[p.~3712, after Corollary~6.8]{Branden2006}.
For integral \(0\le t\le n\), \eqref{eq:mixed-typeD} is the polynomial
\(D_n^{(\ell)}\) of Yang and Zhang with \(\ell=n-t\).
They proved that it has only real zeros~\cite[(28) and Theorem~5.1]{YangZhang2022}.

Their argument extends to every real \(t\le n\). We determine the
exact range of \(t\) for which the polynomial has only real zeros.

The values \(B_n(1)\) and \(A_{n-1}(1)\) count all signed permutations
of \([n]\) and all permutations of \([n-1]\), respectively. Hence
\[
 B_n(1)-2^{n-1}tA_{n-1}(1)=2^n n!-2^{n-1}t(n-1)!
           =2^{n-1}(n-1)!(2n-t).
\]
Thus \(x=1\) is a zero exactly when \(t=2n\), the larger endpoint
in the following proposition. At the smaller endpoint, a reciprocal
pair meets at \(-1\). For odd \(n\), the polynomial already has a
zero at \(-1\) for every \(t\).

\begin{proposition}\label{prop:mixed-typeD-range}
Let \(n\ge2\), \(t\in\mathbb R\), and \(m=\lceil n/2\rceil\).
The polynomial
\(B_n(x)-2^{n-1}t\,xA_{n-1}(x)\)
has only real zeros if and only if
\[
 t\le 2\frac{E_{2m}}{E_{2m-1}}
 \qquad\text{or}\qquad t\ge2n.
\]
\end{proposition}

For \(xA_{n-1}(x)\ne0\), the polynomial in~\eqref{eq:mixed-typeD}
vanishes exactly when
\[
 \frac{B_n(x)}{2^{n-1}xA_{n-1}(x)}=t.
\]
The denominator has degree \(n-1\). Since \(A_{n-1}\) is symmetric
in degree \(n-2\), we have
\[
 x^n\bigl(x^{-1}A_{n-1}(1/x)\bigr)=xA_{n-1}(x).
\]
The numerator \(B_n\) satisfies \(x^nB_n(1/x)=B_n(x)\), so the
quotient has the same value at \(x\) and \(1/x\).
The change \(x\mapsto1/x\) reverses the sign of \((1+x)/(1-x)\).
Thus \(y=((1+x)/(1-x))^2\) represents both members of a reciprocal pair
by a single value. The limits of the quotient at \(x=-1\) and \(x=1\)
give the two endpoints of the real-rootedness range in
Proposition~\ref{prop:mixed-typeD-range}, respectively.

\begin{lemma}\label{lem:mixed-ab-quotient}
Let \(n\ge2\) and \(m=\lceil n/2\rceil\). The quotient
$B_n(x)/(2^{n-1}xA_{n-1}(x))$ is a rational function of
\(y=((1+x)/(1-x))^2\). In lowest terms,
its numerator and denominator in \(y\) both have degree
\(\lfloor n/2\rfloor\), and all their zeros are simple and positive. The numerator
strictly alternates left of the denominator, whose largest zero is \(1\).
The values at \(y=0\) and \(y=\infty\) satisfy
\[
 \lim_{x\to-1}\frac{B_n(x)}{2^{n-1}xA_{n-1}(x)}
 =2\frac{E_{2m}}{E_{2m-1}}<2n
 =\frac{B_n(1)}{2^{n-1}A_{n-1}(1)}.
\]
For even \(n\), direct substitution at \(x=-1\) gives the first
limit. For odd \(n\), \(B_n(x)\) and \(A_{n-1}(x)\) have a common
simple zero at \(-1\), so we cancel \(1+x\) before evaluating there.
\end{lemma}

\begin{proof}
Put \(z=(1+x)/(1-x)\), so \(x=(z-1)/(z+1)\) and \(y=z^2\).
By the definitions in Section~\ref{sec:two-parameter},
\(C_{n-1}(z)=(z+1)^{n-2}A_{n-1}(x)\) and \(H_n(z)=2^{-n}(z+1)^nB_n(x)\).
Since \((z+1)^2x=z^2-1\), we obtain
\[
 \frac{B_n(x)}{2^{n-1}xA_{n-1}(x)}
 =\frac{2H_n(z)}{(z^2-1)C_{n-1}(z)}.
\]
The polynomials \(H_n\) and \((z^2-1)C_{n-1}\) have degree \(n\).
They are even when \(n\) is even and odd when \(n\) is odd.
Indeed, \(C_{n-1}\) has degree \(n-2\) and satisfies
\(C_{n-1}(-z)=(-1)^{n-2}C_{n-1}(z)\).
The operator \(\Lop_1\) reverses parity. If \(h\) has leading term
\(az^j\), then \(\Lop_1 h\) has leading term \((j+1)az^{j+1}\). Since \(H_0=1\) and
\(H_{j+1}=\Lop_1H_j\), it follows that \(H_n\) has degree \(n\)
and leading coefficient \(n!\).
By the product rule and~\eqref{L2C2},
\[
 \begin{aligned}
 \Lop_0\bigl((z^2-1)C_{n-1}\bigr)
 &=(z^2-1)\bigl(2zC_{n-1}+(z^2-1)C'_{n-1}\bigr)\\
 &=(z^2-1)C_n.
 \end{aligned}
\]
To pass from \(n\) to \(n+1\), apply \(\Lop_0\) to
\((z^2-1)C_{n-1}\) and \(\Lop_1\) to \(H_n\).
At \(n=2\), \((z^2-1)C_1=z^2-1\) vanishes at \(z=1\)
and \(H_2=2z^2-1\) vanishes at \(z=1/\sqrt2\).
Apply Theorem~\ref{thm:param-interlace} at each step with
\((\mu,\nu)=(0,1)\). All zeros of \((z^2-1)C_{n-1}\) are simple
and lie in \([-1,1]\), with a zero at each endpoint. All zeros
of \(H_n\) are simple and lie in \((-1,1)\). In increasing order,
the positive zeros of the two polynomials alternate strictly and end
at the zero \(1\) of \((z^2-1)C_{n-1}\).

For odd \(n\), cancel the common simple factor \(z\), which corresponds
to \(1+x=2z/(z+1)\). In either parity the resulting numerator and
denominator are even in \(z\), hence polynomials in \(y=z^2\) of
degree \(\lfloor n/2\rfloor\). Their quotient \(R(y)\) satisfies
$R(z^2)=2H_n(z)/\bigl((z^2-1)C_{n-1}(z)\bigr)$ as an identity of
rational functions. The zeros of the numerator and denominator in \(y\)
are the squares of the positive zeros of the corresponding polynomials
in \(z\). The zeros in \(y\) alternate strictly in increasing order,
beginning with a zero of the numerator and ending with the denominator
zero \(1\).
Neither polynomial in \(y\) vanishes at zero. Their leading coefficients
are those of \(2H_n\) and \((z^2-1)C_{n-1}\), namely \(2n!\) and
\((n-1)!\). Both are positive, and their ratio is \(2n!/(n-1)!=2n\).
Summing over the zeros \(a\) of the reduced denominator, we obtain
from Lemma~\ref{lem:residue-representation}
\[
 R(y)=2n+\sum_a\frac{\rho_a}{y-a},\qquad
 \rho_a>0,\qquad R'(y)=-\sum_a\frac{\rho_a}{(y-a)^2}<0
\]
for real \(y\) outside the poles. Since \(x\to1\) as \(z\to\infty\),
the original Eulerian polynomials give the same limiting value,
\[
 R(\infty)=\frac{B_n(1)}{2^{n-1}A_{n-1}(1)}
           =\frac{2^n n!}{2^{n-1}(n-1)!}=2n.
\]

To find \(R(0)\), we need \(H_{2m}(0)\). The recurrence
\(H_{j+1}=\Lop_1H_j\) with \(H_0=1\) gives
\[
 \sum_{j\ge0}H_j(z)\frac{w^j}{j!}
 =\frac{1}{\cosh w-z\sinh w}.
\]
Indeed, the right-hand side equals one at \(w=0\).
Differentiating with respect to \(w\) gives the same result as applying
\(\Lop_1\) in \(z\), namely
\((z\cosh w-\sinh w)/(\cosh w-z\sinh w)^2\).
Comparing the coefficients of \(w^j/j!\) gives the recurrence for \(H_j\).
At \(z=0\) the series becomes \(\sech w=\sec(iw)\), whose
coefficient of \(w^{2m}/(2m)!\) is \((-1)^mE_{2m}\).

Combining this with the evaluation of \(C_{2m-1}(0)\)
from Section~\ref{sec:Zhang} gives
\[
 H_{2m}(0)=(-1)^mE_{2m},\qquad
 C_{2m-1}(0)=A_{2m-1}(-1)=(-1)^{m-1}E_{2m-1}.
\]
If \(n=2m\), direct substitution gives
\(R(0)=-2H_{2m}(0)/C_{2m-1}(0)\).
If \(n=2m-1\), first cancel the common factor \(z\).
At \(z=0\) the recurrences give \(H_{j+1}(0)=-H_j'(0)\) and
\(C_{j+1}(0)=-C_j'(0)\), so
\[
 R(0)=-\frac{2H'_{2m-1}(0)}{C'_{2m-2}(0)}
     =-\frac{2H_{2m}(0)}{C_{2m-1}(0)}.
\]
In the variable \(x\), we substitute \(x=-1\) for even \(n\).
For odd \(n\), the limit is \(-B_n'(-1)/(2^{n-1}A_{n-1}'(-1))\).
Thus in both cases
\[
 R(0)=2\frac{E_{2m}}{E_{2m-1}},\qquad 0<R(0)<2n,
\]
where the strict upper bound follows from
\(R(0)=2n-\sum_a\rho_a/a\), since each pole \(a\) and residue
\(\rho_a\) is positive.
\end{proof}

\begin{proof}[Proof of Proposition~\ref{prop:mixed-typeD-range}]
Put \(z=(1+x)/(1-x)\), \(y=z^2\), and write
\[
 R(y)=\frac{B_n(x)}{2^{n-1}xA_{n-1}(x)}
     =\frac{2H_n(z)}{(z^2-1)C_{n-1}(z)}.
\]
For odd \(n\), cancel the common factor \(z\) before evaluating
\(R(0)\). The canceled factor \(z\) contributes one copy of the zero at
\(x=-1\) to \(B_n-2^{n-1}t\,xA_{n-1}\), since \(1+x=2z/(z+1)\).

Multiplying~\eqref{eq:mixed-typeD} by \(2^{1-n}(z+1)^n\) gives
\begin{equation}\label{eq:mixed-transform}
 2^{1-n}(z+1)^n\bigl(B_n(x)-2^{n-1}t\,xA_{n-1}(x)\bigr)
 =2H_n(z)-t(z^2-1)C_{n-1}(z).
\end{equation}
The recurrence~\eqref{eq:Hrec} gives \(H_{j+1}(1)=H_j(1)\) and
\(H_{j+1}(-1)=-H_j(-1)\). Since \(H_0=1\), the right side of
\eqref{eq:mixed-transform} has values \(2\) at \(z=1\) and
\(2(-1)^n\) at \(z=-1\). Thus each finite zero in \(z\) gives
a finite, nonzero value of \(x=(z-1)/(z+1)\). A zero at \(x=1\)
corresponds to \(z=\infty\).

Write \(h(y)\) and \(f(y)\) for the polynomials obtained from
\(2H_n(z)\) and \((z^2-1)C_{n-1}(z)\), respectively, by dividing
both by \(z\) when \(n\) is odd and then setting \(y=z^2\).
Thus \(R=h/f\) in lowest terms. The right side of~\eqref{eq:mixed-transform}
is \(h(z^2)-tf(z^2)\) for even \(n\) and \(z(h(z^2)-tf(z^2))\)
for odd \(n\). At a pole \(a\) of \(R\), \(h(a)-tf(a)=h(a)\ne0\).

By Lemma~\ref{lem:mixed-ab-quotient}, \(h\) strictly alternates left
of \(f\), whose largest zero is \(1\). For each real \(t\),
Lemma~\ref{lem:remaining} gives one simple solution of \(R(y)=t\)
in each of the \(\lfloor n/2\rfloor-1\) gaps between successive poles.
Each solution \(y\) in a gap satisfies \(0<y<1\). Its two square
roots lie in \((-1,1)\) and give two distinct negative zeros of
\(B_n-2^{n-1}t\,xA_{n-1}\). Both are simple, since \(2z\ne0\)
at either square root and \(dx/dz=2/(z+1)^2\ne0\).
The gaps contribute \(2(\lfloor n/2\rfloor-1)\) zeros. For odd \(n\),
the canceled factor \(z\) gives one further zero at \(x=-1\).
The gap solutions and the factor \(z\) for odd \(n\) account for
\(n-2\) zeros of the original polynomial, counted with multiplicity.
Since \(f(0)\ne0\), \(y=0\) solves \(R(y)=t\) exactly when \(t=R(0)\).

To the left of its smallest pole, \(R\) decreases from \(2n\) to
\(-\infty\). On \((1,\infty)\), \(R\) decreases from \(+\infty\) to \(2n\).
For \(t<R(0)\), the remaining solution of \(R(y)=t\) lies between
zero and the smallest pole. The solution is negative for \(R(0)<t<2n\)
and exceeds one for \(t>2n\). In all three ranges \(R'(y)<0\) at the
solution, so the solution is simple.
By~\eqref{eq:cayley-quadratic-factor}, the corresponding zeros satisfy
\(y(1-x)^2=(1+x)^2\), so
\[
 (1-y)x^2+2(1+y)x+(1-y)=0.
\]
The discriminant is \(16y\), and the two zeros have product one
and sum \(-2(1+y)/(1-y)\). For \(t<R(0)\), the discriminant is positive
and the sum is negative, so both zeros are negative. For \(R(0)<t<2n\),
the discriminant is negative. The zeros are conjugates with product one and therefore have modulus one. For \(t>2n\), the discriminant
and the sum are positive, so both zeros are positive. In all three
ranges the quadratic has two distinct zeros. Its values at \(-1\)
and \(1\) are \(-4y\) and \(4\), respectively, and \(y\ne0\).
Thus at either zero,
\[
 \frac{dy}{dx}=\frac{4(1+x)}{(1-x)^3}\ne0.
\]
Since \(y\) is a simple solution of \(R(y)=t\), both zeros are simple
in the original polynomial. Their common value of \(y\) differs from
the values obtained in the gaps.

At \(t=R(0)\), we have
\[
 2H_n(z)-R(0)(z^2-1)C_{n-1}(z)
 =(z^2-1)C_{n-1}(z)\bigl(R(z^2)-R(0)\bigr).
\]
For even \(n\), the factor \((z^2-1)C_{n-1}(z)\) has value
\(-C_{n-1}(0)\ne0\) at \(z=0\). For odd \(n\), dividing this factor
by \(z\) gives value \(-C'_{n-1}(0)\ne0\) there.
The partial fraction expansion in the proof of
Lemma~\ref{lem:mixed-ab-quotient} gives $R'(0)<0$, so
$(R(z^2)-R(0))/z^2$ has the nonzero limit $R'(0)$ as $z\to0$.
Thus \(h(z^2)-R(0)f(z^2)=f(z^2)(R(z^2)-R(0))\) vanishes to order exactly two
at \(z=0\). Restoring the factor \(z\) raises the multiplicity to
three for odd \(n\). Since \(x+1=2z/(z+1)\) and the factor
\((z+1)^n\) in~\eqref{eq:mixed-transform} is nonzero at \(z=0\),
the multiplicity at \(x=-1\) is two for even \(n\) and three for odd \(n\).

At \(t=2n=R(\infty)\), the partial fractions give
$\lim_{y\to\infty}y(R(y)-2n)=\sum_a\rho_a>0$, where \(\rho_a\)
is the residue of \(R\) at its pole \(a\).
The denominator \(f\) has degree
\(\lfloor n/2\rfloor\), so \(h(y)-2nf(y)=f(y)(R(y)-2n)\) has degree
exactly \(\lfloor n/2\rfloor-1\). Substituting \(y=z^2\) and restoring
the factor \(z\) for odd \(n\) gives degree \(n-2\) in \(z\) on the right
of~\eqref{eq:mixed-transform}. Since \(1-x=2/(z+1)\),
equation~\eqref{eq:mixed-transform} gives
\[
 \frac{B_n(x)-2^n n xA_{n-1}(x)}{(1-x)^2}
 =2^{n-3}\frac{2H_n(z)-2n(z^2-1)C_{n-1}(z)}{(z+1)^{n-2}}.
\]
The numerator and denominator on the right both have degree
\(n-2\), so their quotient has a finite, nonzero limit as
\(z\to\infty\). Thus \(x=1\) has multiplicity exactly two,
as in Lemma~\ref{lem:cayley-mult}. For odd \(n\), the zero
\(-1\) contributed by \(z\) is simple, since \(R(0)\ne2n\).
At both endpoints the zeros from the gaps are simple and negative.
Thus all \(n\) zeros are real exactly when \(t\le R(0)\) or \(t\ge2n\).
Lemma~\ref{lem:mixed-ab-quotient} gives \(R(0)=2E_{2m}/E_{2m-1}\).
\end{proof}

\begin{remark*}
Let \(n\ge2\) and \(m=\lceil n/2\rceil\). For
\(B_n(x)-2^{n-1}t\,xA_{n-1}(x)\), the proof of
Proposition~\ref{prop:mixed-typeD-range} gives the multiplicities
at both endpoints. At \(t=2E_{2m}/E_{2m-1}\), the zero \(-1\) is
double for even \(n\) and triple for odd \(n\). At \(t=2n\), the
zero \(1\) is double. At either endpoint the remaining zeros are simple
and negative. If \(t\) differs from both endpoints, all zeros are simple.
\end{remark*}

At \(n=3\) we have
\(B_3(x)-4t xA_2(x)=(1+x)\bigl(x^2+(22-4t)x+1\bigr)\).
The combination equals \((1+x)^3\) at \(t=5\) and
\((1+x)(x-1)^2\) at \(t=6\). For \(5<t<6\), the quadratic factor has
discriminant \((22-4t)^2-4<0\).

Theorem~\ref{thm:combined-hurwitz} and the Hermite--Biehler theorem
give another proof that \eqref{eq:mixed-typeD} has only negative zeros
for \(t\le2E_n/E_{n-1}\). Apply Theorem~\ref{thm:combined-hurwitz}
at index \(n+1\) with \((a,b,c)=(0,1,-t)\). The polynomial
\((1+u)A_n(u)-tuA_{n-1}(u)\) is weakly Hurwitz stable exactly in
that range.
Multiplying by \((1+u)^n\) preserves weak Hurwitz stability.
The following expansion identifies the even part of the product
with the polynomial in Proposition~\ref{prop:mixed-typeD-range}
evaluated at \(x=u^2\).
\[
 \begin{aligned}
 &(1+u)^n\bigl((1+u)A_n(u)-tuA_{n-1}(u)\bigr)\\
 &\qquad=B_n(u^2)-2^{n-1}t u^2 A_{n-1}(u^2)
          +u\bigl(2^nA_n(u^2)-tB_{n-1}(u^2)\bigr).
 \end{aligned}
\]
For integers \(0\le t\le n\), the identity follows from
\cite[(29)--(30)]{YangZhang2022} with \(\ell=n-t\). Both sides are affine
in \(t\), so agreement at
\(t=0\) and \(t=n\) proves the identity for every real \(t\).

In the stated range, the factor \((1+u)A_n(u)-tuA_{n-1}(u)\) is
monic with constant term \(1\). Each real zero
contributes a factor \(u+a\) with \(a\ge0\).
Each conjugate pair \(-a\pm bi\), with \(a\ge0\) and \(b\ne0\),
contributes a factor \(u^2+2au+a^2+b^2\). Hence
\((1+u)A_n(u)-tuA_{n-1}(u)\) has nonnegative coefficients.
Multiplication by \((1+u)^n\) gives a product with
constant term \(1\) and linear coefficient at least \(n\), so the
even and odd parts are both nonzero.
The Hermite--Biehler theorem for weakly Hurwitz stable polynomials
therefore gives only nonpositive zeros for
\(B_n(x)-2^{n-1}t\,xA_{n-1}(x)\), the polynomial in \(x=u^2\)
formed from the even powers. Its constant term is \(1\), so all its
zeros are negative
\cite[Theorems~1.5 and~5.3]{YangZhang2022}.
Yang and Zhang's Theorem~5.3 gives weak Hurwitz stability of the
factor for \(t\le n\), using \(k=-t\) and their type~\(A\) polynomial
of rank \(n-1\), which is \(A_n\) in our notation. For even \(n\),
the range \(t\le2E_n/E_{n-1}\) reaches the smaller endpoint in
Proposition~\ref{prop:mixed-typeD-range}.

To apply Proposition~\ref{prop:mixed-typeD-range} at \(t=n\), we
compare \(n\) with the smaller endpoint. A count of alternating permutations
shows that \(n\) lies strictly below this endpoint for every \(n\ge3\).

\begin{corollary}\label{cor:mixed-typeD-at-n}
For \(n\ge2\) and \(m=\lceil n/2\rceil\), we have
$n\le 2E_{2m}/E_{2m-1}$, with equality only for \(n=2\).
For \(n\ge3\), the type~\(D\)
Eulerian polynomial
\(D_n(x)=B_n(x)-n2^{n-1}xA_{n-1}(x)\)
has \(n\) simple negative zeros. For \(n=2\), \(D_2(x)=(1+x)^2\).
\end{corollary}

\begin{proof}
Consider the alternating permutations \(w_1\cdots w_{2m}\) of
\([2m]\) that satisfy \(w_{2m}>w_1\). The cyclic inequalities alternate,
so rotation through any even number of positions preserves the condition.
Each orbit under these rotations has \(m\) distinct
permutations. The largest entry lies in an even position, so exactly
one permutation in each orbit ends in \(2m\). Its prefix is an alternating
permutation of \([2m-1]\). Conversely, appending \(2m\) to any such
prefix gives an alternating permutation satisfying \(w_{2m}>w_1\).
Thus \(mE_{2m-1}\) alternating permutations of \([2m]\) satisfy
\(w_{2m}>w_1\).

For \(m\ge2\), prefix \(2m-1,2m\) to any alternating permutation
of \([2m-2]\). The resulting permutation is alternating, but its
last entry is smaller than its first, so it is excluded by \(w_{2m}>w_1\).
Hence
\(E_{2m}>mE_{2m-1}\), so \(n\le2m<2E_{2m}/E_{2m-1}\) for \(n\ge3\).
Thus \(t=n\) lies strictly below the smaller endpoint for
\(n\ge3\). In this range the proof of Proposition~\ref{prop:mixed-typeD-range}
gives \(n\) simple negative zeros. For \(n=2\), we have \(E_1=E_2=1\) and
\(D_2(x)=1+6x+x^2-4x=(1+x)^2\).
\end{proof}

For \(n=4\), the \(2^3\,4!=192\) signed permutations with an even
number of negative entries have type~\(D\) descent polynomial
\(D_4(x)=B_4(x)-32xA_3(x)=1+44x+102x^2+44x^3+x^4\).
For example, \(w=(1,-3,-2,4)\) has \(w_0=-w_2=3\). Its descents
are at \(0\) and \(1\), so it contributes \(x^2\).
Here \(t=4\) lies below the smaller endpoint \(2E_4/E_3=5\), so
all four zeros are simple and negative.

\section{Proof of Theorem~\ref{thm:eulerian-fixed-interlacing}}
\label{sec:eulerian-fixed-interlacing}

\label{sec:proof-c}

For \(n\ge3\) and \(r=\lfloor n/2\rfloor\), Corollary~\ref{cor:zhang-boundaries}
shows that \(A_n-cxA_{n-2}\) has only negative zeros exactly when
\(c\le E_{2r+1}/E_{2r-1}\). At equality, \(-1\) is a multiple zero.
The same bound on \(c\) governs when \((1+x)A_{n-2}\) interlaces
\(A_n-cxA_{n-2}\). At \(c=0\),
Corollary~\ref{cor:typeA-two-step} gives the order of the zeros of
\(A_{n-2}\) and \(A_n\) on each side of \(-1\). To compare the
zeros as \(c\) varies, we divide \(A_n-cxA_{n-2}\) by \((1+x)A_{n-2}\)
and examine the residues.

\begin{proof}[Proof of Theorem~\ref{thm:eulerian-fixed-interlacing}]
Fix \(n\ge3\) and put \(r=\lfloor n/2\rfloor\). Write
\[
 c_n=\frac{E_{2r+1}}{E_{2r-1}},\qquad
 Q(x)=\frac{A_n(x)-cxA_{n-2}(x)}{(1+x)A_{n-2}(x)}.
\]
The contribution of the term \(-cxA_{n-2}\) to \(Q\) is
$-cx/(1+x)=-c+c/(1+x)$.
The residue at \(-1\) therefore differs from its value at \(c=0\)
by \(c\), while all other residues stay fixed. At a removable
singularity the residue is zero.

The numerator \(A_n-cxA_{n-2}\) is monic of degree \(n-1\), since
\(cxA_{n-2}\) has degree at most \(n-2\). The denominator
\((1+x)A_{n-2}\) is monic of degree \(n-2\).
For even \(n\), Corollary~\ref{cor:typeA-two-step} gives the
factor \(1+x\) in both \(A_n\) and \(A_{n-2}\), so
\[
 Q(x)=\frac{\bigl(A_n(x)-cxA_{n-2}(x)\bigr)/(1+x)}{A_{n-2}(x)}.
\]
The denominator \(A_{n-2}\) now has only simple zeros, one of
which is \(-1\). For odd \(n\), all zeros of the original denominator
\((1+x)A_{n-2}\) are already simple.

At \(c=0\), Corollary~\ref{cor:typeA-two-step} gives
one zero of \(A_{n-2}\) in every gap between consecutive zeros
of \(A_n\) on the same side of \(-1\).
For odd \(n\), the factor \(1+x\) places a denominator zero in the
gap between the nearest zeros of \(A_n\) on either side of \(-1\).
For even \(n\), cancel the common factor \(1+x\). The denominator
\(A_{n-2}\) then has a simple zero at \(-1\) in the middle gap
between the zeros of the numerator \(A_n/(1+x)\). Thus
\[
 \begin{aligned}
 (1+x)A_{n-2}&\prec A_n &&(n\text{ odd}),\\
 A_{n-2}&\prec A_n/(1+x) &&(n\text{ even}).
 \end{aligned}
\]
At \(c=0\), the reduced numerator and denominator of \(Q\) are monic,
and the denominator strictly interlaces the numerator. Fix a zero
\(t\) of the denominator. If \(k\) denominator zeros lie to the
right of \(t\), exactly \(k+1\) numerator zeros lie there. Thus the
numerator at \(t\) has sign \((-1)^{k+1}\), and the derivative of
the denominator at \(t\) has sign \((-1)^k\). The residue of
\(Q\) at \(t\) is the numerator's value divided by the denominator's
derivative, so the residue is negative.

For arbitrary \(c\), we compute the residue of \(Q\) at \(-1\).
As \(x\to-1\), we have \(((1+x)/(1-x))^2\to0\), so
\eqref{rhotau} and~\eqref{eq:Zhang-endpoints} give
$A_n(x)/(xA_{n-2}(x))\to c_n$. For even $n$ the common factor
\(1+x\) cancels in this quotient.
The residue is therefore
\[
 \lim_{x\to-1}(1+x)Q(x)
 =\lim_{x\to-1}x\left(\frac{A_n(x)}{xA_{n-2}(x)}-c\right)
 =c-c_n.
\]
At every other denominator zero, \(-c+c/(1+x)\) is regular.
The residue there has the same negative value as at \(c=0\).

For necessity, suppose that \((1+x)A_{n-2}\preceq A_n-cxA_{n-2}\).
When \(n\) is even, cancel one common factor \(1+x\). In the
alternating list of zeros this deletes two adjacent entries at \(-1\),
one from each polynomial. The remaining entries still alternate.
After this cancellation for even \(n\), the denominator has only
simple zeros in either parity.
Write \(t_1,\ldots,t_m\) for the denominator zeros and
\(a_1,\ldots,a_{m+1}\) for the numerator zeros, both in
nondecreasing order and counted with multiplicity. The assumed interlacing gives
$a_1\le t_1\le a_2\le\cdots\le t_m\le a_{m+1}$.
If \(t_i\) is a common zero, the residue there is zero. Otherwise
\(a_i<t_i<a_{i+1}\). There are \(m+1-i\) numerator zeros and
\(m-i\) denominator zeros to the right of \(t_i\). Since both
polynomials are monic, the numerator and the derivative of the
denominator have opposite signs at \(t_i\). The residue is therefore
negative. The nonpositive residue at \(-1\) gives \(c-c_n\le0\).

Suppose now that \(c\le c_n\). After the initial cancellation of
\(1+x\) for even \(n\), the denominator has a simple zero at \(-1\)
in both parities. The residue of \(Q\) at \(-1\) is \(c-c_n\).
Thus the remaining factor \(1+x\) cancels exactly when \(c=c_n\).
At each other zero of the denominator, \(Q\) has a negative residue
and hence a pole.

Since the numerator and denominator of \(Q\) are monic and their
degrees differ by one, the polynomial part of \(Q\) has slope one.
Let \(t_j\) be the
remaining poles and \(\lambda_j<0\) their residues. Differentiating the
partial fraction expansion gives
\[
 Q'(x)=1-\sum_j\frac{\lambda_j}{(x-t_j)^2}\ge1
 \qquad(x\in\mathbb R,\ x\ne t_j\text{ for all }j).
\]
If there are no poles, \(Q\) is linear and has one simple real zero.
Otherwise, at each pole \(Q\) tends to \(+\infty\) from the left and
to \(-\infty\) from the right. Thus \(Q\) increases from \(-\infty\) to \(+\infty\) on each interval
cut out by the poles, including both outer intervals. Each interval
contains exactly one zero, which is simple because \(Q'(x)\ge1\).
The reduced denominator has one simple zero at each pole, and the
reduced numerator has degree one greater. Thus these intervals account
for all zeros of the reduced numerator. In either case the reduced denominator
strictly interlaces the reduced numerator. To restore a common factor
\(x-a\), insert two adjacent copies of \(a\) in the combined list in
nondecreasing order. Use a numerator zero followed by a denominator zero
if the insertion is at the beginning or after a denominator zero. After
a numerator zero, use the reverse order. Each insertion preserves weak
interlacing. Restoring all common factors gives
\((1+x)A_{n-2}\preceq A_n-cxA_{n-2}\).

At \(c=c_n\), the cancellation at \(-1\) removes one factor \(1+x\)
for odd \(n\) and two in total for even \(n\). All other denominator
zeros remain poles, so
\[
 \gcd\!\bigl(A_n(x)-c_nxA_{n-2}(x),(1+x)A_{n-2}(x)\bigr)
 =
 \begin{cases}
  1+x,& n\text{ odd},\\
  (1+x)^2,& n\text{ even}.
 \end{cases}
\]
By Corollary~\ref{cor:zhang-boundaries}, the zero \(-1\) of
\(A_n-c_nxA_{n-2}\) is double for odd \(n\) and triple for even \(n\).
\end{proof}

Let \(n\ge3\) and \(r=\lfloor n/2\rfloor\). For \(c<E_{2r+1}/E_{2r-1}\),
the residue of \((A_n-cxA_{n-2})/((1+x)A_{n-2})\) at \(-1\) is
strictly negative, as are the residues at all other denominator zeros.
When \(n\) is odd, every zero of \((1+x)A_{n-2}\) is simple.
The nonzero residues exclude any common zero, so
\((1+x)A_{n-2}\) strictly interlaces \(A_n-cxA_{n-2}\). When \(n\) is even,
the only common zero is \(-1\). Removing one factor \(1+x\)
from each polynomial gives
\[
 A_{n-2}\prec\frac{A_n-cxA_{n-2}}{1+x}.
\]
Since \(E_3/E_1=2\) and the alternating-permutation count in the
Introduction gives \(E_{2r+1}/E_{2r-1}\ge2\) for every \(r\ge1\),
the interlacing in Theorem~\ref{thm:eulerian-fixed-interlacing}
holds for all \(n\ge3\) exactly when \(c\le2\).

\subsection{The two symmetric summands}

When \(b=0\) and \(0<c\le2\), the second summand in~\eqref{sd} is
\(c(1+x)A_{n-2}\), so Theorem~\ref{thm:eulerian-fixed-interlacing}
gives interlacing of the two summands for every \(n\ge3\).

For \(d\ge1\), every polynomial \(p\) of degree at most \(d\)
has a unique expression \(p=f+xg\), where \(f\) is symmetric in
degree \(d\) and \(g\) is symmetric in degree \(d-1\). The ordered
pair \((f,g)\) is the symmetric decomposition of \(p\) in
degree \(d\)~\cite[Lemma~2.1]{BrandenSolus2021}.
Write \(f_n=A_n-cxA_{n-2}\) and \(g_n=bA_{n-1}+c(1+x)A_{n-2}\)
for the summands in~\eqref{sd}. Since \(x^{j-1}A_j(1/x)=A_j(x)\)
for \(j\ge1\), the reversal of \(A_n+bA_{n-1}+cA_{n-2}\) in degree
\(n-1\) is
\[
 A_n+bxA_{n-1}+cx^2A_{n-2}
 = f_n+xg_n.
\]
Since \(f_n\) and \(g_n\) are symmetric in degrees \(n-1\) and
\(n-2\), respectively, the symmetric decomposition of the
reversal in degree \(n-1\) is \((f_n,g_n)\).

For \(n\ge3\), \(0\le c\le1\), and \(2c\le b\le1+c\),
both summands in~\eqref{sd} have nonnegative coefficients.
If \(g_n\) is nonzero, it interlaces \(f_n\).

\begin{proof}
Athanasiadis's theorem applies to nonnegative linear combinations of
\(t^{n-k}(1+t)^k\) for \(0\le k\le n\). The terms with \(k=0,1,2\) give
\[
 \begin{aligned}
 t^n+bt^{n-1}+ct^{n-2}
 ={}&(1-b+c)t^n+(b-2c)t^{n-1}(1+t)\\
    &+ct^{n-2}(1+t)^2.
 \end{aligned}
\]
The three coefficients on the right are nonnegative exactly when
\(0\le c\le1\) and \(2c\le b\le1+c\).
The Eulerian transformation sends \(t^j\) to \(xA_j(x)\) for \(j\ge1\).
By~\eqref{sd}, the image of \(t^n+bt^{n-1}+ct^{n-2}\) is
\(x(A_n+bA_{n-1}+cA_{n-2})=xg_n+xf_n\).
Since \(g_n\) is symmetric in degree \(n-2\), \(xg_n\) is symmetric
in degree \(n\). The polynomial \(f_n\) is symmetric in degree \(n-1\),
so the symmetric decomposition of the image in degree \(n\) is \((xg_n,f_n)\).
Athanasiadis's theorem applied to \(xg_n+xf_n\) gives nonnegative
coefficients for \(xg_n\) and \(f_n\), hence also for \(g_n\)
\cite[Theorem~1.1(b) and Example~3.3]{Athanasiadis2025}.
Suppose first that \(g_n\ne0\). Then \(b+c>0\) since \(b\ge2c\ge0\).
The leading coefficient of \(g_n\) is \(b+c\), so \(f_n\) and \(xg_n\)
both have degree \(n-1\). Athanasiadis's relation \(f_n\preceq xg_n\)
therefore says that \(f_n\) alternates left of \(xg_n\).
In particular, all zeros of \(f_n\) and \(g_n\) are real. Since both
polynomials have nonnegative coefficients and their constant terms
are \(1\) and \(b+c>0\), respectively, their zeros are negative.
The combined list of zeros of \(f_n\) and \(xg_n\), in decreasing
order, starts with the simple zero \(0\) of \(xg_n\). Removing \(0\)
leaves the zeros of \(g_n\) between
consecutive zeros of \(f_n\), with equality allowed. Thus
\(g_n\preceq f_n\).
If \(g_n=0\), its constant term gives \(b+c=0\). Since
\(b\ge2c\ge0\) in the range under consideration, this forces
\(b=c=0\) and \(f_n=A_n\).
\end{proof}

The bounds \(0\le c\le1\) and \(2c\le b\le1+c\) ensure that the
three coefficients in the basis expansion are nonnegative. We now
ask for which pairs \((b,c)\) with \(c\ge0\) the summands
in~\eqref{sd} satisfy \(g_n\preceq f_n\) for all \(n\ge3\).
We require both polynomials to be nonzero and to have
nonnegative coefficients.

\begin{theorem}\label{thm:ordinary-symmetric}
Fix \(b\in\mathbb R\) and \(c\ge0\), with \((b,c)\ne(0,0)\).
The interlacing
\[
 bA_{n-1}(x)+c(1+x)A_{n-2}(x)
 \preceq A_n(x)-cxA_{n-2}(x)
\]
holds between nonzero polynomials with nonnegative coefficients for
every \(n\ge3\) if and only if \(b\ge0\) and \(0\le c\le2\).
\end{theorem}

When \(b=c=0\), we have \(g_n=0\) and \(f_n=A_n(x)\), whose zeros
are simple and negative.

Theorem~\ref{thm:eulerian-fixed-interlacing} shows that
\((1+x)A_{n-2}\preceq A_n-cxA_{n-2}\) for \(0\le c\le2\).
To include the term \(bA_{n-1}\), we first locate the zeros of
\(A_n-cxA_{n-2}\) relative to those of \(A_{n-1}\).

\begin{lemma}\label{lem:ordinary-middle-interlaces}
Let \(n\ge3\), \(0\le c\le2\), and \((n,c)\ne(3,2)\).
Then \(A_n(x)-cxA_{n-2}(x)\) has \(n-1\) simple negative zeros, and
\[
 \begin{aligned}
 A_{n-1}(x)&\prec A_n(x)-cxA_{n-2}(x),\\
 bA_{n-1}(x)+c(1+x)A_{n-2}(x)&\prec A_n(x)-cxA_{n-2}(x)
 \qquad(b>0).
 \end{aligned}
\]
\end{lemma}

\begin{proof}
Write \(f_n(x)=A_n(x)-cxA_{n-2}(x)\) and
\(g_n(x)=bA_{n-1}(x)+c(1+x)A_{n-2}(x)\).
To determine the sign of \(f_n\) at each zero of \(A_{n-1}\),
we bound the terms of the partial fraction expansion of
\(A_{n-2}/A_{n-1}\) at \(x=1\) and then apply the Eulerian recurrence.

Let \(u_1<\cdots<u_{n-2}<0\) be the zeros of \(A_{n-1}\). Since
\(A_{n-2}\prec A_{n-1}\) and both polynomials are monic,
\(A_{n-2}(u_i)\) and \(A_{n-1}'(u_i)\) both have sign \((-1)^{n-2-i}\). The partial fraction expansion is
\[
 \frac{A_{n-2}(x)}{A_{n-1}(x)}
 =\sum_{i=1}^{n-2}\frac{A_{n-2}(u_i)}{A_{n-1}'(u_i)}\,
                       \frac{1}{x-u_i}.
\]
Set \(x=1\). Since \(A_j(1)=j!\), this gives
\[
 \sum_{i=1}^{n-2}
    \frac{A_{n-2}(u_i)}{(1-u_i)A_{n-1}'(u_i)}
 =\frac{(n-2)!}{(n-1)!}=\frac1{n-1}.
\]
Each summand is positive, so each is at most \(1/(n-1)\).

At \(u_i\), the Eulerian recurrence gives
\(A_n(u_i)=u_i(1-u_i)A_{n-1}'(u_i)\). Thus
\[
 \begin{aligned}
 \frac{A_n(u_i)-cu_iA_{n-2}(u_i)}
      {u_i(1-u_i)A_{n-1}'(u_i)}
 &=1-c\frac{A_{n-2}(u_i)}{(1-u_i)A_{n-1}'(u_i)}\\
 &\ge 1-\frac{c}{n-1}.
 \end{aligned}
\]
Since \(c\le2\), the lower bound is positive for \(n\ge4\).
For \(n=3\), the hypothesis \((n,c)\ne(3,2)\) gives \(c<2\),
so the bound is again positive. As \(u_i(1-u_i)<0\), we obtain
\(\sgn f_n(u_i)=-\sgn A'_{n-1}(u_i)=(-1)^{n-1-i}\).
These signs alternate. The monic polynomial \(f_n\) has degree
\(n-1\), so its sign near \(-\infty\) is \((-1)^{n-1}\). At the
other end, \(f_n(0)=1\). There is therefore a zero in
$(-\infty,u_1)$, one in $(u_{n-2},0)$, and one in each gap
$(u_i,u_{i+1})$ for $1\le i<n-2$.
These \(n-1\) zeros exhaust the degree. They are simple and
negative, and \(A_{n-1}\prec f_n\).

By Theorem~\ref{thm:eulerian-fixed-interlacing}
with \(r=\lfloor n/2\rfloor\) and \(c\le2\le E_{2r+1}/E_{2r-1}\),
we also have \((1+x)A_{n-2}\preceq f_n\).
Let \(u\) be a zero of \(f_n\). The polynomial \(A_{n-1}\) has as
many zeros to the right of \(u\) as \(f_n\) does, so
\(A_{n-1}(u)/f_n'(u)>0\). If \((1+u)A_{n-2}(u)\ne0\), weak
interlacing gives the same count for \((1+x)A_{n-2}\), with
multiplicity. Thus \((1+u)A_{n-2}(u)\) is either zero or has the
sign of \(f_n'(u)\). Since \(c\ge0\), for every \(b>0\) we obtain
\[
 \frac{g_n(u)}{f_n'(u)}
 =b\frac{A_{n-1}(u)}{f_n'(u)}
   +c\frac{(1+u)A_{n-2}(u)}{f_n'(u)}>0.
\]
The signs of \(g_n\) at consecutive zeros of \(f_n\) are opposite.
Since \(g_n\) has degree \(n-2\) and leading coefficient \(b+c>0\),
it has exactly one zero in each gap. This proves \(g_n\prec f_n\).

\end{proof}

For the excluded pair \((n,c)=(3,2)\), substitution of \(A_1=1\),
\(A_2=1+x\), and \(A_3=1+4x+x^2\) gives
\(bA_2(x)+2(1+x)A_1(x)=(b+2)(1+x)\) and \(A_3(x)-2xA_1(x)=(1+x)^2\).
For \(b\ge0\), the common zero \(-1\) is simple in \(g_3\)
and double in \(f_3\).

\begin{proof}[Proof of Theorem~\ref{thm:ordinary-symmetric}]
Write \(f_n=A_n-cxA_{n-2}\) and \(g_n=bA_{n-1}+c(1+x)A_{n-2}\).
Suppose that both polynomials have nonnegative coefficients and that
\(g_n\preceq f_n\) for all \(n\ge3\). In particular \(g_n(1)\ge0\).
Evaluating at \(x=1\) gives
$g_n(1)=b(n-1)!+2c(n-2)!=(n-2)!\bigl((n-1)b+2c\bigr)$.
If \(b<0\), then \((n-1)b+2c<0\) for all sufficiently large
\(n\), contradicting \(g_n(1)\ge0\). Hence \(b\ge0\). For \(n=3\),
$f_3(x)=A_3(x)-cxA_1(x)=1+(4-c)x+x^2$.
The middle coefficient of \(f_3\) gives \(4-c\ge0\). Since
\(f_3\) has only real zeros, its discriminant gives \((4-c)^2\ge4\).
Together these give \(4-c\ge2\), or \(c\le2\).

Conversely, let \(b\ge0\) and \(0\le c\le2\).
Lemma~\ref{lem:ordinary-middle-interlaces} and the factorization
for \((n,c)=(3,2)\) show that all zeros of \(f_n\) are negative.
Since \(f_n\) is monic, its coefficients are positive.
The coefficients of \(g_n\) are nonnegative term by term, and
\(g_n(0)=b+c>0\) whenever \((b,c)\ne(0,0)\).

For \((n,c)=(3,2)\), we have \(g_3=(b+2)(1+x)\) and \(f_3=(1+x)^2\). Thus \(g_3\preceq f_3\)
for every \(b\ge0\). For \((n,c)\ne(3,2)\) and \(b>0\),
Lemma~\ref{lem:ordinary-middle-interlaces} gives \(g_n\prec f_n\).
When \(b=0<c\),
$g_n=c(1+x)A_{n-2}\preceq A_n-cxA_{n-2}=f_n$
by Theorem~\ref{thm:eulerian-fixed-interlacing}.
Finally, if \(b=c=0\), then \(g_n=0\) and \(f_n=A_n\) has only
simple negative zeros.
\end{proof}

In particular, for \(n\ge3\), \(b\ge0\), and \(0\le c\le2\), the
polynomial \(A_n+bA_{n-1}+cA_{n-2}\) has only real zeros.
Indeed, if \(b=c=0\), the sum is \(A_n\). Otherwise \(b+c>0\), so the
two summands in~\eqref{sd} have degrees \(n-1\) and \(n-2\).
Divide both polynomials by their monic greatest common divisor.
In the combined list this removes adjacent pairs of equal zeros,
one from each polynomial. Within each block of equal zeros the entries
alternate between the two polynomials, so removing these pairs preserves
the alternating order and leaves at most one entry in each block.
Every remaining inequality is strict, and
both quotients have positive leading coefficients. Their degrees still
differ by one, so their sum has the degree of the larger quotient.
If the lower-degree quotient is nonconstant, write its zeros as
\(v_1<\cdots<v_d\). At \(v_i\) the sum has the same value as the
higher-degree quotient, which has \(d+1-i\) zeros to the right of \(v_i\).
Hence the sum has sign \((-1)^{d+1-i}\) at \(v_i\). Its sign is
\((-1)^{d+1}\) for all sufficiently negative \(x\) and \(+1\) for all
sufficiently positive \(x\).
There is a zero in each of the \(d-1\) bounded gaps and one in each
outer interval. These \(d+1\) zeros exhaust the degree of the sum.
If the lower-degree quotient is constant, the sum of the two
quotients is linear.
Restoring the common factors proves real-rootedness.
This includes \(c=2\), beyond the range \(c\le1\) supplied
by the nonnegative expansion used in
\cite[Example~3.3]{Athanasiadis2025}.

The case \(c=2\) also gives an example from chromatic polynomials.
Let \(G\) be the disjoint union of a triangle and \(n-3\) isolated
vertices, where \(n\ge3\). In a proper coloring with \(t\) colors,
the triangle uses three distinct colors, and each isolated vertex
can have any color. Thus the chromatic polynomial of \(G\) is
\(\chi_G(t)=t^{n-2}(t-1)(t-2)\). For \(j\ge1\), the Eulerian
transformation sends \(t^j\) to \(xA_j(x)\). Hence it sends
\((-1)^n\chi_G(-t)=t^{n-2}(t+1)(t+2)\) to
$x\bigl(A_n(x)+3A_{n-1}(x)+2A_{n-2}(x)\bigr)$.
The polynomial \((-1)^n\chi_G(-t)\) has the expansion
\[
 t^{n-2}(t+1)(t+2)=2t^{n-2}(1+t)^2-t^{n-1}(1+t).
\]
The polynomials \(t^{n-k}(1+t)^k\), \(0\le k\le n\), have distinct
lowest degrees, so the expansion in this basis is unique.
The coefficient of \(t^{n-1}(1+t)\) is \(-1\), so Athanasiadis's
nonnegativity hypothesis fails.
At \(b=3\) and \(c=2\), Theorem~\ref{thm:ordinary-symmetric} gives
\(3A_{n-1}+2(1+x)A_{n-2}\preceq A_n-2xA_{n-2}\).
The interlacing is strict for \(n\ge4\) by
Lemma~\ref{lem:ordinary-middle-interlaces}. The two polynomials sum to
\(A_n+3A_{n-1}+2A_{n-2}\), so the image has only real zeros for every
\(n\ge3\).

\section{Proof of Theorem~\ref{thm:largest-intro}}
\label{sec:three-term}

For \(n\ge3\), Theorem~\ref{thm:ordinary-symmetric} gives
real-rootedness of \(A_n+bA_{n-1}+cA_{n-2}\) for \(b\ge0\) and
\(0\le c\le2\) by interlacing the two summands in~\eqref{sd}.
For \(b=c=0\) the combination is \(A_n\). We now fix \(n\) and ask
how large \(c\) can be if the combination is to have only real zeros
for every real \(b\). For ordinary Eulerian polynomials, the first two
bounds on \(c\) are \(2\) at \(n=3\) and \(12-6\sqrt3\) at \(n=4\).

Fix one of the five families in Table~\ref{tab:main-bounds} and an
index \(n\) for which \(F_{n-1}\) is nonconstant. At each zero \(r\)
of \(F_{n-1}\), the term \(bF_{n-1}(r)\) vanishes for every \(b\).
For \(F_n+bF_{n-1}+(c-kx)F_{n-2}\), the next lemma shows that
real-rootedness for every \(b\) is equivalent to
\(c\le-F_n(r)/F_{n-2}(r)+kr\) at each such zero. For
Theorem~\ref{thm:largest-intro}, we take \(k=0\) and show that the
largest zero of \(F_{n-1}\) gives the smallest bound.

\subsection{A criterion for three consecutive polynomials}

For \(n\ge1\), the polynomial \(F_n\) has degree \(n-1\) for \(A_n,S_n,N_n\)
and degree \(n\) for \(B_n,M_n\). Whenever \(F_{n-1}\) is
nonconstant, the polynomials \(F_{n-2},F_{n-1},F_n\) have
consecutive degrees. The next lemma gives a criterion for
real-rootedness for every value of the middle coefficient.

\begin{lemma}\label{lem:free-middle-general}
Let \(F_{n-2},F_{n-1},F_n\) be real polynomials of consecutive
degrees, in increasing order, with positive leading coefficients
and \(F_{n-2}\prec F_{n-1}\). Fix \(c,k\in\mathbb R\).
The polynomial
\[
 F_n(x)+bF_{n-1}(x)+(c-kx)F_{n-2}(x)
\]
has only real zeros for every \(b\in\mathbb R\) if and only if
\[
 c\le \min_{F_{n-1}(r)=0}
 \left(-\frac{F_n(r)}{F_{n-2}(r)}+kr\right).
\]
If the inequality is strict, \(F_{n-1}\) strictly interlaces this
polynomial for every real \(b\).
\end{lemma}

\begin{proof}
Put \(h=F_n+(c-kx)F_{n-2}\) and \(R=h/F_{n-1}\). Since \(h\)
has the leading term of \(F_n\), the polynomial part of \(R\)
is linear with slope \(\lc(F_n)/\lc(F_{n-1})>0\).
At a zero \(r\) of \(F_{n-1}\), the residue of \(R\) is
\begin{equation}\label{eq:three-term-residue}
 \gamma_r
 =\frac{h(r)}{F_{n-1}'(r)}
 =\frac{F_{n-2}(r)}{F_{n-1}'(r)}
   \left(c-kr+\frac{F_n(r)}{F_{n-2}(r)}\right).
\end{equation}
The factor \(F_{n-2}(r)/F_{n-1}'(r)\) is positive because
\(F_{n-2}\) and \(F_{n-1}(x)/(x-r)\) have positive leading
coefficients and equally many zeros to the right of \(r\).
Thus the upper bound on \(c\) is equivalent to \(\gamma_r\le0\)
at every zero \(r\) of \(F_{n-1}\).

Suppose \(h+bF_{n-1}\) has only real zeros for every real \(b\).
Since \(F_{n-1}\) has only real zeros, every nonzero real linear
combination of \(h\) and \(F_{n-1}\) has only real zeros.
Obreschkoff's theorem gives \(F_{n-1}\preceq h\)
\cite[Theorem~2.4]{Branden2006}. At each zero \(r\) of \(F_{n-1}\),
either \(h(r)=0\) or \(h\) has one more zero to the right of \(r\)
than \(F_{n-1}\) does. In the latter case \(h(r)\) and
\(F_{n-1}'(r)\) have opposite signs. Thus \(\gamma_r\le0\) in both cases.

Conversely, suppose all residues are nonpositive. Let \(C\) be the
product of the common factors \(x-r\) for which \(\gamma_r=0\), so
\(R=(h/C)/(F_{n-1}/C)\).
Away from the remaining poles, partial fractions give
\[
 R'(x)=\frac{\lc(F_n)}{\lc(F_{n-1})}
       -\sum_{\gamma_r<0}\frac{\gamma_r}{(x-r)^2}>0.
\]
At a pole, \(R\) has limit \(+\infty\) from the left and
\(-\infty\) from the right. The positive slope gives
\(R(x)\to\pm\infty\) as \(x\to\pm\infty\). Hence \(R(x)=-b\)
has one simple solution in each of the \(s+1\) intervals cut out
by the \(s\) remaining poles. When \(s=0\), the interval is \(\mathbb R\).
Away from the poles, these solutions are precisely the zeros of
\(h/C+bF_{n-1}/C=(F_{n-1}/C)(R+b)\), and they are simple. At each
pole, \(F_{n-1}/C\) vanishes and \(h/C\) is nonzero. The original
polynomial is \(C(h/C+bF_{n-1}/C)\), so each of the
\(\deg F_{n-1}-s\) factors of \(C\) adds one zero counted with
multiplicity, giving \(\deg F_{n-1}+1\) real zeros in total.
With strict inequality on \(c\), all zeros of \(F_{n-1}\) remain
poles, so the intervals give strict interlacing.

\hypertarget{proof:free-middle-multiplicity}{If \(\gamma_r=0\)},
then \(r\) is a simple zero of \(C\) and \(R\) is defined there.
The other factor \(h/C+bF_{n-1}/C\) has a simple zero at \(r\)
exactly when \(b=-R(r)\), in which case \(r\) is a double zero of
the original polynomial. Since \(F_{n-1}\) has a simple zero at \(r\),
evaluation after cancellation gives \(R(r)=h'(r)/F_{n-1}'(r)\).
\end{proof}

\begin{remark*}
At equality in Lemma~\ref{lem:free-middle-general},
\(F_n+bF_{n-1}+(c-kx)F_{n-2}\) and \(F_{n-1}\) share exactly the
zeros that attain the minimum, for every real \(b\). Each such zero
\(r\) is double in the combination precisely when
\[
 b=-\frac{F_n'(r)+(c-kr)F_{n-2}'(r)-kF_{n-2}(r)}
          {F_{n-1}'(r)}.
\]
For every other \(b\), \(r\) is simple. All zeros outside this
set of common zeros are simple.
\end{remark*}

The criterion above asks for a bound at every zero of the middle
polynomial. Recall the operator
\[
 T_d^{a,b}f=(a+(d+b)x)f+x(1-x)f'.
\]
When \(g=T_d^{a,b}f\) and \(h=T_{d+1}^{\alpha,\beta}g\), evaluating
\(f/g\) at \(0\) and \(1\) gives an inequality in \(c\) and \(k\)
that does not depend on the zeros of \(f\).

\begin{proposition}\label{prop:uniform-two-step-bound}
Let \(a,b,\alpha,\beta>0\) and \(d\ge0\). Let \(f\) be a
degree-\(d\) polynomial with positive leading coefficient and simple
negative zeros. Put \(g=T_d^{a,b}f\) and \(h=T_{d+1}^{\alpha,\beta}g\).
If \(c\) and \(k\) are nonnegative and satisfy \((b+d)c+ak\le a(a+b+d)\),
then \(h(x)+tg(x)+(c-kx)f(x)\) has only real zeros for every
\(t\in\mathbb R\). The bound is sharp when \(f\) ranges over all
such polynomials of degree \(d\).
\end{proposition}

\begin{proof}
At each zero \(r\) of \(f\), we have \(g(r)=r(1-r)f'(r)\).
Since \(r(1-r)<0\), the values \(g(r)\) alternate in sign along the
zeros of \(f\). The signs at \(0\) and near \(-\infty\) follow from
\(g(0)=af(0)>0\) and the leading coefficient \(b\lc(f)>0\).
These signs give one zero of \(g\) in each interval cut out by the
zeros of \(f\) in \((-\infty,0)\). This accounts for all \(d+1\)
zeros of \(g\), so they are simple and \(f\prec g\).
The same argument applied to \(g\) gives \(g\prec h\), with all
zeros of \(h\) simple and negative.

At every zero \(r\) of \(g\), strict interlacing gives \(f(r)/g'(r)>0\).
In \(h=T_{d+1}^{\alpha,\beta}g\), the term containing \(\alpha,\beta\)
vanishes at \(r\). Thus \(h(r)=r(1-r)g'(r)\) and \(f(r)/(-h(r))>0\).
The derivative term in \(g=T_d^{a,b}f\) vanishes at \(0\) and \(1\),
so \(g(0)=af(0)\) and \(g(1)=(a+b+d)f(1)\).
Evaluating the partial fraction expansion
\[
 \frac{f(x)}{g(x)}=\sum_{g(r)=0}\frac{f(r)}{(x-r)g'(r)}
\]
at \(x=0\) and \(x=1\), respectively, gives
\[
 \sum_{g(r)=0}\frac{f(r)}{-rg'(r)}=\frac1a,
 \qquad
 \sum_{g(r)=0}\frac{f(r)}{(1-r)g'(r)}=\frac1{a+b+d}.
\]
Using \(h(r)=r(1-r)g'(r)\), we rewrite the difference
\(f(0)/g(0)-f(1)/g(1)\) and the value \(f(1)/g(1)\), respectively, as
\[
 \sum_{g(r)=0}\frac{f(r)}{-h(r)}=\frac{b+d}{a(a+b+d)},
 \qquad
 \sum_{g(r)=0}\frac{rf(r)}{h(r)}=\frac1{a+b+d}.
\]
Consequently,
\[
 \sum_{g(r)=0}\frac{(c-kr)f(r)}{-h(r)}
 =\frac{(b+d)c+ak}{a(a+b+d)}\le1.
\]
Every summand is nonnegative, so \((c-kr)f(r)/(-h(r))\le1\)
at each zero \(r\) of \(g\). Equivalently,
\(c\le-h(r)/f(r)+kr\). Lemma~\ref{lem:free-middle-general} now gives
real-rootedness for every real \(t\).

We prove sharpness. For \(d\ge1\), take
\[
 f_R(x)=\prod_{j=1}^d\left(1+\frac{x}{jR}\right),
 \qquad R>0.
\]
For \(d=0\), take \(f_R=1\). Put \(g_R=T_d^{a,b}f_R\) and
\(h_R=T_{d+1}^{\alpha,\beta}g_R\). Then
\[
 \frac{g_R(x)}{f_R(x)}
 =a+(d+b)x+x(1-x)\sum_{j=1}^d\frac1{x+jR}.
\]
As \(R\to\infty\), the quotient \(g_R/f_R\) and its derivative converge
uniformly on compact intervals to \(a+(d+b)x\) and \(d+b\), respectively.
Fix a closed interval with \(-a/(d+b)\) in its interior.
For all sufficiently large \(R\), the quotient has opposite signs at the
endpoints and a positive derivative throughout the interval. Its unique
zero \(r_R\) there tends to \(-a/(d+b)\), since the limiting linear
polynomial has opposite signs on either side of this point.
On the same interval \(f_R\to1\) uniformly. Thus \(f_R(r_R)\to1\) and
\(g_R'(r_R)=f_R(r_R)(g_R/f_R)'(r_R)\to d+b\).
Since \(h_R(r_R)=r_R(1-r_R)g_R'(r_R)\), the limits above give
\[
 -h_R(r_R)\longrightarrow\frac{a(a+b+d)}{b+d},
 \qquad
 c-kr_R\longrightarrow c+\frac{ak}{b+d}.
\]
Together with \(f_R(r_R)\to1\), this gives
\[
 (c-kr_R)\frac{f_R(r_R)}{-h_R(r_R)}
 \longrightarrow
 \frac{(b+d)c+ak}{a(a+b+d)}.
\]
If \((b+d)c+ak>a(a+b+d)\), the limit exceeds one. For all
sufficiently large \(R\), the residue inequality therefore fails at \(r_R\).
Lemma~\ref{lem:free-middle-general} then gives a real \(t\) for which
\(h_R+t g_R+(c-kx)f_R\) has a nonreal zero. For \(d=0\), the same
calculation is exact with \(f_R=1\). This proves sharpness.
\end{proof}

For example, take \(c=0\) and \(f=A_{n-2}(x;q)\), with \(q>0\) and \(n\ge3\),
so that \(a=\alpha=q\), \(b=\beta=1\) and \(d=n-3\) in
Proposition~\ref{prop:uniform-two-step-bound}. By~\eqref{eq:cycle-insertion},
\(A_{m+1}(x;q)=T_{m-1}^{q,1}A_m(x;q)\) for \(m\ge1\), so
\(g=A_{n-1}(x;q)\) and \(h=A_n(x;q)\). Thus
\(A_n(x;q)+tA_{n-1}(x;q)-kxA_{n-2}(x;q)\) has only real zeros
for every real \(t\) whenever \(0\le k\le q+n-2\).
For \(k<0\), every zero \(r\) of \(g\) is negative and
\(f(r)/(-h(r))>0\). Thus \((-kr)f(r)/(-h(r))<0\), and the same
residue criterion applies. Since \(q+n-2\ge q+1\) for \(n\ge3\),
the bound \(k\le q+1\) suffices for all \(n\ge3\).
Theorem~\ref{thm:cycle-free-middle} shows that this uniform bound is sharp.

Returning to the five families in Table~\ref{tab:main-bounds}, we
take \(k=0\) in Lemma~\ref{lem:free-middle-general}. For each family
and each fixed \(n\) with \(F_{n-1}\) nonconstant, this gives the
exact upper bound on \(c\) for \(F_n+bF_{n-1}+cF_{n-2}\) to have
only real zeros for every real \(b\).
For example, the middle polynomial of \(A_3+bA_2+cA_1\) is
\(A_2=1+x\). Its zero \(-1\) gives the upper bound
\(-A_3(-1)/A_1(-1)=2\). At \(c=2\) the factorization
$A_3+bA_2+2A_1=(x+1)(x+b+3)$
shows directly that \(-1\) is double exactly when \(b=-2\).

For \(n=4\), the two zeros of \(A_3\) give two candidate bounds.
Since \(A_4=(1+x)(A_3+6x)\), we have
\(-A_4(r)/A_2(r)=-6r\) at either zero of \(A_3\).
Thus the minimum is \(-6(-2+\sqrt3)\), attained at the larger
zero of \(A_3\). For the three Eulerian families, the next lemma
writes \(-F_n(r)/F_{n-2}(r)\) at a zero \(r\) of \(F_{n-1}\) as a
constant plus a sum of squares with positive weights. The sum is
smallest at the largest zero of \(F_{n-1}\).

\begin{lemma}\label{lem:largest-recurrence}
\label{prop:eulerian-largest-recurrence}
Suppose that \(F_{n-2}\) is a real polynomial with constant
term one and only simple negative zeros. For \(i=n-2,n-1\),
let \(a_i,s_i>0\) and define
\begin{equation}\label{eq:eulerian-degree-recurrence}
 F_{i+1}=\bigl(1+(a_i+s_i\deg F_i)x\bigr)F_i
             +s_ix(1-x)F_i'.
\end{equation}
Then \(F_{n-2}\prec F_{n-1}\prec F_n\), and all three polynomials
have only simple negative zeros. Let \(\rho\) be the largest zero of \(F_{n-1}\).
If \(a_{n-2}\ge1\) and \(s_{n-2}\ge1\), then every other zero \(r\) of
\(F_{n-1}\) satisfies
\[
 -\frac{F_n(\rho)}{F_{n-2}(\rho)}
 <
 -\frac{F_n(r)}{F_{n-2}(r)}.
\]
\end{lemma}

\begin{proof}
Let \(d=\deg F_{n-2}\), and let \(u\) run over its zeros. Since
\(F_{n-2}'/F_{n-2}=\sum_u1/(x-u)\), dividing the first recurrence by
\(F_{n-2}\) gives
\begin{equation}\label{eq:eulerian-degree-quotient}
 \frac{F_{n-1}(x)}{F_{n-2}(x)}
 =1+a_{n-2}x+s_{n-2}\sum_u \frac{x(1-u)}{x-u},
 \qquad
 x+\frac{x(1-x)}{x-u}=\frac{x(1-u)}{x-u}.
\end{equation}
Since \(x=(x-u)+u\), we have
\[
 \frac{x(1-u)}{x-u}=(1-u)+\frac{u(1-u)}{x-u}.
\]
Thus the polynomial part of \(F_{n-1}/F_{n-2}\) is
$a_{n-2}x+1+s_{n-2}\sum_u(1-u)$,
with slope \(a_{n-2}>0\). The residue at each zero \(u\) of
\(F_{n-2}\) is \(s_{n-2}u(1-u)<0\). Away from these poles,
\[
 \left(\frac{F_{n-1}}{F_{n-2}}\right)'(x)
 =a_{n-2}-s_{n-2}\sum_u\frac{u(1-u)}{(x-u)^2}>0
 \qquad(F_{n-2}(x)\ne0,\ x\in\mathbb R).
\]
Thus the quotient increases strictly
from \(-\infty\) to \(+\infty\) on each of the \(d+1\) intervals
cut out by the zeros of \(F_{n-2}\). Each interval contains one simple
zero of \(F_{n-1}\). The rightmost interval contains \(0\), where the
quotient equals one. Since the quotient is increasing, its zero in
that interval is negative. Every other interval lies in \((-\infty,0)\).
At a zero \(u\) of \(F_{n-2}\), direct substitution gives
$F_{n-1}(u)=s_{n-2}u(1-u)F_{n-2}'(u)\ne0$.
The negative zeros and constant term one of \(F_{n-2}\) give
\(\lc(F_{n-2})>0\). The recurrence gives \(F_{n-1}(0)=1\).
In degree \(d+1\), the two terms containing \(s_{n-2}d\) cancel.
The leading coefficient is
\[
 \lc(F_{n-1})
 =(a_{n-2}+s_{n-2}d)\lc(F_{n-2})
   -s_{n-2}d\lc(F_{n-2})
 =a_{n-2}\lc(F_{n-2})>0.
\]
Hence the \(d+1\) zeros just found exhaust the degree of
\(F_{n-1}\), proving \(F_{n-2}\prec F_{n-1}\). For \(d=0\), the quotient is
\(1+a_{n-2}x\), whose only zero is \(-1/a_{n-2}<0\).
Now \(F_{n-1}\) has constant term one and only simple negative
zeros, so applying the argument to \(F_n/F_{n-1}\) gives
\(F_{n-1}\prec F_n\) and shows that the zeros of \(F_n\) are
simple and negative.

To compare the ratios, let \(r\) be a zero of \(F_{n-1}\). Strict interlacing gives
\(F_{n-2}(r)\ne0\), and \eqref{eq:eulerian-degree-quotient} yields
\[
 a_{n-2}r+s_{n-2}\sum_u\frac{r(1-u)}{r-u}=-1.
\]
The second recurrence gives
\(F_n(r)=s_{n-1}r(1-r)F_{n-1}'(r)\).
Since \(F_{n-1}(r)=0\), we have
\((F_{n-1}/F_{n-2})'(r)=F_{n-1}'(r)/F_{n-2}(r)\).
To obtain the sum of squares, add \((F_{n-1}/F_{n-2})(r)=0\)
to \(r(r-1)(F_{n-1}/F_{n-2})'(r)\).
Indeed, the summand associated with \(u\) satisfies
\[
 \frac{x(1-u)}{x-u}
 +x(x-1)\frac{d}{dx}\frac{x(1-u)}{x-u}
 =\frac{x(1-u)\bigl(x-u-u(x-1)\bigr)}{(x-u)^2}
 =\left(\frac{x(1-u)}{x-u}\right)^2 .
\]
For the term \(1+a_{n-2}x\) in~\eqref{eq:eulerian-degree-quotient},
the same calculation gives
\(1+a_{n-2}x+x(x-1)a_{n-2}=1+a_{n-2}x^2\). Thus
\eqref{eq:eulerian-degree-quotient} gives
\begin{equation}\label{eq:eulerian-degree-squares}
 \begin{aligned}
 -\frac{F_n(r)}{s_{n-1}F_{n-2}(r)}
 &=r(r-1)\left(\frac{F_{n-1}}{F_{n-2}}\right)'(r)
                +\frac{F_{n-1}(r)}{F_{n-2}(r)}\\
 &=1+a_{n-2}r^2+s_{n-2}\sum_u
              \left(\frac{r(1-u)}{r-u}\right)^2.
 \end{aligned}
\end{equation}
Now assume \(a_{n-2}\ge1\) and \(s_{n-2}\ge1\). If \(d=0\), \(F_{n-1}\) has just one
zero, $r=-1/a_{n-2}$, so the minimum is
$-F_n(r)/F_{n-2}(r)=s_{n-1}(1+1/a_{n-2})$.
Now suppose \(d>0\). Every zero \(u\) of \(F_{n-2}\) satisfies \(u<\rho<0\).
Hence \(\rho\) and each fraction \(\rho(1-u)/(\rho-u)\) are negative. Thus
\[
 a_{n-2}|\rho|+s_{n-2}\sum_u \left|\frac{\rho(1-u)}{\rho-u}\right|=1.
\]
There are at least two positive summands in this equality, so each
is less than one. Since the weights \(a_{n-2},s_{n-2}\) are at least one,
each of the absolute values is also less than one. Hence
\[
 a_{n-2}\rho^2+s_{n-2}\sum_u\left(\frac{\rho(1-u)}{\rho-u}\right)^2
 <a_{n-2}|\rho|+s_{n-2}\sum_u \left|\frac{\rho(1-u)}{\rho-u}\right|=1.
\]
At every other zero \(r\) of \(F_{n-1}\), the largest zero
\(u\) of \(F_{n-2}\) lies between \(r\) and \(\rho\). Hence
\[
 \frac{r(1-u)}{r-u}-1=\frac{(-u)(1-r)}{u-r}>0,
 \qquad s_{n-2}\left(\frac{r(1-u)}{r-u}\right)^2>1.
\]
The right-hand side of~\eqref{eq:eulerian-degree-squares} is therefore
less than \(2\) at \(\rho\) and greater than \(2\) at \(r\).
Multiplying by \(s_{n-1}>0\) gives
\[
 -\frac{F_n(\rho)}{F_{n-2}(\rho)}<2s_{n-1}
 < -\frac{F_n(r)}{F_{n-2}(r)},
\]
so no other zero of \(F_{n-1}\) attains the value at \(\rho\).
\end{proof}

\subsection{The Narayana polynomials}

The ordinary and type~\(B\) Narayana polynomials are
\begin{equation}\label{eq:narayana-coefficients}
 \begin{aligned}
 N_n(x)&=\frac1n\sum_{k=0}^{n-1}\binom nk\binom n{k+1}x^k
       &&(n\ge1),\\
 M_n(x)&=\sum_{k=0}^{n}\binom nk^2x^k
       &&(n\ge0).
 \end{aligned}
\end{equation}
Take \((F_j)=(N_j)\) or \((F_j)=(M_j)\). We derive three-term
recurrences from~\eqref{eq:narayana-coefficients}.
Evaluating them at a zero \(r\) of \(F_{n-1}\) expresses
\(-F_n(r)/F_{n-2}(r)\) as a positive multiple of \((1-r)^2\).
Write \(d=\deg F_n\), so \(d=n-1\) for \(F_n=N_n\)
and \(d=n\) for \(F_n=M_n\). Formula~\eqref{eq:narayana-coefficients} shows that
both polynomials are monic and symmetric in their respective degrees,
with constant term one.

Write \(a_{n,k}=[x^k]F_n\), with coefficients outside the degree range
taken to be zero. For \(d\ge1\) and \(0\le k\le d\), cancellation
in~\eqref{eq:narayana-coefficients} gives
\[
 \frac{a_{n-1,k}}{a_{n,k}}=\frac{(d-k)(n-k)}{dn},
 \qquad
 \frac{a_{n-1,k-1}}{a_{n,k}}=\frac{k(n-d+k)}{dn}.
\]
The two numerators are \((n-1-k)(n-k)\) and
\(k(k+1)\) for \(N_n\), and \((n-k)^2\) and \(k^2\) for \(M_n\).
Subtracting the second ratio from the first gives \((d-2k)/d\). Hence
\(d(a_{n,k}-a_{n-1,k}+a_{n-1,k-1})=2ka_{n,k}\), which gives
\begin{equation}\label{eq:narayana-derivative}
 2xF_n'=d\bigl(F_n-(1-x)F_{n-1}\bigr).
\end{equation}
The same two ratios give
\[
 (2n-d)(a_{n,k}-a_{n-1,k}-a_{n-1,k-1})
 =2k a_{n-1,k}+2(d-k)a_{n-1,k-1}.
\]
Indeed, dividing both sides by \(a_{n,k}\) gives \(2k(d-k)(2n-d)/(dn)\).
Multiply the identity with factor \(2n-d\) by \(x^k\) and sum over
\(0\le k\le d\). The first term on the right gives \(2xF_{n-1}'\).
Shifting the index by one in the second gives
\(2(d-1)xF_{n-1}-2x^2F_{n-1}'\). Hence
\[
 (2n-d)\bigl(F_n-(1+x)F_{n-1}\bigr)
 =2(d-1)xF_{n-1}+2x(1-x)F_{n-1}'.
\]
When \(d\ge2\), equation~\eqref{eq:narayana-derivative} at index
\(n-1\) gives \(2xF_{n-1}'=(d-1)(F_{n-1}-(1-x)F_{n-2})\).
Substitution gives
\[
 2(d-1)xF_{n-1}+2x(1-x)F_{n-1}'
 =(d-1)(1+x)F_{n-1}-(d-1)(1-x)^2F_{n-2}.
\]
Collecting the two terms in \((1+x)F_{n-1}\) gives
\[
 (2n-d)F_n=(2n-1)(1+x)F_{n-1}-(d-1)(1-x)^2F_{n-2},
\]
since \((2n-d)+(d-1)=2n-1\). Taking \(d=n-1\) and \(d=n\), respectively,
gives
\begin{equation}\label{eq:nar-three}
 \begin{aligned}
 (n+1)N_n&=(2n-1)(1+x)N_{n-1}-(n-2)(1-x)^2N_{n-2}
                  &&(n\ge3),\\
 nM_n&=(2n-1)(1+x)M_{n-1}-(n-1)(1-x)^2M_{n-2}
                  &&(n\ge2).
 \end{aligned}
\end{equation}

Each recurrence gives strict interlacing and simple negative zeros by
induction, starting with \(N_1=M_0=1\) and \(N_2=M_1=1+x\).
Assume that \(F_{n-2}\prec F_{n-1}\) and that the zeros of \(F_{n-1}\)
are \(r_1<\cdots<r_{d-1}<0\). The monic polynomial \(F_{n-2}\) has
\(d-1-i\) zeros to the right of \(r_i\), so
\[
 (2n-d)F_n(r_i)=-(d-1)(1-r_i)^2F_{n-2}(r_i),
 \qquad \sgn F_n(r_i)=(-1)^{d-i}.
\]
The signs of \(F_n(r_i)\) alternate, giving a zero of \(F_n\)
in each of the \(d-2\) intervals between consecutive zeros of
\(F_{n-1}\). Near \(-\infty\)
the sign of \(F_n\) is \((-1)^d\), opposite to its sign at \(r_1\).
At the other end, \(F_n(r_{d-1})<0<F_n(0)\). Thus the outer
intervals \((-\infty,r_1)\) and \((r_{d-1},0)\) each contain a
zero as well. These \(d\) distinct zeros account for the degree of
\(F_n\), so they are simple and negative.

\subsection{The sharp bound for each index}

The recurrences for \(A_j\) and \(B_j\) are already in the form
of Lemma~\ref{lem:largest-recurrence}. For the second-order Eulerian
polynomials, the same form follows by inserting the two copies of
the largest letter.

Gessel and Stanley~\cite{GesselStanley1978} defined a Stirling
permutation of \(\{1,1,\ldots,j,j\}\) to be a permutation of this
multiset in which every entry between the two copies of a letter
is larger than that letter. Let \(S_j(x)\) count these permutations
by descents between adjacent entries. Counting the final position as
a descent gives their second-order Eulerian polynomial \(xS_j(x)\)
for \(j\ge1\). See~\cite{MaQiYehYeh2026} for recent progress on this topic.

The two copies of \(j+1\) in a Stirling permutation of order \(j+1\)
are consecutive. Deleting them gives a Stirling permutation
of order \(j\). Conversely, the pair can be inserted into any of
the \(2j+1\) positions in the shorter permutation, between or at the
ends of the entries.
Inserting the pair between two entries that form a descent replaces
that descent with one from \(j+1\) to the next entry.
For a permutation with \(k\) descents, there are \(k\) such positions.
Insertion after the last entry also preserves
the descent count. Each of the other \(2j-k\) positions creates one
descent. Summing over the insertions into a permutation with \(k\)
descents gives
\[
 (k+1)x^k+(2j-k)x^{k+1}
 =(1+2jx)x^k+x(1-x)\frac{d}{dx}x^k.
\]
Summation over all Stirling permutations of order \(j\) gives
\begin{equation}\label{eq:second-order-recurrence}
 S_{j+1}=(1+2jx)S_j+x(1-x)S_j',\qquad S_0=S_1=1.
\end{equation}
If \(\deg S_j=j-1\), the coefficient of \(x^j\) in the
recurrence is \((2j-(j-1))\lc(S_j)=(j+1)\lc(S_j)>0\).
Since \(S_1=1\), it follows by induction that \(\deg S_j=j-1\)
for \(j\ge1\). In~\eqref{eq:eulerian-degree-recurrence} we
therefore take \(a_j=j+1\) and \(s_j=1\).

\begin{proof}[Proof of Theorem~\ref{thm:largest-intro}]
The middle polynomial \(F_{n-1}\) is nonconstant for \(n\ge3\) in
the families \(A_n,S_n,N_n\) and for \(n\ge2\) in the families
\(B_n,M_n\). To apply Lemma~\ref{lem:free-middle-general}, we
identify the zero \(r\) of \(F_{n-1}\) that minimizes
\(-F_n(r)/F_{n-2}(r)\).

For the ordinary Eulerian polynomials, \eqref{eq:Arec} has the form
\eqref{eq:eulerian-degree-recurrence} with
\(\deg A_j=j-1\) and \(a_j=s_j=1\) for \(j\ge1\).
For type~\(B\), \eqref{eq:Brec} has the same form with
\(\deg B_j=j\), \(a_j=1\), and \(s_j=2\) for \(j\ge0\).
For \(S_j\), equation~\eqref{eq:second-order-recurrence}
has parameters \(a_j=j+1\) and \(s_j=1\) for \(j\ge1\).

In all three cases \(a_j,s_j\ge1\). Starting with
\(A_1=S_1=1\) and \(B_0=1\), Lemma~\ref{lem:largest-recurrence}
gives strict interlacing and simple negative zeros.
Let \(\rho\) be the largest zero of \(F_{n-1}\). At every other zero
\(r\) of \(F_{n-1}\), Lemma~\ref{lem:largest-recurrence} gives
\(-F_n(r)/F_{n-2}(r)>-F_n(\rho)/F_{n-2}(\rho)\).
Thus only \(\rho\) attains the minimum for each of the three Eulerian families.

For the Narayana polynomials, substitute a zero \(r\) of the middle
polynomial into the corresponding recurrence in~\eqref{eq:nar-three}.
Strict interlacing gives \(F_{n-2}(r)\ne0\), so we may divide by \(F_{n-2}(r)\) to obtain
\begin{equation}\label{eq:nar-sharp}
 \begin{aligned}
 -\frac{N_n(r)}{N_{n-2}(r)}&=\frac{n-2}{n+1}(1-r)^2
          &&\bigl(N_{n-1}(r)=0\bigr),\\
 -\frac{M_n(r)}{M_{n-2}(r)}&=\frac{n-1}{n}(1-r)^2
          &&\bigl(M_{n-1}(r)=0\bigr).
 \end{aligned}
\end{equation}
In each case let \(\rho\) be the largest zero of the middle polynomial.
Since \(r\le\rho<0\), we have \((1-r)^2\ge(1-\rho)^2\) with equality
only at \(r=\rho\). Thus for both Narayana families the minimum in
Lemma~\ref{lem:free-middle-general} is attained only at \(\rho\).
For all five families it now follows from that lemma that
\(F_n+bF_{n-1}+cF_{n-2}\) has only real zeros
for every real \(b\) exactly when \(c\le -F_n(\rho)/F_{n-2}(\rho)\).
\end{proof}

\begin{proposition}\label{prop:five-family-full-interlacing}
Let \((F_j)\) be one of the five polynomial sequences
\(A_j\), \(B_j\), \(S_j\), \(N_j\), and \(M_j\).
Take \(n\ge3\) for \(A_j\), \(S_j\), and \(N_j\), and
\(n\ge2\) for \(B_j\) and \(M_j\). Let \(\rho\) be the largest
zero of \(F_{n-1}\). Fix \(c\in\mathbb R\). Then
\(c\le-F_n(\rho)/F_{n-2}(\rho)\) if and only if, for every \(s\ge1\)
and all real numbers \(-1-c\le b_1\le\cdots\le b_s\), the row
\[
 \bigl(F_n(x)+b_jF_{n-1}(x)+cF_{n-2}(x)\bigr)_{j=1}^{s}
\]
is fully interlacing. Every fully interlacing row of this form
satisfies \(b_j\ge-1-c\) for \(1\le j\le s\).
\end{proposition}

\begin{proof}
The denominator \(F_{n-2}(\rho)\) is nonzero because
\(F_{n-2}\prec F_{n-1}\).
Start with the smallest allowed middle coefficient, \(b=-1-c\).
The corresponding polynomial is \(H=F_n-(1+c)F_{n-1}+cF_{n-2}\).
The leading term of \(H\) is that of \(F_n\), and \(H(0)=0\)
since \(F_j(0)=1\). Assume first that \(c\le-F_n(\rho)/F_{n-2}(\rho)\).
By the proof of Theorem~\ref{thm:largest-intro}, only \(\rho\)
attains the minimum of \(-F_n(r)/F_{n-2}(r)\) over the zeros \(r\)
of \(F_{n-1}\). If \(c\) is strictly below this minimum,
Lemma~\ref{lem:free-middle-general} gives \(F_{n-1}\prec H\).

At equality, put
\[
 R(x)=\frac{F_n(x)+cF_{n-2}(x)}{F_{n-1}(x)}.
\]
By~\eqref{eq:three-term-residue} with \(k=0\), the residue at a zero \(r\)
of \(F_{n-1}\) has the sign of \(c+F_n(r)/F_{n-2}(r)\). Since \(\rho\)
uniquely attains the minimum of \(-F_n(r)/F_{n-2}(r)\), the residue
vanishes at \(\rho\) and is negative at every other zero of \(F_{n-1}\).
Thus only \(x-\rho\) cancels. The remaining poles
are the zeros of \(F_{n-1}/(x-\rho)\), all smaller than \(\rho\).
Write \(\gamma_r<0\) for the residue at each remaining pole \(r\).
The polynomial part of \(R\) is linear, with slope
\(\lc(F_n)/\lc(F_{n-1})>0\). Thus away from the poles,
\[
 R'(x)=\frac{\lc(F_n)}{\lc(F_{n-1})}
       -\sum_{\substack{F_{n-1}(r)=0\\r\ne\rho}}
          \frac{\gamma_r}{(x-r)^2}>0.
\]
At each pole \(r\), the term \(\gamma_r/(x-r)\) gives the limits
\(+\infty\) from the left and \(-\infty\) from the right. The positive
slope gives \(R(x)\to\pm\infty\) as \(x\to\pm\infty\). Hence
\(R(x)=1+c\) has exactly one simple solution in each interval cut
out by the poles. If there are no poles, this interval is \(\mathbb R\).
Since \(H=F_{n-1}(R-1-c)\), these solutions are the zeros of
\(H/(x-\rho)\). Thus \(F_{n-1}/(x-\rho)\prec H/(x-\rho)\).

Both \(\rho\) and \(0\) lie in the rightmost interval. There \(R\)
is strictly increasing, so \(R(\rho)<R(0)=1+c\). The rightmost
solution is \(0\). All other solutions lie below \(\rho\).
At \(\rho\), the factor \(H/(x-\rho)=(F_{n-1}/(x-\rho))(R-1-c)\)
is nonzero. Thus restoring \(x-\rho\) gives \(H\) a simple zero at
\(\rho\), immediately before \(0\). Since \(\rho\) is the largest zero
of \(F_{n-1}\), we have \(F_{n-1}\preceq H\), with \(\rho\) as the
only common zero.

In both cases the interlacing puts every zero of \(H\) other than
its largest at or below \(\rho<0\). Since \(H(0)=0\), its largest
zero is \(0\). Both \(H\) and \(F_{n-1}\) have positive leading
coefficients and nonpositive zeros, so their coefficients are
nonnegative. In our notation, the interlacing relation in Example~3.7(b)
of~\cite{AthanasiadisWagner2026} is \(\preceq\), so it allows the
common zero \(\rho\) at equality. The example gives the fully interlacing column
\((F_{n-1},H)^{\mathsf T}\).
Reflecting it across the reverse diagonal gives the fully
interlacing row \((H,F_{n-1})\)
by~\cite[Theorem~4.1]{AthanasiadisWagner2026}.

The constant matrix
\[
 C=
 \begin{pmatrix}
 1&\cdots&1\\
 1+c+b_1&\cdots&1+c+b_s
 \end{pmatrix}
\]
has nonnegative entries, and its minor on columns \(i<j\)
is \(b_j-b_i\ge0\). Thus \(C\) is totally nonnegative and
hence fully interlacing by~\cite[Example~3.7(a)]{AthanasiadisWagner2026}.
By~\cite[Theorem~4.3(b)]{AthanasiadisWagner2026}, the product
\((H,F_{n-1})C\) is fully interlacing. Its \(j\)th entry is
\(H+(1+c+b_j)F_{n-1}=F_n+b_jF_{n-1}+cF_{n-2}\).

Conversely, take \(s=2\), \(b_1=-1-c\), and \(b_2=-c\).
If the row is fully interlacing for every allowed choice of
middle coefficients, then \((H,H+F_{n-1})\) is fully interlacing.
Thus \(H+F_{n-1}\preceq H\). By Obreschkoff's theorem
\cite[Theorem~2.4]{Branden2006}, every nonzero real linear combination
of \(H\) and \(H+F_{n-1}\) has only real zeros. In particular,
\(H+(1+c+b)F_{n-1}=F_n+bF_{n-1}+cF_{n-2}\) has only real
zeros for every real \(b\). Theorem~\ref{thm:largest-intro} gives
the upper bound on \(c\).

Finally, every entry of a fully interlacing row has nonnegative
coefficients. The \(j\)th entry has constant term \(1+b_j+c\), so
\(b_j\ge-1-c\).
\end{proof}

The next example attains the upper bound on \(c\) in
Proposition~\ref{prop:five-family-full-interlacing} by weighting
noncrossing partitions. Reversing the order of the polynomials in
the resulting row can destroy full interlacing.

\Needspace{6\baselineskip}
\begin{example}
The ordinary Narayana polynomial \(N_j(x)\) counts noncrossing
partitions of \([j]\), with \(x^{k-1}\) for a partition into \(k\)
blocks. For a partition of \([3]\), assign weight \(t\ge0\) when
\(2\) and \(3\) lie in the same block. Assign weight one if
neither adjacent pair has both letters in the same block.
In the remaining case, assign weight zero.

For either adjacent pair, merging its elements when they are in
one block gives a bijection with the partitions of \([2]\) and
preserves the number of blocks. Only the one-block partition has
both adjacent pairs in the same block. Thus the partitions with
neither pair in a single block contribute \(N_3-2N_2+N_1\) by
inclusion--exclusion, and those with \(2\) and \(3\) in the same
block contribute \(tN_2\). The weighted polynomial is
\[
 N_3(x)+(t-2)N_2(x)+N_1(x)=(1+x)(x+t).
\]
The largest zero of \(N_2(x)=1+x\) is \(\rho=-1\), and
\(-N_3(\rho)/N_1(\rho)=1\). Thus \(c=1\) attains the bound in
Proposition~\ref{prop:five-family-full-interlacing} at \(n=3\).
Since \(b=t-2\ge-2=-1-c\), the row
\(\bigl((1+x)(x+t_j)\bigr)_{j=1}^{s}\) is fully interlacing for
\(0\le t_1\le\cdots\le t_s\).
Taking \(t_1=0\) and \(t_2=1\) gives
\(\bigl(x(1+x),(1+x)^2\bigr)\), with a double zero at \(-1\)
in the second entry. The order of the two entries matters.
In the infinite matrix defining full interlacing,
the submatrix in rows \(0,1\) and columns \(2,3\) is
\(\left(\begin{smallmatrix}1&2\\0&1\end{smallmatrix}\right)\),
with determinant \(1\). Reversing the row interchanges these
two columns and gives determinant \(-1\).

For \(c>1\), the allowed coefficient \(b=-1\) gives
\(N_3-N_2+cN_1=x^2+2x+c\), whose discriminant is \(4(1-c)<0\).
\end{example}

Before finding a bound on \(c\) that works for all \(n\), we describe
the zeros of the three-term combination at a fixed index.
Fix one of the five families and an index \(n\) for
which \(F_{n-1}\) is nonconstant, and let \(\rho\) be its largest zero.
If \(c<-F_n(\rho)/F_{n-2}(\rho)\), Lemma~\ref{lem:free-middle-general}
gives one simple zero of \(F_n+bF_{n-1}+cF_{n-2}\) in each interval
cut out by the zeros of \(F_{n-1}\), for every real \(b\).
When \(c=-F_n(\rho)/F_{n-2}(\rho)\), only the factor \(x-\rho\)
cancels in \((F_n+cF_{n-2})/F_{n-1}\). The derivative calculation
after cancellation in
\hyperlink{proof:free-middle-multiplicity}{the proof of Lemma~\ref*{lem:free-middle-general}}
with \(k=0\) shows that \(\rho\) is the only possible multiple zero.
It is double exactly when
\[
 b=-\frac{F_n'(\rho)+cF_{n-2}'(\rho)}{F_{n-1}'(\rho)}.
\]

\subsection{Uniform bounds for Eulerian polynomials}

By Theorem~\ref{thm:largest-intro}, the bound at a fixed index \(n\)
is \(-F_n(\rho)/F_{n-2}(\rho)\), where \(\rho\) is the largest zero
of \(F_{n-1}\). The bound valid for all \(n\) is the infimum of these
bounds, taken over the indices for which \(F_{n-1}\) is nonconstant.
For the three Eulerian families, low-degree factorizations let us
compare the value at \(n=4\) with those at
\(n=3,5\) for \(A_n\), \(n=2,3\) for \(B_n\), and \(n=3\) for \(S_n\).

\smallskip\noindent\emph{Ordinary Eulerian polynomials.}
For \(n=3\), substitute the zero \(-1\) of \(A_2=1+x\). This gives
\(-A_3(-1)/A_1(-1)=2\). For the next two indices use
\[
 \begin{aligned}
 A_3&=x^2+4x+1,& A_4&=(1+x)(A_3+6x),\\
 A_4&=(1+x)(x^2+10x+1),&
 A_5-16xA_3&=(1+x^2)(x^2+10x+1).
 \end{aligned}
\]
At the largest zero \(\rho=-2+\sqrt3\) of \(A_3\), the first line
gives \(-A_4(\rho)/A_2(\rho)=-6\rho\). At the largest zero
\(\rho=-5+2\sqrt6\) of \(A_4\), the second line gives
\(-A_5(\rho)/A_3(\rho)=-16\rho\). Both denominators are nonzero by strict interlacing. The bounds for
\(n=4,5\) are therefore
$6/(2+\sqrt3)$ and $16/(5+2\sqrt6)$,
respectively. The bound at \(n=4\) is less than \(2\), the bound at
\(n=3\). To show that it is also less than the bound at \(n=5\),
cross-multiply and divide by two. It remains to check
\(6\sqrt6<1+8\sqrt3\). Squaring gives \(23<16\sqrt3\),
which follows from \(23^2<16^2\cdot3\).

\smallskip\noindent\emph{Type \(B\) Eulerian polynomials.}
From \(B_0=1\), \(B_1=1+x\), and \(B_2=1+6x+x^2\), we have
\(-B_2(-1)/B_0(-1)=4\). The next two factorizations give the ratios
at the largest zeros of \(B_2\) and \(B_3\):
\[
 \begin{aligned}
 B_3&=(1+x)(B_2+16x)=(1+x)(x^2+22x+1),\\
 B_4-60xB_2&=(x^2+22x+1)(x^2-6x+1).
 \end{aligned}
\]
At the largest zero \(\rho=-3+2\sqrt2\) of \(B_2\), the first identity
gives \(-B_3(\rho)/B_1(\rho)=-16\rho\). At the largest zero
\(\rho=-11+2\sqrt{30}\) of \(B_3\), the second gives \(-B_4(\rho)/B_2(\rho)=-60\rho\).
Thus the bound at \(n=4\) is smaller than the one at \(n=3\)
precisely when
$60/(11+2\sqrt{30})<16/(3+2\sqrt2)$.
Cross-multiplying and dividing by four gives
\(1+30\sqrt2<8\sqrt{30}\). Both sides are positive, so it suffices
to square and check
\(60\sqrt2<119\), which follows from \(\sqrt2<3/2\).
Also \(16/(3+2\sqrt2)<4\), since \(3+2\sqrt2>4\).

\smallskip\noindent\emph{Second-order Eulerian polynomials.}
From \(S_1=1\), \(S_2=1+2x\), and \(S_3=1+8x+6x^2\), we have
\(-S_3(-1/2)/S_1(-1/2)=3/2\). At the larger zero
\(\rho=(-4+\sqrt{10})/6\) of \(S_3\), the identity
$S_4-10xS_2=(1+4x)S_3$
gives \(-S_4(\rho)/S_2(\rho)=-10\rho=5(4-\sqrt{10})/3\).
This is less than the bound \(3/2\) for \(n=3\), since
\(31^2<10^2\cdot10\) gives \(31<10\sqrt{10}\).

For each Eulerian family, the smallest bound among the indices checked
occurs at \(n=4\). It remains to check \(n\ge6\) for \(A_n\) and
\(n\ge5\) for \(B_n\) and \(S_n\). For these indices, keep only the
constant term and the square corresponding to the largest zero of
\(F_{n-2}\) in~\eqref{eq:eulerian-degree-squares}.

In each of these three families the recurrence has the form
\[
 F_{j+1}=(1+\lambda_jx)F_j+sx(1-x)F_j',
 \qquad F_j(0)=1.
\]
In the notation of~\eqref{eq:eulerian-degree-recurrence}, \(s=s_j\)
is independent of \(j\) and \(\lambda_j=a_j+s\deg F_j\). Since
\(a_j>0\) and \(F_j\) is nonconstant at the indices considered below,
\(\lambda_j>s\). Write \(-1/u_j\) for the largest zero of \(F_j\).
By strict interlacing, \(0<u_j<u_{j+1}\). Set \(n=j+2\), so the
largest zero of \(F_{n-1}\) is \(\rho=-1/u_{j+1}\).
The fraction in~\eqref{eq:eulerian-degree-squares} corresponding
to the largest zero \(-1/u_j\) of \(F_j\) is
\[
 \frac{\rho(1+1/u_j)}{\rho+1/u_j}
 =-\frac{u_j+1}{u_{j+1}-u_j}.
\]
Since \(s_{n-1}=s_{n-2}=s\), multiply~\eqref{eq:eulerian-degree-squares}
by \(s\) and retain the constant term and this square.
This gives
\begin{equation}\label{eq:eulerian-uniform-square}
 -\frac{F_{j+2}(\rho)}{F_j(\rho)}\ge s+s^2\left(\frac{u_j+1}{u_{j+1}-u_j}\right)^2.
\end{equation}

For each family, choose \(t>0\) so that \(s+s^2t^2\) exceeds its
bound at \(n=4\). By~\eqref{eq:eulerian-uniform-square}, it suffices
to prove \((u_j+1)/(u_{j+1}-u_j)>t\) at the remaining indices.
We derive a lower bound on \(u_j\) that implies this inequality
and then prove that bound by induction.
Factor \(F_j(x)=\prod_w(1+wx)\). The parameters \(w\) are positive,
with largest value \(u_j\). Their sum is \(\sum_w w=F_j'(0)\).
In particular, \(F_j(-1/u_{j+1})>0\).
Evaluate the recurrence at the largest zero \(-1/u_{j+1}\) of
\(F_{j+1}\) and divide by \(F_j(-1/u_{j+1})\). Using
\(F_j'/F_j=\sum_w w/(1+wx)\) gives
\[
 0=1-\frac{\lambda_j}{u_{j+1}}
      -\frac{s(u_{j+1}+1)}{u_{j+1}}\sum_w\frac{w}{u_{j+1}-w}.
\]
Since \(w\le u_j<u_{j+1}\), we have
\[
 \frac{u_{j+1}+1}{u_{j+1}-w}
 \le \frac{u_{j+1}+1}{u_{j+1}-u_j}
 =1+\frac{u_j+1}{u_{j+1}-u_j}.
\]
Rearranging the recurrence at the zero and using \(\sum_w w=F_j'(0)\) gives
\[
 u_{j+1}=\lambda_j+s\sum_w w\frac{u_{j+1}+1}{u_{j+1}-w}
 \le \lambda_j+sF_j'(0)\left(1+\frac{u_j+1}{u_{j+1}-u_j}\right).
\]
If \((u_j+1)/(u_{j+1}-u_j)\le t\), then
\[
 (1+1/t)u_j+1/t\le u_{j+1}\le\lambda_j+s(1+t)F_j'(0).
\]
The leftmost expression exceeds the rightmost one precisely when
\begin{equation}\label{eq:eulerian-gap-test}
 u_j>stF_j'(0)+\frac{t\lambda_j-1}{1+t}.
\end{equation}
Thus~\eqref{eq:eulerian-gap-test} implies
\((u_j+1)/(u_{j+1}-u_j)>t\).

Write \(h_j=stF_j'(0)+(t\lambda_j-1)/(1+t)\). For the induction step
assume \(u_j>h_j\). Since \(s+1>0\), it suffices to show that
\[
 u_{j+1}-(s+1)u_j\ge h_{j+1}-(s+1)h_j.
\]
Comparing coefficients of \(x\) in the recurrence gives
\(F_{j+1}'(0)=(s+1)F_j'(0)+\lambda_j\). Substitution into the definition
of \(h_{j+1}\) gives
\[
 h_{j+1}-(s+1)h_j=st\lambda_j+
 \frac{t\lambda_{j+1}-1-(s+1)(t\lambda_j-1)}{1+t}
 =\frac{st^2\lambda_j+t(\lambda_{j+1}-\lambda_j)+s}{1+t}.
\]
To bound \(u_{j+1}-(s+1)u_j\) below, retain the term \(w=u_j\)
in the sum for \(u_{j+1}\). Since the other terms are positive,
\[
 u_{j+1}\ge\lambda_j+su_j\frac{u_{j+1}+1}{u_{j+1}-u_j}.
\]
Multiplying by \(u_{j+1}-u_j>0\) and collecting the terms in \(u_ju_{j+1}\) gives
\[
 u_{j+1}^2-(s+1)u_ju_{j+1}\ge\lambda_j u_{j+1}-(\lambda_j-s)u_j.
\]
Dividing by \(u_{j+1}>0\) gives
\[
 u_{j+1}-(s+1)u_j\ge\lambda_j-(\lambda_j-s)\frac{u_j}{u_{j+1}}>s,
\]
where the strict inequality uses \(\lambda_j>s\) and \(u_j/u_{j+1}<1\).
In particular \(u_j/u_{j+1}<1/(s+1)\). Since \(\lambda_j-s>0\), substituting this upper bound for
\(u_j/u_{j+1}\) gives
\[
 u_{j+1}-(s+1)u_j>\lambda_j-\frac{\lambda_j-s}{s+1}
            =\frac{s(\lambda_j+1)}{s+1}.
\]
Thus the induction step follows whenever
\begin{equation}\label{eq:eulerian-gap-induction}
 \frac{s(\lambda_j+1)}{s+1}
 \ge \frac{st^2\lambda_j+t(\lambda_{j+1}-\lambda_j)+s}{1+t}.
\end{equation}

The following table gives \(t\) and the initial index \(j_0\) for
each family. Since \(n=j+2\), these indices correspond to \(n=6\)
for \(A\) and \(n=5\) for \(B\) and \(S\). First we check that
\(s+s^2t^2\) exceeds the bound at \(n=4\). For the induction, we
verify~\eqref{eq:eulerian-gap-test} at \(j_0\) and
\eqref{eq:eulerian-gap-induction} for \(j\ge j_0\).
\[
 \begin{array}{c|c|c|c|c|c|c}
 F & j_0&s&\lambda_j&t&u_{j_0}&F_{j_0}'(0)\\ \hline
 A&4&1&j&11/14&5+2\sqrt6&11\\
 B&3&2&2j+1&3/7&11+2\sqrt{30}&23\\
 S&3&1&2j&2/3&4+\sqrt{10}&8
 \end{array}
\]
The three rows use \(A_4,B_3,S_3\), respectively. For each of these
polynomials, the last two columns give the negative reciprocal of its
largest zero and its coefficient of \(x\). Substituting the listed values of \(s,t\), the comparisons
of \(s+s^2t^2\) with the three bounds at \(n=4\) become
\[
 \frac{317}{196}>\frac6{2+\sqrt3},\qquad
 \frac{134}{49}>\frac{60}{11+2\sqrt{30}},\qquad
 \frac{13}{9}>\frac{5(4-\sqrt{10})}{3}.
\]
Cross multiplication reduces these inequalities to
\(317\sqrt3>542\), \(134\sqrt{30}>733\), and \(15\sqrt{10}>47\).
The squares of the left sides exceed the squares of the right sides
by \(7703,1391,41\), respectively.

At the initial index, the three instances of~\eqref{eq:eulerian-gap-test}
are
\[
 \begin{aligned}
 A:\qquad 5+2\sqrt6
   &> \frac{11}{14}\cdot11+\frac{11\cdot4-14}{25}
    = \frac{689}{70},\\
 B:\qquad 11+2\sqrt{30}
   &> \frac67\cdot23+\frac{3\cdot7-7}{10}
    = \frac{739}{35},\\
 S:\qquad 4+\sqrt{10}
   &> \frac23\cdot8+\frac{2\cdot6-3}{5}
    = \frac{107}{15}.
 \end{aligned}
\]
For the first inequality, \(6\cdot7^2>17^2\) gives
\(5+2\sqrt6>69/7=690/70\). For the second,
\(30\cdot4^2>21^2\) gives \(11+2\sqrt{30}>43/2>739/35\).
The third follows from \(15^2\cdot10=2250>2209=47^2\).

Finally, subtract the right side of~\eqref{eq:eulerian-gap-induction}
from the left. In the three rows the differences are
\[
 \begin{aligned}
 A:\quad &\frac{j+1}{2}
    -\frac{121j+154+196}{350}
     =\frac{54j-175}{350}>0 &&(j\ge4),\\
 B:\quad &\frac{2(2j+2)}3
    -\frac{18(2j+1)+42+98}{70}
     =\frac{86j-97}{105}>0 &&(j\ge3),\\
 S:\quad &\frac{2j+1}{2}
    -\frac{8j+12+9}{15}
     =\frac{14j-27}{30}>0 &&(j\ge3).
 \end{aligned}
\]
Consequently~\eqref{eq:eulerian-gap-test} holds for every
\(j\ge j_0\) in each row. Thus \((u_j+1)/(u_{j+1}-u_j)>t\),
and~\eqref{eq:eulerian-uniform-square} puts the bound at index
\(j+2\) strictly above \(s+s^2t^2\), hence above the bound at \(n=4\).
Together with the small indices, this proves that in all three
families the unique minimum occurs at \(n=4\). Rationalizing the
denominators gives \(6/(2+\sqrt3)=12-6\sqrt3\) and
\(60/(11+2\sqrt{30})=60(\sqrt6-\sqrt5)^2\).

\subsection{Uniform bounds for Narayana polynomials}

For \(N_n\) the bound at index \(n\) is \((n-2)(1-\rho)^2/(n+1)\),
and for \(M_n\) it is \((n-1)(1-\rho)^2/n\), by~\eqref{eq:nar-sharp}.
In each formula, \(\rho\) is the largest zero of the preceding
polynomial in the same family. We will show that the minima occur
at \(n=5\) for \(N_n\) and \(n=6\) for \(M_n\). We compute those
bounds exactly. For the remaining indices, we use the first three
nonconstant coefficients to bound \(|\rho|\) from below.

Take \(F_{n-1}=N_{n-1}\) or \(M_{n-1}\). Its constant term is one,
and every zero \(r\) is negative. Thus
\(F_{n-1}(x)=\prod_r(1+x/|r|)\) and
\[
 \frac{F_{n-1}'(-x)}{F_{n-1}(-x)}
 =\sum_r\frac1{|r|-x}.
\]
For \(j\ge1\), the coefficient of \(x^{j-1}\) in the expansion at
zero is \(\sum_r |r|^{-j}\). All summands are positive, and the one
corresponding to the largest zero is \(|\rho|^{-j}\).
To compute the first three power sums, multiply the quotient identity
by \(F_{n-1}(-x)\) and compare coefficients of \(1,x,x^2\).
The corresponding coefficients of \(F_{n-1}(-x)\) are
\(1,-F_{n-1}'(0),F_{n-1}''(0)/2\). This gives
\[
 \begin{aligned}
 \sum_r |r|^{-1}&=F_{n-1}'(0),\\
 \sum_r |r|^{-2}
 &=F_{n-1}'(0)\sum_r |r|^{-1}-F_{n-1}''(0),\\
 \sum_r |r|^{-3}
 &=F_{n-1}'(0)\sum_r |r|^{-2}
   -\tfrac12F_{n-1}''(0)\sum_r |r|^{-1}+\tfrac12F_{n-1}'''(0).
 \end{aligned}
\]
Substituting the lower power sums into the last two identities and
using \(|\rho|^{-j}\le\sum_r |r|^{-j}\) gives
\begin{equation}\label{eq:narayana-root-estimates}
 \begin{aligned}
 |\rho|^{-1}&\le F_{n-1}'(0),\\
 |\rho|^{-2}&\le F_{n-1}'(0)^2-F_{n-1}''(0),\\
 |\rho|^{-3}&\le F_{n-1}'(0)^3
          -\tfrac32 F_{n-1}'(0)F_{n-1}''(0)+\tfrac12F_{n-1}'''(0).
 \end{aligned}
\end{equation}

\smallskip\noindent\emph{Ordinary Narayana polynomials.}
For \(n=3,4,5\), the polynomials \(N_{n-1}(x)\) are \(1+x\),
\(1+3x+x^2\), and \((1+x)(1+5x+x^2)\), respectively.
Their largest zeros are \(-1\), \((-3+\sqrt5)/2\), and
\((-5+\sqrt{21})/2\), respectively. Substitution into
\((n-2)(1-\rho)^2/(n+1)\) gives
\[
 \begin{aligned}
 n=3:&\quad \tfrac14(1+1)^2=1,\\
 n=4:&\quad \tfrac25\left(\frac{5-\sqrt5}{2}\right)^2=3-\sqrt5,\\
 n=5:&\quad \tfrac12\left(\frac{7-\sqrt{21}}2\right)^2
                  =\frac{35-7\sqrt{21}}4=\frac7{5+\sqrt{21}} .
 \end{aligned}
\]
To compare all other bounds with the one for \(n=5\), observe that
\(32^2=1024<1029=7^2\cdot21\). Thus
$7/(5+\sqrt{21})<7/(5+32/7)=49/67$.
The bounds for \(n=3,4\) are \(1\) and \(3-\sqrt5\). Both exceed
\(3/4>49/67\), since \(5<(9/4)^2\).
We show that every bound for \(n\ge6\) also exceeds \(49/67\).

For \(n\ge9\), the linear coefficient of \(N_{n-1}\) is
\(\binom{n-1}1\binom{n-1}2/(n-1)=(n-1)(n-2)/2\). The first
line of~\eqref{eq:narayana-root-estimates} gives
\(|\rho|\ge2/[(n-1)(n-2)]\). Dropping the positive term
\(|\rho|^2\) from \((1-\rho)^2=1+2|\rho|+|\rho|^2\) gives
\[
 \begin{aligned}
 \frac{n-2}{n+1}(1-\rho)^2
 &>\frac{n-2}{n+1}\left(1+\frac4{(n-1)(n-2)}\right)\\
 &=\frac34+\frac{(n-3)(n-9)}{4(n-1)(n+1)}\ge\frac34 .
 \end{aligned}
\]

For \(n=6,7,8\) we use the third, second, and first power sums in
\eqref{eq:narayana-root-estimates}, respectively. The binomial
coefficients in~\eqref{eq:narayana-coefficients} give
\[
 \begin{aligned}
 N_5:&\quad |\rho|^{-3}\le10^3-3\cdot10\cdot20+3\cdot10
                         =430<(38/5)^3,\\
 N_6:&\quad |\rho|^{-2}\le15^2-2\cdot50=125<12^2,\\
 N_7:&\quad |\rho|^{-1}\le21.
 \end{aligned}
\]
In each row \(\rho\) is the largest zero of the indicated polynomial.
The resulting lower bounds for \(|\rho|\) are \(5/38,1/12,1/21\),
respectively. Substituting each into the ratio in~\eqref{eq:nar-sharp}
gives the three bounds in the middle column. The positive differences in the last
column show that all three ratios exceed \(49/67\).
\[
 \begin{array}{c|c|c}
 n & \text{lower bound for }\dfrac{n-2}{n+1}(1-\rho)^2
       &\text{difference of cross products}\\ \hline
 6 &\tfrac47(1+5/38)^2=1849/2527 &67\cdot1849-49\cdot2527=60\\[2pt]
 7 &\tfrac58(1+1/12)^2=845/1152 &67\cdot845-49\cdot1152=167\\[2pt]
 8 &\tfrac69(1+1/21)^2=968/1323 &67\cdot968-49\cdot1323=29
 \end{array}
\]
Thus the unique minimum is \((35-7\sqrt{21})/4\), at \(n=5\).

\smallskip\noindent\emph{Type \(B\) Narayana polynomials.}
For \(n=6\) we first evaluate \(\frac56(1-\rho)^2\) exactly.
The middle polynomial is
\(M_5=(1+x)(x^4+24x^3+76x^2+24x+1)\).
Set \(x=-t\) in the quartic. Division by \(t^2\) and the identity
\(t^2+t^{-2}=(t+t^{-1})^2-2\) give
$(t+t^{-1})^2-24(t+t^{-1})+74=0$.
Both values \(12\pm\sqrt{70}\) of \(t+t^{-1}\) exceed two.
Each gives two positive reciprocal values of \(t\), one below \(1\)
and one above. Since \(t+t^{-1}\) decreases on \((0,1)\), the
smallest of these four values corresponds to \(12+\sqrt{70}\).
The linear factor gives \(t=1\), so this smallest \(t\) is
\(|\rho|\), where \(\rho\) is the largest zero of \(M_5\). Thus
\(|\rho|+|\rho|^{-1}=12+\sqrt{70}\).
It follows that
\[
 \left(\frac{1+\rho}{1-\rho}\right)^2
 =\frac{|\rho|+|\rho|^{-1}-2}{|\rho|+|\rho|^{-1}+2}
 =\frac{10+\sqrt{70}}{14+\sqrt{70}}.
\]
Multiplying numerator and denominator by \(14-\sqrt{70}\) gives
\((70+4\sqrt{70})/126=(35+2\sqrt{70})/63\).
Since \(0<|\rho|<1\), the positive square root gives
\[
 \frac{2}{1-\rho}
 =1+\frac{1+\rho}{1-\rho}
 =1+\frac{\sqrt{35+2\sqrt{70}}}{3\sqrt7}.
\]
Solving for \(1-\rho\) gives the bound for \(n=6\),
\[
 \frac56(1-\rho)^2=\frac56
 \left(\frac{6\sqrt7}{3\sqrt7+\sqrt{35+2\sqrt{70}}}\right)^2
 =\frac{210}{\bigl(3\sqrt7+\sqrt{35+2\sqrt{70}}\bigr)^2}.
\]

For \(n\ge8\), let \(\rho\) be the largest zero of \(M_{n-1}\).
The coefficients of \(x\) and \(x^2\) in \(M_{n-1}\) are \((n-1)^2\) and \(\binom{n-1}{2}^2\).
The second line of \eqref{eq:narayana-root-estimates} therefore
gives
\[
 |\rho|^{-2}\le(n-1)^4-2\binom{n-1}2^2
 =\tfrac12(n-1)^2(n^2-2)<\tfrac12(n-1)^2n^2.
\]
Thus \(|\rho|>\sqrt2/[n(n-1)]\).
Substituting this into \((n-1)(1+|\rho|)^2/n\) and
discarding the positive square gives
\[
 \frac{n-1}{n}(1-\rho)^2>1-\frac1n+\frac{2\sqrt2}{n^2}
        >1-\frac1n+\frac{11}{4n^2},
\]
where \(2\sqrt2>11/4\) follows by squaring. For every \(n\ge8\), the
lower bound \(1-1/n+11/(4n^2)\) is at least \(235/256\), as the
following factorization shows.
\[
 1-\frac1n+\frac{11}{4n^2}-\frac{235}{256}
 =\frac{(n-8)(21n-88)}{256n^2}\ge0\qquad(n\ge8).
\]
Thus \((n-1)(1-\rho)^2/n>235/256\) for every \(n\ge8\).
It remains to show that the bound at \(n=6\) is below \(235/256\)
and the bounds at \(n=2,3,4,5,7\) are above it.

For \(n=6\), \(\rho\) is the largest zero of \(M_5\). Since \(|\rho|+|\rho|^{-1}=12+\sqrt{70}\) and \(0<|\rho|<1\),
we have \(|\rho|^{-1}=12+\sqrt{70}-|\rho|>12\).
Using \(|\rho|<1/12\) and \(\sqrt{70}>25/3\) in
\(|\rho|^{-1}=12+\sqrt{70}-|\rho|\) sharpens the bound.
The inequality for \(\sqrt{70}\) follows from \(630>625\). Thus
\[
 |\rho|^{-1}>12+\frac{25}{3}-\frac1{12}=\frac{81}{4},
 \qquad
 \tfrac56(1-\rho)^2<\tfrac56(1+4/81)^2=\frac{36125}{39366}
       <\frac{235}{256}.
\]
The last inequality follows from
\(235\cdot39366-256\cdot36125=3010\).

For \(n=2\) we have \(M_1=1+x\), so the bound is
\(\tfrac12(1+1)^2=2\). For \(n=3,4,5\), the first two lines of
\eqref{eq:narayana-root-estimates} give
\[
 \begin{aligned}
 M_2:&\quad |\rho|^{-1}\le4,\\
 M_3:&\quad |\rho|^{-1}\le9,\\
 M_4:&\quad |\rho|^{-2}\le16^2-2\cdot36=184<14^2 .
 \end{aligned}
\]
With \(\rho\) the largest zero in each row, the respective lower
bounds for \(|\rho|\) are \(1/4,1/9,1/14\). The corresponding
lower bounds for \((n-1)(1-\rho)^2/n\) are
\[
 \tfrac23(1+1/4)^2=\tfrac{25}{24},\qquad
 \tfrac34(1+1/9)^2=\tfrac{25}{27},\qquad
 \tfrac45(1+1/14)^2=\tfrac{45}{49}.
\]
They exceed \(235/256\), since
\[
 \frac{25}{24}>\frac{25}{27}>\frac{45}{49}>\frac{235}{256},
 \qquad 25\cdot49-45\cdot27=10,\quad 45\cdot256-49\cdot235=5.
\]
For the remaining index \(n=7\), the coefficients of \(x,x^2,x^3\)
in \(M_6\) are \(36,225,400\). The third power sum
in~\eqref{eq:narayana-root-estimates} therefore gives
\[
 |\rho|^{-3}\le36^3-3\cdot36\cdot225+3\cdot400
 =23556=(86/3)^3-44/27,
\]
and so \(|\rho|>3/86\). We obtain
\[
 \frac67(1-\rho)^2>\frac67(89/86)^2>\frac{235}{256},
 \qquad 6\cdot89^2\cdot256-7\cdot86^2\cdot235=236.
\]
Thus the bound for \(M_n\) attains its minimum only at \(n=6\).

The three Eulerian minima and the two Narayana minima give the five
entries in Table~\ref{tab:main-bounds}.
Fix one family, and for each \(n\) for which \(F_{n-1}\) is nonconstant
let \(\rho\) denote the largest zero of \(F_{n-1}\). If \(c\) is at or below
the tabulated bound for this family, then \(c\le-F_n(\rho)/F_{n-2}(\rho)\) for
every such \(n\). By Theorem~\ref{thm:largest-intro},
\(F_n+bF_{n-1}+cF_{n-2}\) has only real zeros for these
indices and all \(b\in\mathbb R\). Conversely, if \(c\) exceeds the tabulated
bound, take the unique index \(n\) at which the minimum is attained.
At this index \(c>-F_n(\rho)/F_{n-2}(\rho)\), so the same theorem
gives a real \(b\) for which \(F_n+bF_{n-1}+cF_{n-2}\) has a nonreal zero.
By Proposition~\ref{prop:five-family-full-interlacing}, the same
bounds are sharp for full interlacing of all rows with
\(-1-c\le b_1\le\cdots\le b_s\), for every \(n\) with
\(F_{n-1}\) nonconstant.

We can now return to the permutations in~\eqref{inexAn}.
For fixed \(1\le i<j\le n\), the excedance polynomial of permutations
with \(\pi(i)\ne i\) and \(\pi(j)\ne j\) is
\(A_n-2A_{n-1}+A_{n-2}\). The coefficient of \(A_{n-2}\) is \(c=1\), and
\(1<12-6\sqrt3\) since \(11^2>6^2\cdot3\). Thus for every
\(n\ge3\), Lemma~\ref{lem:free-middle-general} gives
\(A_{n-1}\prec A_n-2A_{n-1}+A_{n-2}\).
Every nontrivial cycle has an excedance at its least letter.
Thus the identity, which fixes both \(i\) and \(j\), is the only
permutation with no excedances. The restricted polynomial therefore
has zero constant term. A permutation in the restricted class with
one excedance has a single nontrivial cycle containing both \(i\)
and \(j\). Write that cycle starting with its least letter.
The remaining letters must occur in decreasing order to avoid
a second excedance. Thus each subset
containing \(i\) and \(j\) determines exactly one such permutation.
The coefficient of \(x\) is \(2^{n-2}\), so \(0\) is a simple zero.
By strict interlacing, all other zeros lie to the left of the largest zero
of \(A_{n-1}\) and are negative. For \(n=4\) the polynomial is
\(x(x^2+9x+4)\).

\subsection{A combination of Narayana polynomials}
\label{sec:narayana-interlacing}

For \(n\ge4\), the Fomin--Zelevinsky cluster complex of type
\(D_n\) has \(h\)-polynomial \(M_n-nxN_{n-1}\). With the coefficient of
\(M_n\) set to one, Br\"and\'en's Theorem~7.1~\cite{Branden2006}
asserts that all zeros of \(M_n+\beta n xN_{n-1}\) are real and
simple for \(\beta>-2\). Chen, Yang, and
Zhang~\cite[Theorem~1.1 and footnote~1]{CYZ2016pre} found a
counterexample at \(n=4\) and \(\beta=-19/10\), where
\(M_4-\tfrac{38}{5}xN_3\) has a nonreal conjugate pair.
They proved that \(M_n+\beta n xN_{n-1}\) has only real zeros for
\(\beta\ge-1\). We determine its full real-rootedness range for each \(n\ge2\).

\begin{theorem}\label{thm:narayana-branden-range}
For \(n\ge2\) and real \(\beta\), the polynomial
\(M_n(x)+\beta n\,xN_{n-1}(x)\) has only real zeros if and only if
\[
 \beta\le-4+\frac2n
 \qquad\text{or}\qquad
 \beta\ge-\frac{n-1}{\lfloor n/2\rfloor}.
\]
\end{theorem}

The polynomial \(M_n+\beta n xN_{n-1}\) is symmetric in degree \(n\).
To apply Lemma~\ref{lem:reciprocal-pair}, we verify
\(M_n(r)/N_{n-1}'(r)>0\) at every zero \(r<-1\) of \(N_{n-1}\).
The lemma then identifies the only reciprocal pair that can become
nonreal. The zeros in that pair can coincide only at \(x=\pm1\).
The endpoints of the real-rootedness range follow by evaluating
\(M_n+\beta n xN_{n-1}\) at \(x=1\) and \(x=-1\). For odd \(n\),
use the quotient by \(1+x\) at \(x=-1\).

\begin{lemma}\label{lem:mixed-narayana-residue-sign}
For \(n\ge2\), the polynomial \(N_{n-1}(x)\) strictly interlaces
\(M_{n-1}(x)\). At each zero \(r<-1\) of \(N_{n-1}(x)\), we have
\(M_n(r)/N_{n-1}'(r)>0\).
\end{lemma}
\begin{proof}
From~\eqref{eq:narayana-coefficients} we obtain
\begin{equation}\label{eq:mixed-narayana-identities}
 \begin{aligned}
 M_n&=(1+x)M_{n-1}+2(n-1)xN_{n-1},\\
 n(n-1)N_{n-1}&=(1-x)M_{n-1}'+(n-1)M_{n-1}.
 \end{aligned}
\end{equation}
For the first identity, expand
\(\binom nk^2=(\binom{n-1}k+\binom{n-1}{k-1})^2\).
The two squares give the coefficient of \(x^k\) in
\((1+x)M_{n-1}\), and the cross term gives the coefficient of
\(x^k\) in \(2(n-1)xN_{n-1}\).
For the second, the coefficient of \(x^k\) on the right is
\[
 (k+1)\binom{n-1}{k+1}^2+(n-1-k)\binom{n-1}k^2
 =n\binom{n-1}k\binom{n-1}{k+1}.
\]
Indeed, \((k+1)\binom{n-1}{k+1}=(n-1-k)\binom{n-1}{k}\).
Using this once in each summand gives the common product
\(\binom{n-1}{k}\binom{n-1}{k+1}\) with multipliers
\(n-1-k\) and \(k+1\), whose sum is \(n\).
The result is \([x^k]\,n(n-1)N_{n-1}\).

At each zero \(r\) of \(M_{n-1}\), the second identity gives
\(n(n-1)N_{n-1}(r)=(1-r)M_{n-1}'(r)\).
Since \(1-r>0\), \(N_{n-1}(r)\) has the same sign as
\(M_{n-1}'(r)\). These signs alternate, so
\(N_{n-1}\) has a zero in every gap between consecutive
zeros of \(M_{n-1}\). These \(n-2\) zeros exhaust the degree of
\(N_{n-1}\), and \(N_{n-1}\prec M_{n-1}\). For \(n=2\), this is the
convention \(N_1=1\prec1+x=M_1\).

At a zero \(r\) of \(N_{n-1}\), strict interlacing places one more
zero of \(M_{n-1}\) than of \(N_{n-1}\) to the right of \(r\).
Both polynomials are monic, so the product for \(M_{n-1}(r)\)
has one more negative factor than the product for \(N_{n-1}'(r)\).
Thus \(M_{n-1}(r)/N_{n-1}'(r)<0\). The first identity in
\eqref{eq:mixed-narayana-identities} gives $M_n(r)=(1+r)M_{n-1}(r)$,
so $M_n(r)/N_{n-1}'(r)>0$ when \(r<-1\).
\end{proof}

\begin{lemma}\label{lem:mixed-narayana-endpoint-values}
For \(n\ge2\) and real \(\beta\),
\[
 M_n(1)+\beta nN_{n-1}(1)=\binom{2n-2}{n-1}\left(\beta+4-\frac2n\right).
\]
For even \(n\),
\[
 M_n(-1)-\beta nN_{n-1}(-1)
 =\frac{2n-2+\beta n}{2(n-1)}M_n(-1),\qquad M_n(-1)\ne0.
\]
For odd \(n\), the factor \(1+x\) divides \(M_n(x)+\beta n\,xN_{n-1}(x)\)
and
\[
 \left.\frac{M_n(x)+\beta n\,xN_{n-1}(x)}{1+x}\right|_{x=-1}
 =n(1+\beta/2)M_{n-1}(-1),\qquad M_{n-1}(-1)\ne0.
\]
\end{lemma}
\begin{proof}
At \(x=1\), the summands in \(M_n(1)\) and \((n-1)N_{n-1}(1)\)
are \(\binom nk\binom n{n-k}\) and
\(\binom{n-1}k\binom{n-1}{n-2-k}\), respectively,
by~\eqref{eq:narayana-coefficients}. Vandermonde's identity gives
\[
 M_n(1)=\binom{2n}n,\qquad
 nN_{n-1}(1)=\frac{n}{n-1}\binom{2n-2}{n-2}
           =\binom{2n-2}{n-1}.
\]
The ratio $M_n(1)/(nN_{n-1}(1))$ is $(2n)(2n-1)/n^2=4-2/n$.
Factoring \(\binom{2n-2}{n-1}\) out of \(M_n(1)+\beta nN_{n-1}(1)\)
now gives the stated value at \(1\).

To evaluate at \(-1\), eliminate \(xN_{n-1}\) by the first
identity in~\eqref{eq:mixed-narayana-identities}. We obtain
\[
 2(n-1)\bigl(M_n+\beta n\,xN_{n-1}\bigr)
 =(2n-2+\beta n)M_n-\beta n(1+x)M_{n-1}.
\]
Suppose first that \(n\) is even. The simple negative zeros of
\(M_n\) other than \(-1\) occur in reciprocal pairs. A single
zero at \(-1\) would leave an odd number to pair. Hence \(M_n(-1)\ne0\).
At $x=-1$ the term containing $1+x$ vanishes, leaving
$(2n-2+\beta n)M_n(-1)$ on the right. Thus \(M_n+\beta n xN_{n-1}\)
vanishes at \(-1\) exactly when \(\beta=-2+2/n\).

When \(n\) is odd, symmetry in the odd degrees \(n\) and \(n-2\)
forces \(1+x\) to divide both \(M_n\) and \(N_{n-1}\).
To determine when the combination has a multiple zero at \(-1\),
we evaluate \((M_n+\beta n\,xN_{n-1})/(1+x)\) there.
This quotient value equals the derivative of
\(M_n+\beta n\,xN_{n-1}\) at \(-1\).
With \(F_n=M_n\), \(d=n\), and $x=-1$, identity~\eqref{eq:narayana-derivative}
gives $-2M_n'(-1)=n\bigl(M_n(-1)-2M_{n-1}(-1)\bigr)$.
Since $M_n(-1)=0$, we have $M_n'(-1)=nM_{n-1}(-1)$.
Differentiating the expression for \(2(n-1)(M_n+\beta n xN_{n-1})\)
and setting \(x=-1\), we obtain
\[
 \begin{aligned}
 \frac{(2n-2+\beta n)M_n'(-1)-\beta nM_{n-1}(-1)}{2(n-1)}
 &=\frac{n(2n-2+\beta n)-\beta n}{2(n-1)}M_{n-1}(-1)\\
 &=n(1+\beta/2)M_{n-1}(-1).
 \end{aligned}
\]
The even-degree case gives \(M_{n-1}(-1)\ne0\). Thus for odd \(n\),
the quotient \((M_n+\beta n xN_{n-1})/(1+x)\) vanishes at \(-1\)
exactly when \(\beta=-2\).
\end{proof}

\begin{proof}[Proof of Theorem~\ref{thm:narayana-branden-range}]
The polynomials \(M_n\) and \(N_{n-1}\) are monic of degrees \(n\)
and \(n-2\), respectively, and symmetric in those degrees. Since
\(xN_{n-1}\) has degree \(n-1\) and zero constant term, \(M_n+\beta n\,xN_{n-1}\) is monic of
degree \(n\) with constant term one. It is symmetric in degree \(n\) since
\(x^n(x^{-1}N_{n-1}(1/x))=xN_{n-1}(x)\).
Lemma~\ref{lem:mixed-narayana-residue-sign} gives
\(M_n(r)/N'_{n-1}(r)>0\) at every zero \(r<-1\) of \(N_{n-1}\).
Hence Lemma~\ref{lem:reciprocal-pair} applies to the quotient
\(Q(u)=M_n(x)/(xN_{n-1}(x))\), where \(u=x+x^{-1}\).
We cancel the common factor \(1+x\) when \(n\) is odd.
For \(xN_{n-1}(x)\ne0\), dividing \(M_n+\beta n xN_{n-1}\) by
$xN_{n-1}$ gives $Q(u)+\beta n$. Thus $x$ is a zero of the polynomial
exactly when \(Q(u)=-\beta n\).
At \(x=0\) the combination equals one. At any zero \(r\ne-1\) of
\(N_{n-1}\), the combination has value \(M_n(r)=(1+r)M_{n-1}(r)\ne0\)
by~\eqref{eq:mixed-narayana-identities} and \(N_{n-1}\prec M_{n-1}\).
For odd \(n\), the restored factor \(1+x\) contributes one copy of
the zero \(-1\), whose multiplicity is higher precisely when
\((M_n+\beta n xN_{n-1})/(1+x)\) also vanishes at \(-1\).

All poles of \(Q\) lie below \(-2\), so the rightmost interval
contains \([-2,2]\). On that interval \(Q\) increases from
\(-\infty\) to \(+\infty\). For \(n=2\) we have \(Q(u)=u+4\).
For \(n=3\), canceling \(1+x\) in
\((1+x)(x^2+8x+1)/(x(1+x))\) gives \(Q(u)=u+8\).
For \(n=2,3\), the rightmost interval is therefore the whole real line.
For \(n\ge2\), let \(v\) be the largest solution of \(Q(u)=-\beta n\).
This is the unique solution on the rightmost interval. As \(\beta\)
increases, \(-\beta n\) decreases.
Hence \(v\) decreases strictly, passing through \(2\) and then \(-2\).
The factorization in
Lemma~\ref{lem:reciprocal-pair} contains the factor \(x^2-vx+1\).
Every other solution \(u\) lies to the left of a pole and hence
is less than \(-2\). The quadratic \(x^2-ux+1\) then has positive
discriminant. Its two distinct zeros have product one and negative
sum \(u\), so both are negative.
Neither equals \(-1\), for this would give \(u=-2\). Two
different values of \(u\) cannot give the same zero \(x\), since
\(x\ne0\) determines \(u=x+x^{-1}\).
Thus the quotient of \(M_n+\beta n\,xN_{n-1}\) by \(x^2-vx+1\)
has \(n-2\) simple negative zeros, including \(-1\) when \(n\) is odd.
The discriminant of \(x^2-vx+1\) is \(v^2-4\).
When \(-2<v<2\), its zeros are nonreal conjugates. Their product
is one, so both lie on the unit circle.

All the other factors have only negative zeros. Hence \(v=2\)
if and only if \(M_n(1)+\beta nN_{n-1}(1)=0\).
Lemma~\ref{lem:mixed-narayana-endpoint-values} gives \(\beta=-4+2/n\).
At \(\beta=-4+2/n\) the factor \(x^2-vx+1\) is \((x-1)^2\), so \(1\)
is a double zero.

At \(-1\), the quadratic \(x^2-vx+1\) has value \(2+v\).
All other quadratic factors are nonzero there, since their parameters
are less than \(-2\). Hence \(v=-2\) exactly when
\(M_n+\beta n\,xN_{n-1}\) vanishes at \(-1\) for even \(n\), or when
its quotient by \(1+x\) vanishes there for odd \(n\).
By Lemma~\ref{lem:mixed-narayana-endpoint-values}, these conditions give
\(\beta=-2+2/n\) for even \(n\) and \(\beta=-2\) for odd \(n\).
In both cases \(\beta=-(n-1)/\lfloor n/2\rfloor\) and
\(x^2-vx+1=(x+1)^2\). Thus \(-1\) is a double zero in even degree.
The fixed factor \(1+x\) makes it a triple zero in odd degree.

The quadratic \(x^2-vx+1\) has real zeros precisely when
\(v\ge2\) or \(v\le-2\). Since \(v\) decreases strictly with \(\beta\),
the first inequality holds exactly for \(\beta\le-4+2/n\) and the
second exactly for \(\beta\ge-(n-1)/\lfloor n/2\rfloor\).
If \(\beta<-4+2/n\), then \(v>2\), so \(x^2-vx+1\) has two
distinct positive zeros. If \(\beta>-(n-1)/\lfloor n/2\rfloor\), then
\(v<-2\), so this factor also has two distinct negative zeros.
\end{proof}

\begin{remark*}
If \(\beta<-4+2/n\), then \(M_n(x)+\beta n\,xN_{n-1}(x)\) has two
positive zeros and \(n-2\) negative zeros, all simple.
If \(-4+2/n<\beta<-(n-1)/\lfloor n/2\rfloor\), then
\(M_n(x)+\beta n\,xN_{n-1}(x)\) has \(n-2\) simple negative zeros
and one nonreal conjugate pair on the unit circle. At \(\beta=-4+2/n\), the
zero \(1\) is double. At \(\beta=-(n-1)/\lfloor n/2\rfloor\), the
zero \(-1\) is double for even \(n\) and triple for odd \(n\).
All other zeros at either endpoint are simple and negative.
If \(\beta>-(n-1)/\lfloor n/2\rfloor\), then all \(n\) zeros are
simple and negative.
\end{remark*}

\begin{example}
For \(n=4\), the choice \(\beta=-1\) gives the \(h\)-polynomial
of the type~\(D_4\) cluster complex,
\(M_4(x)-4xN_3(x)=1+12x+24x^2+12x^3+x^4\).
Its four zeros are simple and negative since \(-1>-3/2\).
The two endpoints of the parameter ranges give
\[
 \begin{aligned}
 M_4-14xN_3&=(x-1)^2(x^2+4x+1) &&(\beta=-7/2),\\
 M_4-6xN_3&=(x+1)^2(x^2+8x+1) &&(\beta=-3/2).
 \end{aligned}
\]
Chen, Yang, and Zhang's parameter \(\beta=-19/10\) lies strictly
between these two endpoints. The two nonreal zeros of
\(M_4-\tfrac{38}{5}xN_3\) therefore lie on the unit circle.
To compute the quadratic factors, put \(u=x+x^{-1}\).
Then
\[
 x^{-2}\bigl(M_4(x)-\tfrac{38}{5}xN_3(x)\bigr)
 =u^2+\tfrac{42}{5}u+\tfrac{56}{5}.
\]
The two roots in \(u\) are $u_\pm=(-21\pm\sqrt{161})/5$.
Since \(11<\sqrt{161}<13\), we have \(u_-<-2\) and
\(-2<u_+<0\). The corresponding factors of
\(M_4(x)-\tfrac{38}{5}xN_3(x)\) are
\[
 \begin{aligned}
 x^2-u_-x+1&=x^2+\frac{21+\sqrt{161}}5x+1,\\
 x^2-u_+x+1&=x^2+\frac{21-\sqrt{161}}5x+1.
 \end{aligned}
\]
The zeros of \(x^2-u_-x+1\) are real and distinct because \(u_-<-2\).
Their sum is negative and their product is one, so both are negative.
Since \(-2<u_+<0\), the zeros of \(x^2-u_+x+1\) are nonreal
conjugates. Their product is one, so both lie on the unit circle.
\end{example}

For \(n\ge4\), the inequality \(\lfloor n/2\rfloor<n-1\) gives
\(-1>-(n-1)/\lfloor n/2\rfloor\). Thus \(\beta=-1\) lies in the
range where all \(n\) zeros are simple and negative. This recovers
Br\"and\'en's conclusion for the type~\(D_n\) cluster \(h\)-polynomial
\(M_n(x)-nxN_{n-1}(x)\)~\cite[Corollary~7.2]{Branden2006}.

\section{\texorpdfstring{\(q\)}{q}-Eulerian polynomials}
\label{cycle-weighted-Eulerian}\label{sec:cycle-interlacing}

Recall that \(A_n(x;q)\) counts permutations of \([n]\) by excedances,
with weight \(q\) per cycle, including fixed points. For a positive
integer \(q\), this gives each cycle a choice of \(q\) colors.
For \(p,q>0\), we first compare \(A_n(x;p)\) with
\(A_{n+1}(x;q)\) and with \(A_n(x;q)\). We then compare
\(A_n(x;q)\) with the type~\(B\) Eulerian polynomial
\(B_n(x;p)\), where \(p\) weights each negative entry of a signed
permutation. Finally, we study three-term combinations and the
Eulerian transformation.
Brenti's recurrence~\cite[Proposition~7.2]{Brenti2000} is
\begin{equation}\label{eq:cycle-insertion}
 A_{n+1}(x;q)=(q+nx)A_n(x;q)+x(1-x)A_n'(x;q)
 \qquad(n\ge0).
\end{equation}
To see this, let \(\pi\in\mathfrak S_n\) have \(k\) excedances.
Making \(n+1\) a singleton cycle multiplies the weight by \(q\).
Inserting \(n+1\) after \(i\) in a cycle replaces
\(i\mapsto\pi(i)\) by \(i\mapsto n+1\mapsto\pi(i)\).
After insertion, the position \(i\) is an excedance, the new position
\(n+1\) is not, and all other positions keep their excedance status.
Thus insertion after any of the \(k\) excedances preserves the number
of excedances. Insertion after any of the other \(n-k\) positions adds
one excedance.
The total contribution from \(\pi\) is
$q^{\operatorname{cyc}(\pi)}\bigl((q+k)x^k+(n-k)x^{k+1}\bigr)$.
Summing over \(\pi\) gives~\eqref{eq:cycle-insertion}.

Only the identity permutation has no excedances, so \(A_n(0;q)=q^n\).
For \(q>0\) and \(n\ge1\), the leading term of \(A_n(x;q)\) is
\(qx^{n-1}\). Indeed, \(A_1(x;q)=q\). If \(A_n(x;q)\) has
leading term \(qx^{n-1}\), then~\eqref{eq:cycle-insertion} gives
the coefficient of \(x^n\) in \(A_{n+1}(x;q)\) as
\(nq-(n-1)q=q\).

For an exponential generating function of \(A_n(x;q)\), see
\cite[(10)]{MaMaYehYeh2024} with \(p=1\).
Brenti~\cite[Theorem~7.5]{Brenti2000} proved that \(A_n(x;q)\) has
only simple negative zeros for positive rational \(q\).
For real \(q>0\), see Br\"and\'en~\cite[Theorem~6.1]{Branden2006}
and Ma and Wang~\cite[Theorem~2 and Proposition~4(a)]{MaWang2008}.

To prove that \(A_n(x;q)\) strictly interlaces \(A_{n+1}(x;q)\), we
use~\eqref{eq:cycle-insertion}. Start with \(A_1(x;q)=q\) and
\(A_2(x;q)=q(x+q)\). If \(A_n(x;q)\) has simple negative zeros
\(u_1<\cdots<u_{n-1}\), then $A_{n+1}(u_i;q)=u_i(1-u_i)A_n'(u_i;q)$.
At \(u_i\), the derivative \(A_n'(u_i;q)\) has sign
\((-1)^{n-1-i}\), and \(u_i(1-u_i)<0\). Thus
\(A_{n+1}(u_i;q)\) has sign \((-1)^{n-i}\).
The signs alternate from one \(u_i\) to the next. The sign of
\(A_{n+1}(x;q)\) is \((-1)^n\) near \(-\infty\) and positive at \(0\),
where \(A_{n+1}(0;q)=q^{n+1}\). Thus each of the \(n\) intervals
cut out of \((-\infty,0)\) by the \(u_i\) contains a zero of
\(A_{n+1}(x;q)\). Since \(A_{n+1}(x;q)\) has degree \(n\), each
interval contains exactly one simple zero.

\subsection{Interlacing for two cycle weights}

For fixed \(q>0\), which weights \(p>0\) make \(A_n(x;p)\)
strictly interlace \(A_{n+1}(x;q)\) for every \(n\ge2\)?

\begin{theorem}\label{thm:cycle-two-weights}
Let \(p,q>0\). Then
\[
 A_n(x;p)\prec A_{n+1}(x;q)\quad\text{for every }n\ge2
 \quad\Longleftrightarrow\quad q\le p\le q+1.
\]
\end{theorem}

For sufficiency, we use the recurrence~\eqref{eq:cycle-insertion} to
apply the two-parameter interlacing theorem at each step of the induction.
For necessity, we compare the
constant, linear, and quadratic coefficients of two strictly interlacing
polynomials. The next lemma gives the inequalities we need here and
for type~\(B\) in Section~\ref{sec:independent-AB-weights}.

\begin{firstcoefflemma}
Suppose \(f\prec g\), where \(f\) and \(g\) have simple negative zeros,
positive leading coefficients and \(\deg g=\deg f+1\ge2\). Then
\begin{equation}\label{eq:interlacing-first-coefficients}
 \frac{[x^2]g}{[x]g}<\frac{[x]f}{f(0)}
 <\frac{[x]g}{g(0)}.
\end{equation}
\end{firstcoefflemma}
\begin{proof}
Let \(n=\deg g\) and write $g(x)/g(0)=\prod_{i=1}^n(1+v_i x)$,
where $v_1>\cdots>v_n>0$.
The zeros of \(g\), from right to left, are
\(-1/v_1,\ldots,-1/v_n\). Since \(f\prec g\), each interval
\((-1/v_{i+1},-1/v_i)\) contains one zero \(u\) of \(f\).
The function \(x\mapsto-1/x\) is strictly increasing on
\((-\infty,0)\), so \(v_{i+1}<-1/u<v_i\).
The factorization \(f(x)/f(0)=\prod_{f(u)=0}(1-x/u)\) gives
\([x]f/f(0)=\sum_{f(u)=0}(-1/u)\). Summing the interval bounds
and then using the linear coefficient of \(g/g(0)\) gives
\[
 \sum_{i=2}^n v_i<\frac{[x]f}{f(0)}
 <\sum_{i=1}^{n-1}v_i<\sum_{i=1}^n v_i=\frac{[x]g}{g(0)}.
\]
To bound \(\sum_{k=2}^n v_k\) from below, group the terms
\(v_i v_k\), \(i<k\), in the quadratic coefficient of \(g/g(0)\)
by their index \(k\). Then
\[
 \frac{[x^2]g}{[x]g}
  =\sum_{k=2}^n v_k\,
        \frac{\sum_{i<k}v_i}{\sum_{i=1}^n v_i}
  <\sum_{k=2}^n v_k.
\]
Combining the two inequalities gives~\eqref{eq:interlacing-first-coefficients}.
\end{proof}

\begin{proof}[Proof of Theorem~\ref{thm:cycle-two-weights}]

Suppose that \(q\le p\le q+1\). By~\eqref{eq:cycle-insertion},
\(A_{n+1}(x;p)=T_{n-1}^{p,1}A_n(x;p)\) and \(A_{n+2}(x;q)=T_n^{q,1}A_{n+1}(x;q)\).
The initial polynomials \(A_1(x;p)=p\) and \(A_2(x;q)=q(x+q)\)
satisfy \(A_1(x;p)\prec A_2(x;q)\). If
\(A_n(x;p)\prec A_{n+1}(x;q)\), Theorem~\ref{thm:two-endpoint-interlacing}
in the variable \(x\), as given after~\eqref{eq:cayley-two-endpoint},
applies since \(0\le p-q\le1\). Thus
\(A_{n+1}(x;p)\prec A_{n+2}(x;q)\). Induction proves sufficiency.

To derive the necessary bounds on \(p\) and \(q\), we compute
the coefficients in~\eqref{eq:interlacing-first-coefficients}.
For any fixed \(s>0\), recurrence~\eqref{eq:cycle-insertion} gives
\[
 \frac{A_{n+1}(x;s)}{(1-x)^{n+s+1}}
 =\left(s+x\frac{d}{dx}\right)
   \frac{A_n(x;s)}{(1-x)^{n+s}}.
\]
The operator \(s+x\,d/dx\) multiplies the coefficient of
\(x^j\) by \(s+j\). Starting with the binomial series for
\((1-x)^{-s}\) at \(n=0\) therefore gives
\begin{equation}\label{eq:cycle-weight-series}
 \frac{A_n(x;s)}{(1-x)^{n+s}}
 =\sum_{j\ge0}\binom{s+j-1}{j}(s+j)^n x^j.
\end{equation}
Multiplying by \((1-x)^{n+s}\) and comparing coefficients,
we obtain
\[
 \begin{aligned}
 A_n(0;s)&=s^n,\\
 [x]A_n(x;s)&=s(s+1)^n-(n+s)s^n,\\
 [x^2]A_n(x;s)&=\frac{s(s+1)}2(s+2)^n
   -s(n+s)(s+1)^n+\binom{n+s}{2}s^n.
 \end{aligned}
\]
For fixed \(s>0\), divide the linear and quadratic coefficients by
\((s+1)^n\) and \((s+2)^n\), respectively. In each case every term
except the first tends to zero, so
\[
 [x]A_n(x;s)\sim s(s+1)^n,\qquad
 [x^2]A_n(x;s)\sim\frac{s(s+1)}2(s+2)^n.
\]
Thus the \(n\)th roots of the constant, linear and quadratic
coefficients tend to \(s,s+1,s+2\), respectively.

Now suppose that \(A_n(x;p)\prec A_{n+1}(x;q)\) for every \(n\ge2\).
By~\eqref{eq:interlacing-first-coefficients},
\[
 \frac{[x^2]A_{n+1}(x;q)}{[x]A_{n+1}(x;q)}
 <\frac{[x]A_n(x;p)}{A_n(0;p)}
 <\frac{[x]A_{n+1}(x;q)}{A_{n+1}(0;q)}.
\]
Taking \(n\)th roots and using the coefficient limits with
\((n+1)/n\to1\) gives
\[
 \frac{q+2}{q+1}\le\frac{p+1}{p}\le\frac{q+1}{q}.
\]
Subtracting one from each term gives
\(1/(q+1)\le1/p\le1/q\). Since \(p,q>0\), this is precisely
\(q\le p\le q+1\).
\end{proof}

We now compare the cycle weights \(p\) and \(q\) at the same size \(n\).

\begin{corollary*}
Let \(p,q>0\). The polynomial \(A_n(x;p)\) strictly alternates left
of \(A_n(x;q)\) for every \(n\ge2\) if and only if \(q<p\le q+1\).
\end{corollary*}

\begin{proof}
First compare the weights \(q\) and \(q+1\). Differentiate
\eqref{eq:cycle-weight-series} with respect to \(x\), taking \(s=q\).
The coefficient of \(x^j\)
in the derivative is
\((j+1)\binom{q+j}{j+1}(q+j+1)^n=q\binom{q+j}{j}(q+j+1)^n\).
By~\eqref{eq:cycle-weight-series} with \(s=q+1\), the derivative equals
\(qA_n(x;q+1)/(1-x)^{n+q+1}\). Multiplying by \((1-x)^{n+q+1}\) gives
\[
 qA_n(x;q+1)=(1-x)A_n'(x;q)+(n+q)A_n(x;q).
\]
At each zero \(u\) of \(A_n(x;q)\), the sign of \(A_n(u;q+1)\)
is therefore the sign of \(A_n'(u;q)\). The values of \(A_n(x;q+1)\)
alternate in sign at the zeros of \(A_n(x;q)\). At the smallest zero of
\(A_n(x;q)\), the sign is opposite to the sign of \(A_n(x;q+1)\)
near \(-\infty\).
Thus \(A_n(x;q+1)\) has one zero to the left of the
smallest zero of \(A_n(x;q)\) and one in each successive gap.

For \(s>0\), write \(u_1(s)<\cdots<u_{n-1}(s)\) for the zeros of
\(A_n(x;s)\). The coefficients of \(A_n(x;s)\) are polynomials in \(s\).
Since the zeros are simple, the implicit function theorem shows that
each \(u_j(s)\) is differentiable.
We show by induction on \(n\) that \(u_j'(s)<0\) for every \(j\).
At \(n=2\) the only zero is \(-s\). Suppose the assertion holds at \(n\).
Away from the zeros of \(A_n(x;s)\), put
\(R(x;s)=A_{n+1}(x;s)/A_n(x;s)\). By~\eqref{eq:cycle-insertion},
\[
 R(x;s)=s+nx+x(1-x)\sum_{j=1}^{n-1}\frac{1}{x-u_j(s)}.
\]
Differentiating with respect to \(s\) at fixed \(x<0\) gives
\[
 \partial_s R(x;s)
 =1+x(1-x)\sum_{j=1}^{n-1}\frac{u_j'(s)}{(x-u_j(s))^2}>0.
\]
Indeed, \(x(1-x)<0\) and \(u_j'(s)<0\), so each term after \(1\)
is positive. Let \(v=v(s)\) be the \(i\)th zero of \(A_{n+1}(x;s)\)
from the left. There are \(n-i\) zeros of \(A_{n+1}\) to its right
and, by strict interlacing, \(n-i\) zeros of \(A_n\) to its right.
By the positive leading coefficients, \(A_{n+1}'(v;s)\) and \(A_n(v;s)\)
both have sign \((-1)^{n-i}\), so
\((\partial_x R)(v;s)=A_{n+1}'(v;s)/A_n(v;s)>0\).
Differentiating \(R(v(s);s)=0\) gives
\[
 v'(s)=-\frac{(\partial_s R)(v;s)}{(\partial_x R)(v;s)}<0.
\]
This completes the induction. For \(q<p<q+1\), we therefore have
\(u_j(q+1)<u_j(p)<u_j(q)\) for \(1\le j\le n-1\).
For \(1\le j<n-1\), the comparison at \(q+1\) now gives
\(u_j(p)<u_j(q)<u_{j+1}(q+1)<u_{j+1}(p)\).
Together with \(u_{n-1}(p)<u_{n-1}(q)\), this proves the stated order.

Conversely, suppose \(A_n(x;p)\) strictly alternates left of \(A_n(x;q)\)
for every \(n\ge2\). At \(n=2\), the two zeros satisfy \(-p<-q\),
so \(p>q\).
For each \(n\), choose \(M>0\) so that \(-M\) lies to the left of all
zeros of both polynomials. The factor \(1+x/M\) adds a leftmost
zero to \(A_n(x;q)\), so the assumed order now gives
\(A_n(x;p)\prec(1+x/M)A_n(x;q)\).
Apply the left inequality in~\eqref{eq:interlacing-first-coefficients}.
As \(M\to\infty\), each coefficient of \((1+x/M)A_n(x;q)\) tends to the
corresponding coefficient of \(A_n(x;q)\), giving
\[
 \frac{[x^2]A_n(x;q)}{[x]A_n(x;q)}
 \le \frac{[x]A_n(x;p)}{A_n(0;p)}.
\]
Taking \(n\)th roots and using the coefficient limits in the proof of
Theorem~\ref{thm:cycle-two-weights} yields
\((q+2)/(q+1)\le(p+1)/p\), or \(p\le q+1\).
\end{proof}

Set one cycle weight equal to \(1\) and denote the other by \(q\).
Since \(A_n(x;1)=A_n(x)\), the corollary gives, for every \(n\ge2\),
\[
\begin{array}{ll}
 A_n(x)\prec A_n(x;q), & 0<q<1,\\[2pt]
 A_n(x;q)\prec A_n(x), & 1<q\le2.
\end{array}
\]
At \(q=1\) the two polynomials coincide.

Let \(p\) and \(q\) be positive integers. Give each cycle of a permutation
of \([n]\) one of \(p\) colors, and each cycle of a permutation of
\([n+1]\) one of \(q\) colors. Counting by excedances gives \(A_n(x;p)\)
and \(A_{n+1}(x;q)\), respectively. By Theorem~\ref{thm:cycle-two-weights}, \(A_n(x;p)\) strictly interlaces
\(A_{n+1}(x;q)\) for every \(n\ge2\) exactly when \(p=q\) or \(p=q+1\).
For \(p=2\) and \(q=1\), each cycle on \([n]\) is red or blue.
The resulting polynomial \(A_n(x;2)\) strictly alternates left of \(A_n(x)\) and strictly
interlaces \(A_{n+1}(x)\) for every \(n\ge2\).

\subsection{Interlacing of Eulerian polynomials of types \texorpdfstring{\(A\)}{A} and \texorpdfstring{\(B\)}{B}}
\label{sec:independent-AB-weights}

Let \(N(w)\) denote the number of negative entries in a signed
permutation \(w=w_1\cdots w_n\). Set \(w_0=0\), and let \(d_B(w)\)
be the number of indices \(0\le i<n\) for which \(w_i>w_{i+1}\).
Following Brenti~\cite{Brenti1994} and Br\"and\'en~\cite[(6.2)]{Branden2006},
put
\[
 B_n(x;p)=\sum_w p^{N(w)}x^{d_B(w)},\qquad B_0(x;p)=1,
\]
where \(w\) ranges over all words \(w_1\cdots w_n\) such that
\((|w_1|,\ldots,|w_n|)\) is a permutation of \([n]\).
At \(p=1\) this gives \(B_n(x)\). At \(p=0\) only ordinary
permutations contribute, so \(B_n(x;0)=A_n(x)\).
Br\"and\'en~\cite[Theorem~6.7]{Branden2006} proved that
\(A_n(x)\prec B_n(x;p)\) for every \(n\ge1\) and \(p>0\).
Since \(A_n(x;1)=A_n(x)\), Br\"and\'en's result gives the case \(q=1\)
of Theorem~\ref{thm:independent-AB-weights}.
We give an operator proof of the case \(p=q=1\) at the end of
Section~\ref{sec:two-parameter}.
For fixed \(p>0\), which cycle weights \(q\) make \(A_n(x;q)\)
strictly interlace \(B_n(x;p)\) for every \(n\ge2\)?

\begin{theorem}\label{thm:independent-AB-weights}
Let \(p,q>0\). Then
\[
 A_n(x;q)\prec B_n(x;p)\quad\text{for every }n\ge2
 \quad\Longleftrightarrow\quad
 \frac{1}{1+p}\le q\le1+\frac{1}{1+p}.
\]
\end{theorem}

\begin{proof}
We use the two-parameter interlacing theorem for sufficiency and
compare the first three coefficients for necessity.

Brenti's identity
\cite{Brenti1994}, also given in~\cite[(6.3)]{Branden2006}, reads
\begin{equation}\label{eq:negative-entry-series}
 \frac{B_n(x;p)}{(1-x)^{n+1}}
 =\sum_{j\ge0}\bigl(1+(1+p)j\bigr)^n x^j .
\end{equation}
The operator \(1+(1+p)x\,d/dx\) multiplies the coefficient of
\(x^j\) by \(1+(1+p)j\). Applying the operator to~\eqref{eq:negative-entry-series}
therefore gives the series for \(B_{n+1}(x;p)/(1-x)^{n+2}\). Differentiating
the quotient on the left and multiplying by \((1-x)^{n+2}\) gives
\[
 (1-x)B_n(x;p)+(1+p)x\bigl((n+1)B_n(x;p)+(1-x)B_n'(x;p)\bigr)
 = B_{n+1}(x;p).
\]
Collecting terms yields
\begin{equation}\label{eq:negative-entry-recurrence}
 B_{n+1}(x;p)=\bigl(1+(n+(n+1)p)x\bigr)B_n(x;p)
               +(1+p)x(1-x)B_n'(x;p).
\end{equation}
Only \(1,2,\ldots,n\) has no type~\(B\) descent, and only
\(-1,-2,\ldots,-n\) has \(n\) such descents. Thus \(B_n(x;p)\)
has constant term one and leading term \(p^n x^n\).

Suppose that
\(\frac{1}{1+p}\le q\le1+\frac{1}{1+p}\).
The recurrences~\eqref{eq:cycle-insertion} and
\eqref{eq:negative-entry-recurrence} give
\[
 \begin{aligned}
 A_{n+1}(x;q)&=T_{n-1}^{q,1}A_n(x;q),\\
 \frac{B_{n+1}(x;p)}{1+p}
 &=T_n^{1/(1+p),\,p/(1+p)}B_n(x;p).
 \end{aligned}
\]
We have \(A_1(x;q)=q\prec 1+px=B_1(x;p)\). Suppose that
\(A_n(x;q)\prec B_n(x;p)\). The differences between the two pairs
of parameters are
\(q-1/(1+p)\in[0,1]\) and \(1-p/(1+p)=1/(1+p)\in(0,1)\).
Theorem~\ref{thm:two-endpoint-interlacing}, in the form given after
\eqref{eq:cayley-two-endpoint}, therefore gives
\(A_{n+1}(x;q)\prec B_{n+1}(x;p)\). Induction proves sufficiency.

Conversely, suppose \(A_n(x;q)\prec B_n(x;p)\) for every \(n\ge2\).
The strict interlacing gives simple real zeros. Both polynomials have
nonnegative coefficients and positive constant terms, so none of
these zeros lies in \([0,\infty)\). Their degrees differ by one.
The coefficient comparison~\eqref{eq:interlacing-first-coefficients}
therefore gives
\[
 \frac{[x^2]B_n(x;p)}{[x]B_n(x;p)}
 <\frac{[x]A_n(x;q)}{A_n(0;q)}
 <\frac{[x]B_n(x;p)}{B_n(0;p)}.
\]
To evaluate the two outer ratios, multiply~\eqref{eq:negative-entry-series}
by \((1-x)^{n+1}\). The first three coefficients are
\[
 \begin{aligned}
 B_n(0;p)&=1,\\
 [x]B_n(x;p)&=(2+p)^n-(n+1),\\
 [x^2]B_n(x;p)&=(3+2p)^n-(n+1)(2+p)^n+\binom{n+1}{2}.
 \end{aligned}
\]
Thus the \(n\)th roots of the linear and quadratic coefficients of
\(B_n(x;p)\) tend to \(2+p\) and \(3+2p\), respectively. For the
middle ratio, the proof of Theorem~\ref{thm:cycle-two-weights} gives
\(A_n(0;q)=q^n\) and \(\bigl([x]A_n(x;q)\bigr)^{1/n}\to q+1\).
Taking \(n\)th roots in the coefficient inequalities and letting
\(n\to\infty\) gives
\[
 \frac{3+2p}{2+p}\le\frac{q+1}{q}\le 2+p.
\]
Subtracting one gives \((1+p)/(2+p)\le1/q\le1+p\).
Since all quantities are positive, this is precisely
\(1/(1+p)\le q\le(2+p)/(1+p)\), as claimed.
\end{proof}

The choice \(q=(1+p)/2\) has a direct interpretation in terms of
signed permutations.
A cycle of the ordinary permutation underlying a signed permutation
is \emph{negative} if an odd number of its letters carry a minus sign.
In \(B_n(x;p)\),
each negative entry has weight \(p\), and \(x\) records type~\(B\) descents.
On the same set of signed permutations, give each negative cycle weight
\(p\) and let \(x\) record the excedances of the underlying permutation. For example,
\((-2,-1,-3)\) has three negative entries and two type~\(B\) descents.
Its underlying permutation has one excedance and two cycles,
\((1\,2)\) and \((3)\), of which only \((3)\) is negative. The
signed permutation contributes \(p^3x^2\) and \(px\) to the two polynomials,
respectively.

For a cycle of length \(j\), flip the sign of its smallest letter.
This pairs each assignment with an even number of minus signs with
one having an odd number. The two assignments have cycle weights \(1\)
and \(p\), so all assignments to this cycle have total weight \(2^{j-1}(1+p)\).
Multiplying over all cycles gives a factor \(2^n\) and a weight
\((1+p)/2\) for each cycle. Thus the excedance polynomial is
\(2^n A_n(x;(1+p)/2)\). For \(p>0\),
Theorem~\ref{thm:independent-AB-weights} now gives
\[
 A_n\!\left(x;\frac{1+p}{2}\right)\prec B_n(x;p)
 \quad\text{for every }n\ge2
 \quad\Longleftrightarrow\quad \sqrt2-1\le p\le\sqrt3.
\]
Indeed, the two bounds on \(q=(1+p)/2\) reduce to
\((1+p)^2\ge2\) and \(p^2\le3\).

For \(p=1\), the two endpoints in Theorem~\ref{thm:independent-AB-weights}
are \(q=1/2\) and \(q=3/2\). Both have direct counting interpretations.
If every cycle must have an even number of negative entries, a cycle of length
\(j\) admits \(2^{j-1}\) sign assignments. If all sign assignments are
allowed and each negative cycle may be marked or left unmarked, each
even assignment gives one choice and each odd assignment gives two.
There are \(2^{j-1}\) assignments of each parity, giving
\(3\cdot2^{j-1}\) choices.
The corresponding excedance polynomials are \(2^nA_n(x;1/2)\) and
\(2^nA_n(x;3/2)\). Keep these two polynomials fixed and vary \(p>0\)
in \(B_n(x;p)\). By Theorem~\ref{thm:independent-AB-weights},
\(2^nA_n(x;1/2)\) strictly interlaces \(B_n(x;p)\) for every \(n\ge2\)
exactly when \(p\ge1\), while \(2^nA_n(x;3/2)\) does so exactly when
\(0<p\le1\).
Thus \(p=1\) is the unique weight for which both comparisons hold
for every \(n\ge2\).

\subsection{Three-term combinations}

At \(q=1\), \(A_n(x;q)-cxA_{n-2}(x;q)\) is the polynomial in
Zhang's conjecture. We now allow any \(q>0\) and add the term
\(bA_{n-1}(x;q)\). We seek the bound on \(c\) that guarantees
real-rootedness uniformly over \(n\ge3\) and \(b\in\mathbb R\).

\begin{theorem}\label{thm:cycle-free-middle}
Fix \(q>0\) and \(c\in\mathbb R\). The polynomial
\[
 A_n(x;q)+bA_{n-1}(x;q)-cxA_{n-2}(x;q)
\]
has only real zeros for all \(n\ge3\) and \(b\in\mathbb R\)
if and only if \(c\le q+1\).
\end{theorem}

\begin{proof}
Let \(n\ge4\) and let \(r_1<\cdots<r_{n-2}<0\) be the zeros of
\(A_{n-1}(x;q)\). Since \(A_{n-2}(x;q)\prec A_{n-1}(x;q)\), apply
Lemma~\ref{lem:free-middle-general} with \(F_j=A_j(x;q)\).
In the lemma's coefficient \(c-kx\), set the constant term to zero
and take \(k\) equal to the parameter \(c\) of this theorem.
By the lemma, the polynomial has only real zeros for all
\(b\in\mathbb R\) exactly when
\(0\le -A_n(r_i;q)/A_{n-2}(r_i;q)+cr_i\) at each \(r_i\).
Substituting \(A_n(r_i;q)=r_i(1-r_i)A_{n-1}'(r_i;q)\)
from~\eqref{eq:cycle-insertion} and dividing by \(r_i<0\) gives
\[
 c\le
 (1-r_i)\frac{A_{n-1}'(r_i;q)}{A_{n-2}(r_i;q)}
 \qquad(1\le i\le n-2).
\]
We use the partial-fraction expansion
\[
 \frac{A_{n-2}(x;q)}{A_{n-1}(x;q)}
 =\sum_{i=1}^{n-2}
   \frac{A_{n-2}(r_i;q)}{A_{n-1}'(r_i;q)}
   \frac{1}{x-r_i}.
\]
Strict interlacing shows that every residue is positive.
At \(x=1\), the \(i\)th summand is the reciprocal of the upper bound
on \(c\) at \(r_i\). By
\eqref{eq:cycle-insertion},
\(A_{n-1}(1;q)=(q+n-2)A_{n-2}(1;q)\), and hence
\[
 \sum_{i=1}^{n-2}
 \frac{A_{n-2}(r_i;q)}
      {(1-r_i)A_{n-1}'(r_i;q)}
 =\frac{1}{q+n-2}.
\]
Since \(n\ge4\), this sum has at least two positive terms. Each term
is therefore strictly smaller than the sum, and taking reciprocals
gives
\[
 (1-r_i)\frac{A_{n-1}'(r_i;q)}{A_{n-2}(r_i;q)}>q+n-2>q+1.
\]
Thus the residue inequalities are strict for \(c\le q+n-2\).
This includes \(c\le q+1\), the range independent of \(n\).
By Lemma~\ref{lem:free-middle-general},
\(A_{n-1}(x;q)\prec A_n(x;q)+bA_{n-1}(x;q)-cxA_{n-2}(x;q)\)
for every real \(b\).

It remains to consider \(n=3\). The recurrence gives
\(A_1(x;q)=q\), \(A_2(x;q)=q(x+q)\), and
\(A_3(x;q)=q\bigl(x^2+(3q+1)x+q^2\bigr)\). Substitution gives
\[
 \begin{aligned}
 &\frac{A_3(x;q)+bA_2(x;q)-cxA_1(x;q)}{q}\\
 &\qquad=(x+q)^2+(b+q+1-c)(x+q)+q(c-q-1).
 \end{aligned}
\]
At \(x=-q\), the quadratic has value \(q(c-q-1)\). For \(c<q+1\),
this value is negative, so its two zeros lie on opposite sides of \(-q\).
For \(c=q+1\), the quadratic is \((x+q)(x+q+b)\).
Thus \(A_2(x;q)\) interlaces the original polynomial for every real \(b\)
when \(c\le q+1\). If \(c>q+1\), taking \(b=c-q-1\) gives
\((x+q)^2+q(c-q-1)>0\) for real \(x\). This proves necessity
and completes the proof.
\end{proof}

For fixed \(n\), the bound on \(c\) can be larger. At \(n=4\) and \(q=1\),
the identity \(A_4=(1+x)A_3+6xA_2\) gives \(A_4(r)=6rA_2(r)\) at
both zeros \(r\) of \(A_3\). The criterion in the proof therefore gives
\(c\le6\) at this index. At \(c=6\), the three-term sum is \((1+x+b)A_3(x)\).

\begin{corollary*}
Let \(q>0\), \(n\ge3\), \(b\in\mathbb R\), and \(c\le q+1\).
The polynomial
\(A_n(x;q)+bA_{n-1}(x;q)-cxA_{n-2}(x;q)\)
has a multiple zero precisely in the case \(n=3\), \(c=q+1\), and \(b=0\).
In this case \(-q\) is a double zero.
\end{corollary*}

\begin{proof}
For \(n\ge4\), the proof of Theorem~\ref{thm:cycle-free-middle}
shows that \(A_{n-1}(x;q)\) strictly interlaces the three-term combination.
For \(n=3\) and \(c<q+1\), the two zeros of the combination lie on
opposite sides of \(-q\).
For \(n=3\) and \(c=q+1\), the polynomial is \(q(x+q)(x+q+b)\).
Its two zeros coincide exactly when \(b=0\).
\end{proof}

\subsection{Full interlacing and a symmetric matrix representation}

For \(n\ge3\) and \(c\le q+1\), the proof of Theorem~\ref{thm:cycle-free-middle}
gives
\[
 A_{n-1}(x;q)\preceq A_n(x;q)+bA_{n-1}(x;q)-cxA_{n-2}(x;q)
 \qquad(b\in\mathbb R).
\]
For \(-q\le b_1\le\cdots\le b_s\), substitute \(b=b_i\) in the
polynomial on the right. Theorem~\ref{thm:cycle-full-interlacing} shows
that the resulting row in the order \(i=1,\ldots,s\) is fully interlacing.

\begin{proof}[Proof of Theorem~\ref{thm:cycle-full-interlacing}]
Assume \(c\le q+1\), and put
\[
 H=A_n(x;q)-qA_{n-1}(x;q)-cxA_{n-2}(x;q).
\]
For \(n\ge4\), the proof of Theorem~\ref{thm:cycle-free-middle}
gives \(A_{n-1}(x;q)\prec H\). For \(n=3\),
\[
 A_2(x;q)=q(x+q),\qquad H=qx(x+2q+1-c).
\]
The negative zero of \(H\) is at most \(-q\), so \(A_2(x;q)\preceq H\).
In all cases, \(H(0)=0\) and interlacing with \(A_{n-1}(x;q)\)
puts all other zeros of \(H\) below zero. Since \(H\) has leading
coefficient \(q>0\), both polynomials have nonnegative coefficients.
By~\cite[Example~3.7(b)]{AthanasiadisWagner2026}, we obtain the
fully interlacing column \((A_{n-1}(x;q),H)^{\mathsf T}\).
Reflection across the reverse diagonal preserves full interlacing by
\cite[Theorem~4.1]{AthanasiadisWagner2026}, giving the row
\((H,A_{n-1}(x;q))\).
The constant matrix
\[
 \begin{pmatrix}
 1&\cdots&1\\
 b_1+q&\cdots&b_s+q
 \end{pmatrix}
\]
has nonnegative entries. Its minor on columns \(i<j\) is
\(b_j-b_i\ge0\). Thus the constant matrix is totally nonnegative.
It follows from~\cite[Example~3.7(a)]{AthanasiadisWagner2026} that the
matrix is fully interlacing.
The product of \((H,A_{n-1}(x;q))\) with this matrix is fully
interlacing by~\cite[Theorem~4.3(b)]{AthanasiadisWagner2026}.
The \(i\)th entry of the product is \(H+(b_i+q)A_{n-1}(x;q)\), which is the
\(i\)th polynomial in the theorem.

Conversely, if \(c>q+1\), the choice \(n=3\) and \(b_1=c-q-1>0\)
gives a polynomial with nonreal zeros by the proof of
Theorem~\ref{thm:cycle-free-middle}. Even the row of length one is
then not fully interlacing. This proves necessity.
Finally, for every \(n\ge3\) the constant term of the \(i\)th entry
is \(q^{n-1}(q+b_i)\), so full interlacing requires \(b_i\ge-q\).
\end{proof}

\begin{example}\label{ex:cycle-full-interlacing}
For \(n\ge3\), forbid the two-cycles \((1\,2)\) and \((1\,3)\) in
permutations of \([n]\). The two cycles share the letter \(1\) and
cannot occur together. For either cycle, the permutations
containing it contribute
\(qxA_{n-2}(x;q)\). Deleting the cycle removes one cycle and one
excedance, and increasing relabeling preserves all other excedances.
Thus forbidding the two cycles subtracts \(2qxA_{n-2}(x;q)\)
from \(A_n(x;q)\).

Now give the singleton cycle \((1)\) weight \(t\ge0\).
Every permutation containing \((1)\) already avoids the two forbidden
cycles. Deleting \((1)\) leaves an arbitrary permutation of the other
\(n-1\) letters, so changing the weight of the singleton \((1)\)
from \(q\) to \(t\) adds \((t-q)A_{n-1}(x;q)\).
The resulting excedance polynomial is
\[
 A_n(x;q)+(t-q)A_{n-1}(x;q)-2qxA_{n-2}(x;q).
\]
By Theorem~\ref{thm:cycle-full-interlacing}, for any
\(0\le t_1\le\cdots\le t_s\), the corresponding row of excedance
polynomials is fully interlacing for \(0<q\le1\),
since \(2q\le q+1\). At \(n=3\) and \(t=q\), the polynomial is
\(q\bigl(x^2+(1+q)x+q^2\bigr)\). The monic quadratic in parentheses
has discriminant \((1-q)(1+3q)\), which is negative for \(q>1\).
Thus \(0<q\le1\) is the exact range
for full interlacing uniformly over \(n\ge3\) and nonnegative
weights \(t_1\le\cdots\le t_s\). At \(t=q\) and \(0<q\le1\),
all zeros are negative. They are simple except at \(n=3\), \(q=1\),
where the polynomial is \((1+x)^2\).
\end{example}

We now fix \(n,q,c\) and vary the middle coefficient \(b\).
The matrix representation
uses the zeros of \(A_{n-1}(x;q)\) and the residues of
\((A_n(x;q)-cxA_{n-2}(x;q))/A_{n-1}(x;q)\).

\begin{proposition}\label{prop:cycle-arrowhead}
Fix \(q>0\), \(n\ge4\), and \(c\le q+1\). Let
\(r_1<\cdots<r_{n-2}<0\) be the zeros of \(A_{n-1}(x;q)\), and define the positive numbers
\[
 v_i=\sqrt{-\frac{A_n(r_i;q)-cr_iA_{n-2}(r_i;q)}
                 {A_{n-1}'(r_i;q)}}\qquad(1\le i\le n-2).
\]
For \(b\in\mathbb R\), let
\[
 M_{n,q,c}(b)=
 \begin{pmatrix}
 c-b-q-\displaystyle\sum_{i=1}^{n-2}(1-r_i)&v_1&\cdots&v_{n-2}\\
 v_1&r_1&&0\\
 \vdots&&\ddots&\\
 v_{n-2}&0&&r_{n-2}
 \end{pmatrix}.
\]
Then the real symmetric arrowhead matrix \(M_{n,q,c}(b)\) satisfies
\[
 A_n(x;q)+bA_{n-1}(x;q)-cxA_{n-2}(x;q)
 =q\det\bigl(xI-M_{n,q,c}(b)\bigr).
\]
If \(\xi_1(b)<\cdots<\xi_{n-1}(b)\) are the zeros of the polynomial
on the left, then
\[
 \xi_1(b)<r_1<\xi_2(b)<\cdots<r_{n-2}<\xi_{n-1}(b),
\]
and
\[
 \frac{d\xi_j(b)}{db}
 =-\left(1+\sum_{i=1}^{n-2}
       \frac{v_i^2}{(\xi_j(b)-r_i)^2}\right)^{-1}\in(-1,0).
\]
\end{proposition}

\begin{proof}
The proof of Theorem~\ref{thm:cycle-free-middle} gives
\(A_{n-1}(x;q)\prec A_n(x;q)-cxA_{n-2}(x;q)\).
Put
\[
 R(x)=\frac{A_n(x;q)-cxA_{n-2}(x;q)}{A_{n-1}(x;q)}.
\]
At \(r_i\), strict interlacing and positive leading coefficients give
opposite signs for \(A_n(r_i;q)-cr_iA_{n-2}(r_i;q)\) and
\(A_{n-1}'(r_i;q)\). The residue of \(R\) at \(r_i\) therefore satisfies
\[
 \frac{A_n(r_i;q)-cr_iA_{n-2}(r_i;q)}{A_{n-1}'(r_i;q)}<0.
\]
By the definition of \(v_i\), this residue is \(-v_i^2\), with \(v_i>0\).

The recurrence~\eqref{eq:cycle-insertion} gives
\[
 \frac{A_n(x;q)}{A_{n-1}(x;q)}
 =q+(n-1)x+x(1-x)
       \sum_{i=1}^{n-2}\frac1{x-r_i}.
\]
Since
\[
 \frac{x(1-x)}{x-r_i}=-x+(1-r_i)+\frac{r_i(1-r_i)}{x-r_i},
\]
the polynomial part of \(A_n(x;q)/A_{n-1}(x;q)\) is
\(x+q+\sum_i(1-r_i)\). The polynomial part of
\(xA_{n-2}(x;q)/A_{n-1}(x;q)\) is one, since both Eulerian
polynomials have leading coefficient \(q\). Adding the simple fractions
with residues \(-v_i^2\) to the polynomial part of \(R\) gives
\[
 R(x)=x+q+\sum_{i=1}^{n-2}(1-r_i)-c
       -\sum_{i=1}^{n-2}\frac{v_i^2}{x-r_i}.
\]
For \(x\ne r_i\) for all \(i\), the Schur complement of the
lower-right diagonal block of \(xI-M_{n,q,c}(b)\) gives
\[
\begin{aligned}
 \det\bigl(xI-M_{n,q,c}(b)\bigr)
 &=\prod_{i=1}^{n-2}(x-r_i)\\
 &\quad\times\left(x+b+q+\sum_{i=1}^{n-2}(1-r_i)-c
       -\sum_{i=1}^{n-2}\frac{v_i^2}{x-r_i}\right).
\end{aligned}
\]
The factor in parentheses is \(R(x)+b\). Since
\(A_{n-1}(x;q)=q\prod_i(x-r_i)\), multiplying the determinant
by \(q\) gives the three-term combination for \(x\ne r_i\).
Both sides are polynomials, so the determinant identity holds for all \(x\).

Differentiating the partial fraction expansion gives
\[
 R'(x)=1+\sum_{i=1}^{n-2}\frac{v_i^2}{(x-r_i)^2}>1.
\]
At each pole \(r_i\), the negative residue \(-v_i^2\) gives the
limits \(+\infty\) from the left and \(-\infty\) from the right.
As \(x\to\pm\infty\), \(R(x)=x+O(1)\). Thus on each interval cut out by
\(r_1,\ldots,r_{n-2}\), \(R\) increases from \(-\infty\) to \(+\infty\).
The equation \(R(x)=-b\) has one simple solution \(\xi_j(b)\) on each
interval. These \(n-1\) solutions exhaust the zeros of
\(A_{n-1}(x;q)(R(x)+b)\), which has degree \(n-1\).
Since \(R'(\xi_j(b))>0\), implicit differentiation of
\(R(\xi_j(b))+b=0\) gives \(\xi_j'(b)=-1/R'(\xi_j(b))\).
\end{proof}

\begin{remark*}
If \(b_1<b_2\), then \(\xi_j(b_2)<\xi_j(b_1)\) for every \(j\).
Only the upper-left entry depends on \(b\). The rank-one perturbation is
\(M_{n,q,c}(b_2)=M_{n,q,c}(b_1)-(b_2-b_1)e_1e_1^{T}\), where \(e_1\)
is the first coordinate vector in \(\mathbb R^{n-1}\).

For \(b\in\mathbb R\), Proposition~\ref{prop:cycle-arrowhead} gives
\(\xi_j(b)<r_j\le r_{n-2}<0\) for \(1\le j\le n-2\).
The constant term is \(q^{n-1}(q+b)\), so \(\xi_{n-1}(-q)=0\).
Since \(\xi_{n-1}(b)\) decreases strictly with \(b\), it is positive
for \(b<-q\) and negative for \(b>-q\).
\end{remark*}

\begin{example}
Take \(n=4\) and \(q=1\) in Example~\ref{ex:cycle-full-interlacing}.
The polynomial for singleton weight \(t\) is
\(x^3+(8+t)x^2+(5+4t)x+t\), so \(b=t-1\) and \(c=2\) in
Proposition~\ref{prop:cycle-arrowhead}.
At \(t=0\) it is \(x(x^2+8x+5)\), and at \(t=1\) it is
\((1+x)(x^2+8x+1)\). By the proposition, each interval cut out by
the zeros of \(A_3(x;1)\) contains one zero of each polynomial.
Since \(b\) increases with \(t\), the zero at \(t=1\) lies to the
left of the zero at \(t=0\) in each interval. Thus the zeros at
\(t=1\) strictly alternate left of those at \(t=0\).
The largest zero decreases from \(0\) to \(-4+\sqrt{15}\).
\end{example}

\subsection{The Eulerian transformation with a cycle weight}

For \(q>0\), define the linear map \(\mathcal A^{(q)}:\mathbb R[t]\to
\mathbb R[x]\) by $\mathcal A^{(q)}(t^j)=x^jA_j(1/x;q)$ for $j\ge0$.
In particular, \(\mathcal A^{(q)}(1)=1\). For \(j\ge1\),
\(x^{j-1}A_j(1/x;q)\) is the reversal of \(A_j(x;q)\) in degree
\(j-1\). The zeros of the reversal are the reciprocals of the zeros of
\(A_j(x;q)\).
At \(q=1\), symmetry of the ordinary Eulerian polynomials gives
\(\mathcal A^{(1)}(t^j)=xA_j(x)\) for \(j\ge1\).
Thus \(\mathcal A^{(1)}\) is the Eulerian transformation studied by
Br\"and\'en and Jochemko~\cite{BrandenJochemko2022}.
Let \(n\ge1\) and let \(f\) be a nonzero nonnegative linear
combination of \(t^{n-k}(1+t)^k\) for \(0\le k\le n\).
Athanasiadis~\cite[Theorem~1.1(a)]{Athanasiadis2025} proved that
\(\mathcal A^{(1)}(f)\) has only real zeros and satisfies
\[
 \mathcal A^{(1)}((1+t)^n)\preceq\mathcal A^{(1)}(f)
          \preceq xA_n(x).
\]
All three polynomials have degree \(n\), so in each comparison the
polynomial on the left alternates left of the one on the right.
For every \(q>0\), each monomial \(t^j\) has a real-rooted image.
Even at \(q=1\), a nonnegative sum of monomials can have an image
with nonreal zeros. For example,
\(\mathcal A^{(1)}(1+t^2)=1+x+x^2\).
In degree two, a nonnegative combination of \(t^2,t(1+t),(1+t)^2\)
with constant term one has linear coefficient at least two.
Thus \(1+t^2\) lies outside this class. We ask which cycle weights give
real-rooted images of all such combinations in all degrees \(n\ge1\).

\begin{theorem}\label{thm:cycle-transform}
Fix \(q>0\), and let the linear map \(\mathcal A^{(q)}:\mathbb R[t]\to\mathbb R[x]\)
be given by \(\mathcal A^{(q)}(t^j)=x^jA_j(1/x;q)\) for \(j\ge0\),
with \(A_0(x;q)=1\). Then \(q\ge1\) if and only if, for every
\(n\ge1\) and all \(a_0,\ldots,a_n\ge0\) with positive sum,
the polynomial
\[
 \mathcal A^{(q)}\!\left(\sum_{k=0}^n a_k t^{n-k}(1+t)^k\right)(x)
 =\sum_{k=0}^n a_k\sum_{j=n-k}^n\binom{k}{n-j}x^jA_j(1/x;q)
\]
has only real zeros.
For \(q\ge1\), the image has only simple nonpositive zeros, and
\(x^{n-1}A_n(1/x;q)\) interlaces the image. Strict interlacing fails
exactly when \(n\ge2\) and \(a_1=\cdots=a_n=0\).
\end{theorem}

For \(q\ge1\), fix a negative zero \(r\) of \(x^{n-1}A_n(1/x;q)\).
By Lemma~\ref{lem:cycle-basis-signs}, the images of
\(t^{n-k}(1+t)^k\) for \(1\le k\le n\) have a common nonzero sign
at \(r\). The image of \(t^n\) vanishes at \(r\). Thus the image of
\(\sum_{k=0}^n a_k t^{n-k}(1+t)^k\) has the common sign at \(r\) whenever
all \(a_k\ge0\) and some \(a_k>0\) for \(k\ge1\). For \(0<q<1\),
Lemma~\ref{lem:cycle-transform-obstruction} uses the first three
coefficients after the constant term to show that the image of
\((1+t)^n+2^nt^n\) has a nonreal zero in every sufficiently large degree.

\begin{lemma}\label{lem:cycle-basis-signs}
Let \(q\ge1\), \(n\ge2\) and \(1\le k\le n\). At every negative zero
\(r\) of \(x^nA_n(1/x;q)\), the value
\[
 \mathcal A^{(q)}\!\left(t^{n-k}(1+t)^k\right)(r)
 =\sum_{j=n-k}^n\binom{k}{n-j}r^jA_j(1/r;q)
\]
has the same nonzero sign as
\(\left.\dfrac{d}{dx}\bigl(x^nA_n(1/x;q)\bigr)\right|_{x=r}\).
When \(k=n\), the term with \(j=0\) is one since \(A_0(x;q)=1\).
\end{lemma}

\begin{proof}
For any \(q>0\), write \(P_j(x)=x^jA_j(1/x;q)\), with
\(P_0=1\). For \(x\ne0\), differentiation gives
$x^{j-1}A_j'(1/x;q)=jP_j(x)-xP_j'(x)$.
Replacing \(x\) by \(1/x\) in~\eqref{eq:cycle-insertion} and
multiplying by \(x^{j+1}\) gives
\(P_{j+1}=(qx+j)P_j+(x-1)(jP_j-xP_j')\). Collecting terms yields
the polynomial identity
\begin{equation}\label{eq:cycle-reversed}
 P_{j+1}=(j+q)xP_j+x(1-x)P_j'\qquad(j\ge0).
\end{equation}
For \(j\ge1\), the leading term \(qx^{j-1}\) of \(A_j(x;q)\)
gives the linear term \(qx\) of \(P_j\). Thus \(P_j(0)=0\) and
\(P_j'(0)=q\), so zero is simple. For \(j\ge1\),
\(A_j(x;q)\prec A_{j+1}(x;q)\). Taking reciprocals reverses both
lists of negative zeros and preserves interlacing, so
\(P_j/x\prec P_{j+1}/x\).
The leading coefficient of \(P_j\) is \(q^j\), since \(A_j(0;q)=q^j\).

Assume \(q\ge1\) and fix a negative zero \(r\) of \(P_n\).
We want the sign of the image of \(t^{n-k}(1+t)^k\) at \(x=r\).
Extend \(\mathcal A^{(q)}\) \(\mathbb R[x]\)-linearly to \(\mathbb R[x][t]\)
by setting \(\mathcal A^{(q)}(F)=\sum_\ell c_\ell(x)P_\ell(x)\) for
\(F(t,x)=\sum_\ell c_\ell(x)t^\ell\). Then
\[
 \mathcal A^{(q)}(F)(r)
 =\sum_\ell c_\ell(r)P_\ell(r)=\mathcal A^{(q)}(F(t,r))(r).
\]
Set \(a=(1-r)^{-1}\), so \(0<a<1\). At \(x=r\), \(t+a(1-x)=1+t\).
Thus replacing the \(k\) factors \(1+t\) by \(t+a(1-x)\) preserves
the value of the image at \(r\). With \(a\) held constant, the images
obtained by adjoining factors \(t+a(1-x)\) one at a time satisfy a
differential recurrence in \(x\). In particular,
\[
 \mathcal A^{(q)}\!\left((t+a(1-x))F(t,x)\right)
   =\sum_\ell c_\ell(x)\bigl(P_{\ell+1}(x)+a(1-x)P_\ell(x)\bigr).
\]
For \(n-k\le j\le n\), let \(f_j\) be the image of
\(t^{n-k}(t+a(1-x))^{j-n+k}\). We show by induction that at every
zero \(s\) of \(P_j\),
\[
 \frac{f_j(s)}{P_j'(s)}\ge0,
\]
with strict inequality when \(s<0\) and \(j>n-k\). Taking \(j=n\)
and \(s=r\) gives the required sign.

In the binomial expansion, the term with \(t^\ell\) uses \(a(1-x)\)
from \(j-\ell\) of the \(j-n+k\) factors, giving
\[
 f_j(x)=\sum_{\ell=n-k}^j\binom{j-n+k}{j-\ell}
          a^{j-\ell}(1-x)^{j-\ell}P_\ell(x),
 \qquad n-k\le j\le n.
\]
Thus \(f_{n-k}=P_{n-k}\).

Multiplying \(t^\ell\) by \(t\) replaces \(P_\ell\) by \(P_{\ell+1}\).
For each summand of \(f_j\), the product rule and~\eqref{eq:cycle-reversed} give
\[
 \begin{aligned}
 &\left((j+q)x+x(1-x)\frac{d}{dx}\right)
       \bigl((1-x)^{j-\ell}P_\ell\bigr)\\
 &\quad=(1-x)^{j-\ell}
    \bigl((j+q)xP_\ell-(j-\ell)xP_\ell+x(1-x)P_\ell'\bigr)\\
 &\quad=(1-x)^{j-\ell}P_{\ell+1}.
 \end{aligned}
\]
Multiply by \(\binom{j-n+k}{j-\ell}a^{j-\ell}\) and sum over
\(n-k\le\ell\le j\). Adding the contribution \(a(1-x)f_j\)
from the other term of \(t+a(1-x)\) gives
\begin{equation}\label{eq:cycle-transform-factor}
 f_{j+1}=\bigl((j+q)x+a(1-x)\bigr)f_j+x(1-x)f_j'
 \qquad(n-k\le j<n).
\end{equation}
In particular, \(f_{j+1}(0)=a f_j(0)\). If \(f_j\) has degree \(j\),
the coefficient of \(x^{j+1}\) in \(f_{j+1}\) is
\((j+q-a-j)\lc(f_j)=(q-a)\lc(f_j)\).
The initial polynomial \(f_{n-k}=P_{n-k}\) has degree \(n-k\),
positive leading coefficient, and nonnegative constant term.
Since \(a>0\) and \(q-a>0\), induction gives \(\deg f_j=j\),
\(\lc(f_j)>0\), and \(f_j(0)\ge0\) for every \(n-k\le j\le n\).

At the initial index \(j=n-k\), the ratios \(f_j(s)/P_j'(s)\) at
the zeros \(s\) of \(P_j\) vanish because \(f_{n-k}=P_{n-k}\).
If \(k=n\), the assertion for \(j=0\) is vacuous since \(P_0=1\).

Suppose \(j<n\) and \(f_j(s)/P_j'(s)\ge0\) at every zero \(s\) of \(P_j\).
Both \(f_j\) and \(P_j\) have degree \(j\). The polynomial \(P_j\)
has leading coefficient \(q^j\) and simple zeros. For \(P_j(x)\ne0\),
partial fractions give
\[
 \frac{f_j(x)}{P_j(x)}
 =\frac{\lc(f_j)}{q^j}
  +\sum_{P_j(s)=0}\frac{f_j(s)}{P_j'(s)}\,\frac1{x-s}.
\]
A common zero of \(f_j\) and \(P_j\) contributes zero to the sum.
Let \(y\) be a negative zero of \(P_{j+1}\). Since \(P_1(x)=qx\)
has no negative zero, \(j\ge1\). The strict interlacing of \(P_j/x\)
and \(P_{j+1}/x\) gives \(P_j(y)\ne0\).
To use the partial fractions for \(f_j/P_j\), we
express \(f_{j+1}(y)/((1-y)P_j(y))\) in terms of that quotient
and its derivative. At \(x=y\), equation~\eqref{eq:cycle-reversed} gives
\((j+q)yP_j(y)=-y(1-y)P_j'(y)\), so
\[
 \frac{(j+q)y+a(1-y)}{1-y}
   = a-y\frac{P_j'(y)}{P_j(y)}.
\]
Divide~\eqref{eq:cycle-transform-factor} at \(x=y\) by
\((1-y)P_j(y)\) and substitute this expression. The quotient rule
and the partial fraction expansion give
\[
 \begin{aligned}
 \frac{f_{j+1}(y)}{(1-y)P_j(y)}
 &=a\frac{f_j(y)}{P_j(y)}
     +y\left.\frac{d}{dx}\frac{f_j(x)}{P_j(x)}\right|_{x=y}\\
 &=a\frac{\lc(f_j)}{q^j}
   +\sum_{P_j(s)=0}\frac{f_j(s)}{P_j'(s)}
            \frac{(a-1)y-as}{(y-s)^2}>0.
 \end{aligned}
\]
For each zero \(s\) of \(P_j\), the factor multiplying
\(f_j(s)/P_j'(s)\) in the last sum is
\[
 \frac{a}{y-s}-\frac{y}{(y-s)^2}=\frac{(a-1)y-as}{(y-s)^2}.
\]
The term \(a\,\lc(f_j)/q^j\) is positive, and the induction
hypothesis gives \(f_j(s)/P_j'(s)\ge0\). Since \(y<0\) and \(s\le0\),
we also have \((a-1)y>0\) and \(-as\ge0\). The strict interlacing and positive leading
coefficients of \(P_j/x\) and \(P_{j+1}/x\) give
\[
 \frac{P_j(y)}{P_{j+1}'(y)}
 =\frac{(P_j/x)(y)}{(P_{j+1}/x)'(y)}>0,
\]
where we used \(P_{j+1}'(y)=y(P_{j+1}/x)'(y)\). Since both ratios
and \(1-y\) are positive, we have
\[
 \frac{f_{j+1}(y)}{P_{j+1}'(y)}
 =\frac{f_{j+1}(y)}{(1-y)P_j(y)}
  \frac{P_j(y)}{P_{j+1}'(y)}(1-y)>0.
\]
At the zero \(y=0\), we have
\(f_{j+1}(0)\ge0\) and \(P_{j+1}'(0)=q>0\). This covers every zero of
\(P_{j+1}\) and completes the induction.

Finally, \(k\ge1\) gives \(n>n-k\), so the claim at \(j=n\)
and \(s=r<0\) gives \(f_n(r)/P_n'(r)>0\). Since \(a(1-r)=1\),
the defining sum for \(f_n\) gives
\[
 f_n(r)=\sum_{\ell=n-k}^n\binom{k}{n-\ell}P_\ell(r)
       =\mathcal A^{(q)}\!\left(t^{n-k}(1+t)^k\right)(r).
\]
This is the asserted sign.
\end{proof}

For \(n=2\) and \(k=1\), the image of \(t(1+t)\) is \(qx(2+qx)\).
At the negative zero \(r=-1/q\) of \(x^2A_2(1/x;q)=qx(1+qx)\),
the image has value \(-1\), and the derivative of \(x^2A_2(1/x;q)\)
has value \(-q\). The two values have the same sign.

To rule out \(0<q<1\), we use a necessary condition on the first
three coefficients after the constant term. Suppose the polynomial
\(1+az+bz^2+cz^3+\cdots\) has nonnegative coefficients and only real
zeros. Then
\begin{equation}\label{eq:cycle-transform-three-coefficients}
 a(a^2-4b)+12c\ge0.
\end{equation}
Indeed, the polynomial is positive on \([0,\infty)\), so each zero
is negative. Since its constant term is one, it factors as
\(\prod_i(1+\lambda_i z)\) with \(\lambda_i>0\).
The coefficients \(a,b,c\) are the first three elementary symmetric
functions of the \(\lambda_i\), and Newton's identities give
\[
 \sum_i\lambda_i=a,\qquad
 \sum_i\lambda_i^2=a^2-2b,\qquad
 \sum_i\lambda_i^3=a^3-3ab+3c.
\]
Using these power sums, we obtain
\[
 \begin{aligned}
 \sum_i\lambda_i(a-2\lambda_i)^2
 &=a^3-4a(a^2-2b)+4(a^3-3ab+3c)\\
 &=a^3-4ab+12c\ge0.
 \end{aligned}
\]
For two factors, \(a^2-4b=(\lambda_1-\lambda_2)^2\).
For \((1+z)^3\), all zeros are real and \(a=b=3\), \(c=1\).
For this cubic, the term \(12c\) is needed to make
\(a(a^2-4b)+12c\) nonnegative.

\begin{lemma}\label{lem:cycle-transform-obstruction}
Let \(0<q<1\). For every sufficiently large \(n\), the polynomial
\[
 \mathcal A^{(q)}\!\left((1+t)^n+2^nt^n\right)(x)
 =1+\sum_{j=1}^n\binom{n}{j}x^jA_j(1/x;q)+2^nx^nA_n(1/x;q)
\]
has a nonreal zero.
\end{lemma}

\begin{proof}
Write \(P_j(x)=x^jA_j(1/x;q)\), with \(P_0=1\).
For \(j\ge1\), the linear coefficient of \(P_j\) is \(q\).
Hence \((1+t)^n\) contributes \(q(2^n-1)\)
to the linear coefficient of the image, and \(2^nt^n\) contributes
\(2^nq\). Substituting \(x=z/2^n\) divides their sum by \(2^n\)
and makes the linear coefficient tend to \(2q\). Each zero \(x_0\)
of the original image corresponds to the zero \(2^nx_0\) of the
scaled image, so the substitution preserves real-rootedness.
Write
\[
 \mathcal A^{(q)}\!\left((1+t)^n+2^nt^n\right)(z/2^n)
 =1+\alpha_nz+\beta_nz^2+\gamma_nz^3+\cdots.
\]
Since \(\alpha_n\to2q>0\), the inequality
\(\alpha_n(\alpha_n^2-4\beta_n)+12\gamma_n\ge0\) will fail for all sufficiently
large \(n\) if \(\gamma_n\to0\) and \(\beta_n\) tends to a value greater
than \(q^2\).
The recurrence~\eqref{eq:cycle-reversed}, valid for every \(q>0\), gives
\[
 \frac{P_{j+1}}{(1-x)^{j+q+1}}
 =x\frac{d}{dx}\frac{P_j}{(1-x)^{j+q}}.
\]
The operator \(x\,d/dx\) multiplies the coefficient of \(x^i\) by \(i\).
Iterating from \(P_0/(1-x)^q=(1-x)^{-q}\) gives, for every \(q>0\),
\begin{equation}\label{eq:cycle-transform-coefficient}
 [x^i]\frac{P_j(x)}{(1-x)^{j+q}}
 =\binom{q+i-1}{i}i^j\qquad(i\ge1,\ j\ge0).
\end{equation}
For \(q>0\) and \(j\ge1\), the quotient \(P_j(x)/(1-x)^{j+q}\)
has zero constant term. By~\eqref{eq:cycle-transform-coefficient},
its coefficients of \(x\) and \(x^2\) are \(q\) and
\(\tfrac12 q(q+1)2^j\), respectively. Multiplying by
\((1-x)^{j+q}=1-(j+q)x+O(x^2)\) gives
\begin{equation}\label{eq:cycle-first-terms}
 [x]P_j=q,\qquad [x^2]P_j=\frac{q(q+1)}2\,2^j-q(j+q).
\end{equation}

For \(0<q<1\) and \(i\ge1\), we have
\(\binom{q+i-1}{i}=\prod_{k=1}^i(q+k-1)/k\le1\).
Both \(P_j\) and \((1-x)^{-j-q}\) have nonnegative coefficients,
and the latter has constant term one. Hence
\eqref{eq:cycle-transform-coefficient} gives \(0\le[x^i]P_j\le i^j\).
The binomial expansion of \((1+t)^n\) gives
\[
 [x^i]\mathcal A^{(q)}((1+t)^n)
 =\sum_{j=0}^n\binom nj[x^i]P_j
 \le\sum_{j=0}^n\binom nj i^j=(1+i)^n
 \qquad(i\ge1).
\]
Adding \(2^nP_n\) to the image of \((1+t)^n\) and substituting
\(x=z/2^n\) gives
\[
 \begin{aligned}
 \alpha_n&=2^{-n}\bigl(q(2^n-1)+2^nq\bigr)=q(2-2^{-n}),\\
 0\le \beta_n-2^{-n}[x^2]P_n&\le(3/4)^n,\\
 0\le \gamma_n&=2^{-3n}[x^3]\mathcal A^{(q)}((1+t)^n)
           +2^{-2n}[x^3]P_n\\
 &\le\frac{4^n}{2^{3n}}+\frac{3^n}{2^{2n}}
  =(1/2)^n+(3/4)^n.
 \end{aligned}
\]
By~\eqref{eq:cycle-first-terms},
\(2^{-n}[x^2]P_n=q(q+1)/2-q(n+q)2^{-n}\). Hence
\[
 \alpha_n\longrightarrow 2q,\qquad
 \beta_n\longrightarrow\frac{q(q+1)}2,\qquad \gamma_n\longrightarrow0.
\]

If the scaled image has only real zeros, its nonnegative coefficients
and constant term one give
\(\alpha_n(\alpha_n^2-4\beta_n)+12\gamma_n\ge0\)
by~\eqref{eq:cycle-transform-three-coefficients}. For each fixed
\(0<q<1\), we have $\alpha_n^2-4\beta_n\to4q^2-2q(q+1)=2q(q-1)$. Hence
\(\alpha_n(\alpha_n^2-4\beta_n)+12\gamma_n\to4q^2(q-1)<0\). For every
sufficiently large \(n\), the scaled image violates the inequality,
so the original image also has a nonreal zero.
\end{proof}

For \(q=1/2\), the image of \((1+t)^4+16t^4\) already has a
nonreal zero. Indeed, substituting \(x=2z\) in the image gives
\[
 \mathcal A^{(1/2)}\!\left((1+t)^4+16t^4\right)(2z)
   =1+31z+281z^2+310z^3+17z^4.
\]
Here \(31(31^2-4\cdot281)+12\cdot310=-1333<0\), in violation
of~\eqref{eq:cycle-transform-three-coefficients}.

\begin{proof}[Proof of Theorem~\ref{thm:cycle-transform}]
Write \(P_j(x)=x^jA_j(1/x;q)\), with \(P_0=1\).
Suppose first that \(q\ge1\). In degree one, the image of
\(a_0t+a_1(1+t)\) is \(a_1+q(a_0+a_1)x\). Since \(a_0+a_1>0\),
it has one simple nonpositive zero. By convention, the nonzero
constant \(P_1/x=q\) strictly interlaces this linear polynomial,
including when \(a_1=0\).

Let \(n\ge2\), and write
\[
 f(t)=\sum_{k=0}^n a_k t^{n-k}(1+t)^k,
 \qquad a_k\ge0,\quad\sum_{k=0}^n a_k>0.
\]
Here \([t^{n-1}]f=\sum_{k=1}^n k a_k\), so \(f\) is a multiple
of \(t^n\) exactly when \(a_1=\cdots=a_n=0\).
Each basis polynomial \(t^{n-k}(1+t)^k\) is monic of degree \(n\),
so \(\lc(f)=\sum_k a_k>0\). Since \(P_j\) has leading term
\(q^jx^j\), \(\mathcal A^{(q)}(f)\) has leading term
\(q^n\lc(f)x^n\).

Suppose \(a_k>0\) for some \(k\ge1\), and let \(r\) be any zero of
\(P_n/x\). Then \(r<0\) and \(\mathcal A^{(q)}(t^n)(r)=P_n(r)=0\).
By linearity,
\[
 \frac{\mathcal A^{(q)}(f)(r)}{P_n'(r)}
 =\sum_{k=1}^n a_k\,
   \frac{\mathcal A^{(q)}\!\left(t^{n-k}(1+t)^k\right)(r)}{P_n'(r)}>0.
\]
Each quotient in the sum is positive by Lemma~\ref{lem:cycle-basis-signs},
and at least one coefficient \(a_k\) in the sum is positive.
Write the zeros of \(P_n/x\) as \(r_1<\cdots<r_{n-1}<0\).
At \(r_i\), the derivative of \(P_n/x\) has sign
\((-1)^{n-1-i}\). Since \(P_n'(r_i)=r_i(P_n/x)'(r_i)\) and
\(r_i<0\), the sign of $\mathcal A^{(q)}(f)(r_i)$ is
$(-1)^{n-i}$ for $1\le i\le n-1$. Thus the image has a zero in
every gap between consecutive \(r_i\). The signs \((-1)^n\) near
\(-\infty\) and \((-1)^{n-1}\) at \(r_1\) give a zero to the left of
\(r_1\). At \(r_{n-1}\), the image is negative. Since \(P_0=1\) and
\(P_j(0)=0\) for \(j\ge1\), the value of the image at \(0\) equals the
input's constant term \(f(0)=a_n\).
If \(a_n>0\), there is a zero in \((r_{n-1},0)\).
If \(a_n=0\), the point \(0\) itself is a zero.
The image has degree \(n\), so the \(n\) distinct nonpositive zeros
are all simple. Their positions give
\(P_n/x\prec\mathcal A^{(q)}(f)\).

The only remaining case is \(f=a_0t^n\), with \(a_0>0\).
Then \(\mathcal A^{(q)}(f)=a_0P_n\) has simple nonpositive zeros
and \(P_n/x\preceq\mathcal A^{(q)}(f)\). Here \(n\ge2\), so the
polynomials have all \(n-1\) zeros of \(P_n/x\) in common.
This proves sufficiency.

For \(0<q<1\), take \(a_0=2^n\), \(a_n=1\), and \(a_k=0\)
for \(0<k<n\). Then \(f(t)=(1+t)^n+2^nt^n\). By
Lemma~\ref{lem:cycle-transform-obstruction},
its image has a nonreal zero for every sufficiently large \(n\).
This proves necessity.
\end{proof}

We can also vary the constant term in each factor \(1+t\).

\begin{proposition}\label{prop:shifted-cycle-basis}
Fix \(q>0\) and \(\lambda\ge0\).
\begin{enumerate}
\item [$(i)$] Suppose \(q\ge1\) and \(0\le\lambda\le1\).
Let \(n\ge0\) and let \(a_0,\ldots,a_n\ge0\) have positive sum. Then
\[
 \mathcal A^{(q)}\!\left(
 \sum_{k=0}^n a_k t^{n-k}(t+\lambda)^k\right)
\]
has only real zeros.
\item [$(ii)$] If \(\lambda>1\), then
\(\mathcal A^{(q)}\!\left(\lambda^{-n}(t+\lambda)^n+2^nt^n\right)\)
has a nonreal zero for every sufficiently large \(n\).
\end{enumerate}
Consequently, for each fixed \(q\ge1\), every nonzero
nonnegative combination of \(t^{n-k}(t+\lambda)^k\),
\(0\le k\le n\), has a real-rooted image in every degree
\(n\ge0\) if and only if \(0\le\lambda\le1\).
\end{proposition}

\begin{proof}
Suppose first that \(q\ge1\) and \(0\le\lambda\le1\).
Since \(t+\lambda=(1-\lambda)t+\lambda(1+t)\), we have
\[
 t^{n-k}(t+\lambda)^k
 =\sum_{j=0}^k\binom{k}{j}(1-\lambda)^{k-j}\lambda^j
       t^{n-j}(1+t)^j.
\]
All coefficients in this expansion are nonnegative, so
Theorem~\ref{thm:cycle-transform} gives the conclusion for \(n\ge1\).
For \(n=0\), the image is the positive constant \(a_0\).

For \(\lambda>1\) and \(q>0\), put
\(f_n(t)=\lambda^{-n}(t+\lambda)^n+2^nt^n\).
Write \(P_j(x)=x^jA_j(1/x;q)\), with \(P_0=1\).
Since \(P_j(x)\) and \((1-x)^{-j-q}\) have nonnegative coefficients
and the latter series has constant term one,
\eqref{eq:cycle-transform-coefficient} gives
\[
 0\le[x^i]P_j(x)\le\binom{q+i-1}{i}i^j
 \qquad(i\ge1,\ j\ge0).
\]
The image of \(f_n\) is
\[
 \mathcal A^{(q)}(f_n)(x)
 =\sum_{j=0}^n\binom nj\lambda^{-j}P_j(x)+2^nP_n(x).
\]
After substituting \(x=z/2^n\), write the scaled image as
\(1+\alpha_nz+\beta_nz^2+\gamma_nz^3+\cdots\).
For \(i=1,2,3\), the contribution of the first sum to the
coefficient of \(z^i\) is at most
\[
 \binom{q+i-1}{i}
 \left(\frac{1+i/\lambda}{2^i}\right)^n,
\]
Since \(1+i/\lambda<1+i\le2^i\), this contribution tends to zero.
The contribution of \(2^nP_n(x)\) after the substitution is
\[
 [z^i]\bigl(2^nP_n(z/2^n)\bigr)=2^{-(i-1)n}[x^i]P_n(x).
\]
For \(i=3\), the bound on \([x^3]P_n\) gives a contribution of
at most \(\binom{q+2}{3}(3/4)^n\). For \(i=1,2\), use
\eqref{eq:cycle-first-terms} with \(j=n\). Thus
\[
 \alpha_n\longrightarrow q,\qquad
 \beta_n\longrightarrow\frac{q(q+1)}2,\qquad \gamma_n\longrightarrow0.
\]

The polynomial in \(z\) has nonnegative coefficients and constant term one.
If it had only real zeros, then~\eqref{eq:cycle-transform-three-coefficients}
would give
\(\alpha_n(\alpha_n^2-4\beta_n)+12\gamma_n\ge0\).
The left-hand side tends to
\(q\left(q^2-2q(q+1)\right)=-q^2(q+2)<0\), so the inequality fails
for all sufficiently large \(n\). This proves part~$(ii)$.
\end{proof}

\begin{remark}
For $0<q<1$, the real-rootedness conclusion of part~$(i)$
of Proposition~\ref{prop:shifted-cycle-basis} also holds
for $0\le\lambda\le q$.
Write the input as
\[
 \sum_{k=0}^n a_k t^{n-k}(t+\lambda)^k,\qquad
 a_k\ge0,\qquad \sum_{k=0}^n a_k>0.
\]
At \(\lambda=0\), all the polynomials \(t^{n-k}(t+\lambda)^k\)
equal \(t^n\), whose image has only real zeros.
If \(a_k=0\) for every \(k\ge1\), the input is a positive multiple
of \(t^n\) and its image is the same multiple of \(x^nA_n(1/x;q)\).
This includes \(n=0\), when the image is the positive constant \(a_0\).
For \(n=1\),
\[
 \mathcal A^{(q)}\!\left(a_0t+a_1(t+\lambda)\right)
 =\lambda a_1+q(a_0+a_1)x.
\]
Its zero is \(-\lambda a_1/[q(a_0+a_1)]\le0\).

Now suppose \(n\ge2\), \(0<\lambda\le q\), and \(a_k>0\) for some
\(k\ge1\). At each negative zero \(r\) of \(x^nA_n(1/x;q)\), take
\(a=\lambda/(1-r)\). The induction in the proof of
Lemma~\ref{lem:cycle-basis-signs} requires \(0<a<1\) and \(a<q\).
Here \(0<a<q<1\), so \(q-a>0\) and \((a-1)y-as>0\) for
\(y<0\) and \(s\le0\). At \(x=r\) we have \(a(1-r)=\lambda\),
so for \(1\le k\le n\) the image of
\(t^{n-k}(t+\lambda)^k\) has the same sign at \(r\) as the derivative
of \(x^nA_n(1/x;q)\). The \(k=0\) image vanishes at \(r\), so the image
of the sum has that same nonzero sign. The image has degree \(n\)
and leading coefficient \(q^n\sum_{k=0}^n a_k>0\). Its constant
term is \(a_n\lambda^n\ge0\),
since \(\mathcal A^{(q)}(1)=1\) and
\(\mathcal A^{(q)}(t^j)(0)=0\) for \(j\ge1\).
Write \(r_1<\cdots<r_{n-1}<0\) for the negative zeros of
\(x^nA_n(1/x;q)\). The image has alternating signs at the \(r_i\),
with sign \((-1)^{n-1}\) at \(r_1\) and negative sign at \(r_{n-1}\).
Hence the image has a zero in each of the \(n-2\) gaps between
consecutive \(r_i\). The sign \((-1)^n\) near \(-\infty\) gives a zero
to the left of \(r_1\). A positive constant term gives a zero in
\((r_{n-1},0)\). If the constant term vanishes, \(0\) is a zero.
These \(n\) distinct
nonpositive zeros exhaust the degree of the image and are simple.
\end{remark}

For the original basis \(t^{n-k}(1+t)^k\), taking \(q=1\)
in Theorem~\ref{thm:cycle-transform} recovers
Athanasiadis's real-rootedness result.
For \(q=1\) we have \(x^{n-1}A_n(1/x)=A_n(x)\), so the ordinary Eulerian
polynomial \(A_n(x)\) strictly interlaces the transformed sum whenever some
\(t^{n-k}(1+t)^k\) with \(k\ge1\) has positive coefficient.

Reversing the image of \(t^{n-k}(1+t)^k\) in degree \(n\) gives
\(\sum_{i=0}^k\binom{k}{i}x^i A_{n-i}(x;q)\).
The reversed polynomials form the array in~\cite[Remark~7.2(c)]{Athanasiadis2025}
with \(h_j(x)=A_j(x;q)\). In their study of arrays associated with
Eulerian and derangement transformations, Liu and Yan~\cite{LiuYan2025}
give sufficient conditions on \((h_j)\) for each row to form an
interlacing sequence.

Fix \(n-k\) letters of \([n]\), and let \(S\) range over the subsets
containing them. Each of the other \(k\) letters may be included or
omitted, so \(t^{n-k}(1+t)^k\) counts these subsets by size.
For each \(S\), give a permutation \(\pi:S\to S\) weight
\(q\) for each cycle and \(x\) for each \(i\in S\) with \(\pi(i)\le i\).
Increasing relabeling of \(S\) preserves the number of cycles
and each comparison \(\pi(i)\le i\). If \(|S|=j\), then exactly
\(j-\excc(\pi)\) positions satisfy \(\pi(i)\le i\). Thus summing
the weights over permutations of \(S\) gives \(x^jA_j(1/x;q)\).
When \(k=n\), the empty subset contributes one. Summing over the subsets
containing the prescribed letters gives
\(\mathcal A^{(q)}(t^{n-k}(1+t)^k)\).
For \(q\ge1\) and \(1\le k\le n\), Theorem~\ref{thm:cycle-transform}
shows that \(x^{n-1}A_n(1/x;q)\) strictly interlaces this image.

For example, take \(n=k=3\) and \(q=1\). The empty subset
contributes one, and the three singletons contribute \(3x\).
On each two-element subset the transposition contributes \(x\) and
the identity contributes \(x^2\). The six permutations of \([3]\)
contribute \(x+4x^2+x^3\). Thus
\[
 \begin{aligned}
 \mathcal A^{(1)}((1+t)^3)
 &=1+3x+3(x+x^2)+(x+4x^2+x^3)\\
 &=1+7x+7x^2+x^3=(1+x)(1+6x+x^2).
 \end{aligned}
\]
For the same input \((1+t)^3\), real-rootedness can fail when \(q<1\).
Take \(q=1/8\) and set \(x=8z\). Recurrence~\eqref{eq:cycle-reversed}
gives the images of \(t,t^2,t^3\), evaluated at \(x=8z\), as
\(z\), \(z+z^2\), and \(z+11z^2+z^3\), respectively.
The binomial expansion therefore gives
\[
 \begin{aligned}
 \mathcal A^{(1/8)}((1+t)^3)(8z)
 &=1+3z+3(z+z^2)+(z+11z^2+z^3)\\
 &=1+7z+14z^2+z^3.
 \end{aligned}
\]
With \(a=7\), \(b=14\), and \(c=1\), the left side
of~\eqref{eq:cycle-transform-three-coefficients} is
\(7(7^2-4\cdot14)+12=-37\).
Thus at \(q=1/8\), a single basis polynomial, \((1+t)^3\), already
has an image with a nonreal zero.

\section*{Acknowledgements}
This work is supported by the National Natural Science Foundation of China
(No.~12071063) and the Natural Science Foundation of Shandong Province of
China (ZR2026MS0047).

%\bibliographystyle{amsplain}
%\bibliography{references}

\end{document}